\documentclass[10pt]{article}

\usepackage{vmargin}   
\setmarginsrb{3.2cm}{2.5cm}{3.2cm}{2.5cm}{1cm}{1cm}{1cm}{1.6cm}   
\usepackage[fleqn]{amsmath}   
\usepackage{amssymb,amsthm}  
\usepackage{wasysym}

\def\onehalf{\frac{1}{2}}
\newcounter{Prob}
\newcounter{Prop}
\newcounter{Theo}
\newcounter{Def}
\newcounter{Lemm} 
   
{\theoremstyle{remark}
\newtheorem{rem}{Remark}[section]}

\def\Lemma#1#2{\refstepcounter{Lemm}
      \bigskip\noindent { \sf $\,$ Lemma \theLemm:} { \rm #1 }\\
          {\em #2}\\
}
\def\Proposition#1#2{\refstepcounter{Prop}
     \bigskip\noindent { \sf $\,$ Proposition \theProp:} { \rm #1 }\\
          {\em #2}\\
}
\def\Theorem#1#2{\refstepcounter{Theo}
     \bigskip\noindent { \sf $\,$ Theorem \theTheo:} { \rm #1 }\\
          {\em #2}\\
}

\def\Definition#1#2{\refstepcounter{Def}
   \bigskip\noindent   { \sf $\,$ Definition \theDef:} { \rm #1 }\\
          {\em #2}\\
}

\def\Problem#1#2{\refstepcounter{Prob}
   \bigskip\noindent   { \sf $\,$ Problem \theProb:} { \rm #1 }\\
          {\em #2}\\
}

\def\Proof#1{
 \smallskip\noindent {\sf \underline{Proof}:} 
  {\rm #1}
  \smallskip
  \hfill $\Box$\\
  }

\def\proof{\medskip\noindent {\sf{\underline{Proof}:\ }}}

\let\qed=\cqfd
 
\def\negvspace{\vspace{-0.5cm}}
  
\def\id{{\rm id}}
\def\End{{\rm End}}

  \def\fB{\mathfrak{b}}  
    
  \def\fH{\mathfrak{h}}  
     
  \def\fN{\mathfrak{n}}  
     
    \def\fU{U}

\title{Quantum Dynamical coBoundary Equation \\ for  finite dimensional simple Lie algebras}

\author{E. Buffenoir\thanks{e-mail: Eric.Buffenoir@lpta.univ-montp2.fr}, Ph. Roche\thanks{e-mail: Philippe.Roche@lpta.univ-montp2.fr}, V. Terras\thanks{e-mail: Veronique.Terras@lpta.univ-montp2.fr} \\
 Laboratoire de Physique Th\'eorique et d'Astroparticules \\ UMR 5207 CNRS\\ 
 Universit\'e Montpellier 2, 34000 Montpellier, France.}  
\date{\today}

\begin{document}

\maketitle  

\begin{abstract} 
For a finite dimensional simple Lie algebra $\mathfrak{g}$, the standard universal solution $R(x)\in U_q(\mathfrak{g})^{\otimes 2}$ of the Quantum Dynamical Yang--Baxter Equation  quantizes the standard trigonometric solution of the Classical Dynamical Yang--Baxter Equation. 
It  can be built from the standard $R$--matrix and from the  solution $F(x)\in U_q(\mathfrak{g})^{\otimes 2}$ of the Quantum Dynamical coCycle Equation as $R(x)=F^{-1}_{21}(x)\,R\,F_{12}(x).$ $F(x)$ can be computed explicitely as an infinite product through the use of an auxiliary linear equation, the ABRR equation.

Inspired by explicit results in the fundamental representation, it has been conjectured  that, in the case where $\mathfrak{g}=sl(n+1)$ $(n\geq 1)$ only, there could exist an element 
$M(x)\in U_q(sl(n+1))$ such that the dynamical gauge transform  $R^J$ 
 of $R(x)$ by $M(x)$,
\begin{equation*}
R^J=M_1(x)^{-1}M_2(xq^{h_1})^{-1}R(x)M_1(xq^{h_2})M_2(x),
\end{equation*} 
does not depend on $x$ and is a universal solution of the Quantum Yang--Baxter Equation.
In the fundamental  representation, $R^J$ corresponds to the standard solution $R$ for $n=1$ and to the  Cremmer--Gervais's one $R^J_{12}=J_{21}^{-1}R_{12}J_{12}$ for $n>1$.   
For consistency, $M(x)$ should therefore satisfy the Quantum Dynamical coBoundary Equation, i.e
\begin{equation*}
F(x)=\Delta(M(x)){ J} M_2(x)^{-1}(M_1(xq^{h_2}))^{-1},
\end{equation*}
in which $J\in U_q(sl(n+1))^{\otimes 2}$ is the universal cocycle associated to the  Cremmer--Gervais's solution.

The aim of this article is to prove this conjecture and to study the properties of the solutions of the Quantum Dynamical coBoundary Equation.
In particular, by introducing new basic algebraic objects which are the building blocks of the Gauss decomposition of $M(x)$, we construct $M(x)$ in $U_q(sl(n+1))$ as an explicit infinite product which converges in every finite dimensional representation. 
We emphasize  the relations between these  basic objects and some Non Standard Loop algebras and exhibit relations with  the dynamical quantum Weyl group.
\end{abstract}

%%%%%%%%%%%%%%%%%%%%%%%%%%%%%%%%%%%%%%%%%%%%%%%%%%%%%%%%%%%%%%%%%%%%%%%%%%%%%%
\section{Introduction}
%%%%%%%%%%%%%%%%%%%%%%%%%%%%%%%%%%%%%%%%%%%%%%%%%%%%%%%%%%%%%%%%%%%%%%%%%%%%%%

%%%%%%%%%%%%%%%%%%%%%%%%%%%%%%%%%%%%%%%%%%%%%%%%%%%%%%%%%%%%%%%%%%%%%%%%%%%%%%
\subsection{History of the problem}

The theory of dynamical quantum groups is nowadays a well established part of mathematics, 
see for example the review by P. Etingof \cite{E}. This theory originated from the notion of Dynamical Yang-Baxter equation  which arose in the work of Gervais-Neveu on  Liouville theory
\cite{GN} and was formalized  first by G. Felder \cite{Fe} who also understood its 
relation with IRF statistical models. The work \cite{ABB} gives some  connections between Dynamical Yang-Baxter Equation and various models in mathematical physics.

\medskip

In 1991, O. Babelon \cite{Bab} found  a universal explicit solution $F(x)$ to
the Quantum Dynamical Yang-Baxter equation in the case where $\mathfrak{g}=sl(2)$.  He obtained this solution $F(x)$ by showing that it is a Quantum Dynamical coBoundary, i.e. $F(x)=\Delta(M(x)) M_2(x)^{-1}(M_1(xq^{h_2}))^{-1}$. The question whether this work can be generalized to any finite dimensional simple Lie algebra is an important problem which has, until now, received only uncomplete solution.

It has been shown in \cite{ABRR} (see also \cite{JKOS}) that, in the case of any finite dimensional simple Lie algebra $\mathfrak{g}$, the  standard solution $F(x)$ of the Dynamical Cocycle Equation of weight zero satisfying an additional condition (of triangularity type) can be obtained  as the  unique solution of a linear equation now called the ABRR equation.
It  implies that $F(x)$ can be expressed as an explicit  infinite product which converges in any finite dimensional representation.
It remained nevertheless to study whether $F(x)$ is a dynamical coboundary and, if it is the case, to construct explicitely $M(x)\in U_q(\mathfrak{g})$. 

For $\mathfrak{g}=sl(2)$, we have shown \cite{BR2} that $M(x)$ can be written as a simple infinite product which simplifies greatly the solution given in \cite{Bab}.

Concerning the more general case $\mathfrak{g}=sl(n+1)$, the first hint 
appears in the article \cite{CG} on the study of Toda field theory:  
the authors of \cite{CG} proved in particular that, in the fundamental 
representation of $sl(n+1)$, the standard solution of the dynamical 
Yang-Baxter equation (computed first in \cite{BG}) can be  dynamically gauged 
through a matrix $M(x)$ to a constant solution of the Yang-Baxter Equation 
which is non standard in the case where $n\geq 2$ and is now called  the 
``Cremmer-Gervais's'' solution. The expression of $M(x)^{-1}$ is especially 
simple: it is a Vandermonde matrix. 

In the classification of Belavin and Drinfeld  \cite{BD}, in which all the non-skewsymmetric classical r-matrices for simple Lie algebras are, up to an isomorphism, classified by a combinatorial object called the Belavin-Drinfeld triple, the  Cremmer-Gervais's solution  is associated to a particular Belavin-Triple of $sl(n+1)$ known as the shift $\tau$.   

Moreover, it was proved in  \cite{BDF} that it is only for $\mathfrak{g}=sl(n+1)$ that  the standard solution of the dynamical Yang-Baxter equation can be  dynamically gauged to a constant solution of Yang-Baxter Equation which, in that case, is a quantization of the r-matrix associated to the shift $\tau.$

Later, O. Schiffman  \cite{Sc} generalized the notion of Belavin-Drinfeld triple and provided a classification of classical dynamical r-matrices up to isomorphism through the use of generalized  Belavin-Drinfeld triple.
Then, P. Etingof, T. Schedler and O. Schiffmann \cite{ESS} managed to quantize explicitely all the previous generalized Belavin-Drinfeld triple. In particular, they obtained a universal expression for the twist $J\in U_q(sl(n+1))^{\otimes 2}$ associated to the shift.
The universal  Cremmer-Gervais's solution is therefore $R^{J}=J_{21}^{-1}R_{12}J_{12}.$ 

\medskip

On the basis of all these works a natural problem to address is to construct a universal coboundary element in $U_q(\mathfrak{g})$ for $\mathfrak{g}=sl(n+1)$
%\footnote{We recall that, for any other finite dimensional simple Lie algebra $\mathfrak{g}$, such a coboundary 
%does not exist (see \cite{BDF} and Theorem~\ref{lemme2} of the present paper).}
, i.e. to solve the equation 
\begin{equation}
F(x)=\Delta(M(x)){ J} M_2(x)^{-1}(M_1(xq^{h_2}))^{-1}
\end{equation}
where $F(x)$ is the standard solution of the dynamical cocycle equation in 
$U_q(sl(n+1))^{\otimes 2}$ and $J$ is the universal twist associated to the shift.

We will solve this problem in the present work.

%%%%%%%%%%%%%%%%%%%%%%%%%%%%%%%%%%%%%%%%%%%%%%%%%%%%%%%%%%%%%%%%%%%%%%%%%%%%%%%
\subsection{Description and content of the present work}

Section~\ref{sectionQUEA} contains results (old and new) on quantum groups. We recall the definition and 
properties of $U_q(\mathfrak g)$ and give emphasis on the quantum Weyl group and the 
explicit expressions for the standard $R$-matrix.
We then recall the definition of a Belavin-Drinfeld triple and analyze the  
$sl(n+1)$ case  with the special triple known as shift.
We give in this case a direct  construction of the universal twist as a simple finite 
product of $q$-exponentials. Note that this construction is purely combinatorial and does not rely on the result of \cite{ESS}.

In Section~\ref{sectionDQG} we recall some results on dynamical quantum groups. We present the quantum dynamical cocycle equation, the linear equation 
and give a summary of the results of Etingof, Schedler, Schiffmann \cite{ESS}.
We then formulate precisely the Dynamical coBoundary Equation and recall the known results in the fundamental representation of
 $sl(n+1).$

Section~\ref{sectionDBP} is the core of our work. We introduce the notion of primitive loop which is constructed 
from any solution of the Dynamical coBoundary Equation. This primitive loop satisfies a reflection equation with the corresponding  $R$-matrix being of  Cremmer-Gervais's type.
We then study the Gauss decomposition of $M(x)$ and shows that additional properties satisfied by these objects  implies the
 Dynamical coBoundary Equation.
 
Section~\ref{sectionConstr} is devoted to the explicit construction of $M(x)$ in the $U_q(sl(n+1))$ case as an infinite product converging 
in every finite dimensional representation.

In Section~\ref{sectionQWG} we analyze the relations between $M(x)$ and the dynamical quantum Weyl group. 

Section~\ref{sectionQRA} contains a construction of $M(x)$ through the 
representation theory of  non-standard reflection algebras, and in particular
through what we call the primitive representations.

%%%%%%%%%%%%%%%%%%%%%%%%%%%%%%%%%%%%%%%%%%%%%%%%%%%%%%%%%%%%%%%%%%%%%%%%%%%%%%
\section{Results on  Quantum Universal Envelopping Algebras}
\label{sectionQUEA}
%%%%%%%%%%%%%%%%%%%%%%%%%%%%%%%%%%%%%%%%%%%%%%%%%%%%%%%%%%%%%%%%%%%%%%%%%%%%%%

%%%%%%%%%%%%%%%%%%%%%%%%%%%%%%%%%%%%%%%%%%%%%%%%%%%%%%%%%%%%%%%%%%%%%%%%%%%%%%
\subsection{Basic Results on Quantized Simple Lie Algebras}

Let $\mathfrak{g}$ be a finite dimensional simple complex Lie algebra. 
We denote by ${\mathfrak h}$ a Cartan subalgebra and by $r$ the rank of 
$\mathfrak{g}.$ 
Let $\Phi$ be the set of roots,  $(\alpha_{i},i=1,\dots,r)$ a choice of simple 
roots and $\Phi^+$ the corresponding set of positive roots.
 
Let $(,)$ be a non zero $ad$-invariant symmetric bilinear form on 
$\mathfrak{g}$. Its restriction on $\mathfrak{h}$ is non degenerate and 
therefore induces a non degenerate symmetric bilinear form on  
$\mathfrak{h}^*$ that we still denote  $(,)$.
We denote by $\Omega$ the element  in $S^2(\mathfrak{g})^{\mathfrak g}$ 
(the vector space of symmetric elements of $\mathfrak{g}^{\otimes 2}$ invariant 
under the adjoint action)  associated to 
$(,)$, and by $\Omega_{\mathfrak h}$ its projection on 
$\mathfrak{h}^{\otimes 2}.$

Let $A=(a_{ij})$ be the Cartan matrix of $\mathfrak{g}$, with elements
$a_{ij}=2\frac{(\alpha_i,\alpha_j)}{(\alpha_i,\alpha_i)}$. 
There exists a unique collection of coprime  positive integers $d_i$ such 
that $d_i a_{ij}=a_{ji}d_j$.
$(,)$ is then uniquely defined by imposing that $(\alpha_i,\alpha_i)=2 d_i.$

If $\alpha\in\mathfrak{h}^{\star}$, we define the element 
$t_{\alpha}\in \mathfrak{h}$ such that 
$(t_{\alpha},h)=\alpha(h)$, $\forall h\in \mathfrak{h}$, 
and denote $\mathfrak{h}_{\alpha}={\mathbb C}t_{\alpha}.$ 
To each root $\alpha$ we associate the coroot $\alpha^{\vee}=
\frac{2}{(\alpha,\alpha)}\alpha$ and denote 
$h_{\alpha}=t_{\alpha^{\vee}},$ therefore
$ h_{\alpha_i}=\frac{1}{d_i}t_{\alpha_i}.$

Let $\lambda^{\alpha_1},\ldots,\lambda^{\alpha_r} \in \mathfrak{h}^{\star}$ 
be the set of fundamental weights, i.e. $\lambda^{\alpha_i}(h_{\alpha_j})=\delta^i_j.$
We will also denote  by $\zeta^{\alpha_i}=t_{\lambda^{\alpha_i}}\in \mathfrak{h}$.

\bigskip 

Let us now define the Hopf algebra $U_q(\mathfrak{g}).$ 
We will assume that $q$ is a complex number with $0<\vert q\vert <1.$ 
We define for each root $q_{\alpha}=q^{(\alpha,\alpha)/2},$ 
as well as $q_{i}=q_{\alpha_{i}}$,
and we denote $[z]_q=\frac{q^z-q^{-z}}{q-q^{-1}},\ z\in {\mathbb C}.$ 

$U_q(\mathfrak{g})$ is the unital associative algebra generated by 
$e_{\alpha_1},\ldots,e_{\alpha_r},f_{\alpha_1},\ldots,f_{\alpha_r}$ and 
$q^h,h\in \mathfrak{h}$, with defining relations:
\begin{align}
  & q^h q^{h'}=q^{h+h'},\quad 
    q^h e_{\alpha_i} q^{-h}=q^{\alpha_i(h)}e_{\alpha_i},\quad 
    q^h f_{\alpha_i} q^{-h}=q^{-\alpha_i(h)}f_{\alpha_i},\quad
    \forall  h,h'\in \mathfrak{h},
  \\
  & [e_{\alpha_i},f_{\alpha_j}]=\delta_{ij}
    \frac{q_i^{h_{\alpha_i}}-q_i^{-h_{\alpha_i}}}{q_i-q_i^{-1}},
  \\
  & \sum_{k=0}^{1-a_{ij}}\frac{(-1)^k}{[k]_{q_i}!\, [1-a_{ij}-k]_{q_i}!}
    \,e_{\alpha_i}^{1-a_{ij}-k}\,e_{\alpha_j}\,e_{\alpha_i}^k=0,
  \label{serree}\\
  & \sum_{k=0}^{1-a_{ij}}\frac{(-1)^k}{[k]_{q_i}!\, [1-a_{ij}-k]_{q_i}!}
    \,f_{\alpha_i}^{1-a_{ij}-k}\,f_{\alpha_j}\,f_{\alpha_i}^k=0.
  \label{serref}
\end{align}
$U_q(\mathfrak{g})$ is a Hopf algebra with coproduct:
\begin{equation}
  \Delta(q^h)=q^h\otimes q^h,
  \quad
  \Delta(e_{\alpha_i})= e_{\alpha_i}\otimes q_i^{h_{\alpha_i}}
  +1\otimes e_{\alpha_i},
  \quad
  \Delta(f_{\alpha_i})=f_{\alpha_i}\otimes 1
  +q_i^{-h_{\alpha_i}}\otimes f_{\alpha_i}.\label{coproductsln}
\end{equation}

Let us now define different notions associated to the  polarisation of 
$U_q({\mathfrak g}).$
We denote by $U_q(\fB_+)$ (resp. $U_q(\fB_-)$) the
algebra generated by $q^h$, $h\in\fH$, $e_{\alpha_i}$, $i=1,\dots,r$
(resp. $q^h$, $f_{\alpha_i}$, $i=1,\dots,r$),
and by $\fU_q(\fN_+)$ (resp. $\fU_q(\fN_-)$) the subalgebra of 
$U_q(\mathfrak{g})$ generated
by $e_{\alpha_i},\ i=1,\dots,r$ (resp. $f_{\alpha_i},\ i=1,\dots,r$).
We have  $U_q(\fB_+)=\fU_q(\fN_+)\otimes U_q({\mathfrak h})$ as a vector space 
as well as  $U_q(\fB_-)= U_q({\mathfrak h})\otimes \fU_q(\fN_-)$.
We denote by $\iota_{\pm}:U_q(\fB_{\pm})\rightarrow U_q({\mathfrak h})$ 
the projections on the zero-weight subspaces. $\iota_{\pm}$ are morphisms of algebra and we define  the ideals   $U^{\pm}_q({\mathfrak g})=ker \iota_{\pm}$. Endly, we denote by
 ${ C}_q({\mathfrak h})$ the centralizer in $U_q(\mathfrak{g})$ 
of the subalgebra 
$U_q({\mathfrak h})$, i.e. the subalgebra of zero-weight elements of 
$U_q(\mathfrak{g}).$ 

\bigskip

We now define a completion of $U_q(\mathfrak g)$ which enables us to define
elements (such as $R$) which are expressed in $U_q(\mathfrak{g})$ as an 
infinite series or infinite product, but which evaluation in each finite 
dimensional representation is well defined. 
We denote ${\bf Rep}_{U_q(\mathfrak{g})}$  the category of finite dimensional 
representations of $U_q(\mathfrak{g})$,  
its objects are finite dimensional $U_q(\mathfrak{g})$-modules and 
arrows are interwiners. 
We define ${ \bf Vect}$ to be  the category  of vector spaces.
There is a forgetful functor ${\bf U}:{\bf Rep}_{U_q(\mathfrak{g})}\rightarrow
{ \bf Vect}.$
We now define $\End(\bf{U})$ to be the set of natural transformations from 
${\bf U}$ to ${\bf U}$ which preserve the addition, 
i.e. an element $a\in \End({\bf U})$ is a family of endomorphisms 
$a_V\in \End(V)$ such that:
\begin{alignat}{2}
   &\text{(naturality)} &\quad 
     &\text{for all }  f:V\rightarrow W \text{ intertwiner then }
          a_W\circ f=f\circ a_V,
     \label{naturality}\\
   & \text{(additivity)} &\quad 
     &a_{V\oplus W}=a_{V}\oplus a_{W}\label{additivity}.
\end{alignat}
$\End(\bf{U})$ is naturally endowed with a structure of algebra.
We have a canonical homomorphism of algebra
$U_q(\mathfrak{g})\rightarrow \End(\bf{U})$ which associates to each element 
$a\in U_q(\mathfrak{g}) $ the natural transformation whose $V$ component is 
the element $\pi_V(a).$  We denote  $(U_q(\mathfrak{g}))^c=\End(\bf{U})$ this completion.

We can extend this construction to define a completion of 
$U_q(\mathfrak{g})^{\otimes n}$ as follows.
We define $ {\bf U}^{\otimes n}$ to be the functor from the category 
${\bf Rep}_{U_q(\mathfrak{g})}^{\times n} $ to  
${ \bf Vect}$ which associates to an n-uplet  $(V_1,\ldots,V_n)$ the vector space 
$\otimes_{i=1}^n V_i$.
We define $\End({\bf U}^{\otimes n})$ to be the set of natural transformations 
from ${\bf U}^{\otimes n}$ to ${\bf U}^{\otimes n}$ which are additive in 
each entry.
An element $a$ in $\End({\bf U}^{\otimes n})$ is therefore a family of 
endomorphisms $a_{V_1,\ldots,V_n}\in \End(\otimes_{i=1}^n V_i)$ satisfying the axioms of naturality 
and additivity in  each entry.
$\End({ \bf U}^{\otimes n})$  is naturally endowed with a structure of algebra and we denote as well 
$(U_q(\mathfrak{g})^{\otimes n})^c=\End({ \bf U}^{\otimes n}).$

The coproduct of $U_q(\mathfrak{g})$ therefore defines a morphism of algebras 
from $(U_q(\mathfrak{g}))^c $ to  $(U_q(\mathfrak{g})^{\otimes 2})^c$  which associates to the 
element $a\in (U_q(\mathfrak{g}))^c $ the element 
$\Delta(a)$ 
where $\Delta(a)_{V,W}=a_{V\otimes W}.$

We now define the corresponding completions of $U_q^{\pm}(\mathfrak g),\ U_q(\mathfrak{h}),\ 
{ C}_q({\mathfrak h}).$

Let $V$ be a finite dimensional $U_q(\mathfrak{g})$-module, it is $\mathfrak{h}$ semisimple and we have $V=\bigoplus_{\lambda\in {\mathfrak h}^*}V[\lambda].$ 
An element $a\in \End(V)$ is said strictly upper triangular if 
$aV[\lambda]\subset \oplus_{\lambda'> \lambda} V[\lambda'].$
An element $a\in  \End(V)$ is said to be zero-weight if $aV[\lambda]\subset V[\lambda].$ If $a$ is zero-weight and such that the  restriction of $a$ to $V[\lambda]$ is proportional to $\id_{V[\lambda]},$ $a$ is said to be diagonal.

We define $(U_q^{+}(\mathfrak{g}))^c$ as being the subspace of elements $a\in (U_q(\mathfrak{g}))^c$ such that $a_V$ is strictly upper triangular for all finite dimensional module $V.$
The analog definition holds for strictly lower triangular and this notion defines the subspace $(U_q^{-}(\mathfrak{g}))^c.$
Note that $(U_q^{\pm}(\mathfrak{g}))^c$ are subalgebras (without unit) of 
$(U_q(\mathfrak{g}))^c$ and that the canonical homomorphisms can be restricted to homomorphisms of algebras 
$U_q^{\pm}(\mathfrak{g})\rightarrow (U_q^{\pm}(\mathfrak{g}))^c.$

We define $(U_q(\mathfrak{h}))^c$ as being the 
subalgebra  of elements $a\in (U_q(\mathfrak{g}))^c$ such that $a_V$ is diagonal for all  finite dimensional module $V.$
We also define $({ C}_q({\mathfrak h}))^c$ the subalgebra of  elements $a\in (U_q(\mathfrak{g}))^c$  such that  $a_V$ is zero-weight  for all  finite dimensional module $V.$
As a result we define the subalgebras $(U_q^{}(\mathfrak{b}_\pm))^c=(U_q(\mathfrak{h}))^c\oplus(U_q^{\pm}(\mathfrak{g}))^c.$

\bigskip

$U_q({\mathfrak g})$ is a quasitriangular Hopf algebra with an $R$-matrix 
$R\in (U_q(\mathfrak{b}_+)\otimes U_q(\mathfrak{b}_-))^c$, 
called the standard $R$-matrix, 
which satisfies the following quasitriangularity axioms:
\begin{align}
  &(\Delta \otimes id)(R^{(\pm)})=R^{(\pm)}_{13}\,R^{(\pm)}_{23},
  \qquad
  (id \otimes \Delta)(R^{(\pm)})=R^{(\pm)}_{13}\,R^{(\pm)}_{12},
              \label{quasitriangularity}\\
  &R^{(\pm)}\, \Delta(a) =\Delta'(a)\, R^{(\pm)},\ 
                \forall a\in  U_q(\mathfrak{g}).
\label{quasitriangularity2}
\end{align}
Here we have used the notation  $R^{(+)}=R_{12},\ R^{(-)}=R_{21}^{-1}$.
The explicit expression of $R$ in terms of the root system will be recalled 
further (see Eq.~\eqref{R=KR}).

Moreover $U_{q}({\mathfrak g})$ is a ribbon Hopf algebra, which means
that it exists an invertible element $v\in(U_q(\mathfrak{g}))^c $ such that
\begin{align}
  &\mbox{$v$ is a central element,}\nonumber\\
  & v^2=uS(u), \quad \epsilon(v)=1, \quad S(v)=v, \label{v}\\
  &\Delta(v)=(R_{21}R_{12})^{-1}(v\otimes v).\label{propertyv}
\end{align}
Here, $u$ is the element $u=\sum_i S(b_i)a_i\in (U_q(\mathfrak{g}))^c$, 
where $S$ is the antipode 
and $R=\sum_i a_i \otimes b_i$. It satisfies the properties
\begin{align}
  &S^2(x)=uxu^{-1},\ \forall x \in U_{q}{(\mathfrak{g})},\\
  &\Delta(u)=(R_{21}R_{12})^{-1}(u\otimes u).\label{propertyu}
\end{align}
In this framework, the element $\mu=uv^{-1}$ is a group-like element that we 
choose as follows:
\begin{equation*}
   \mu=q^{2 t_\rho}=q^{2\sum_i \zeta^{\alpha_i}} 
    \quad \text{with}\quad \rho=\frac{1}{2}\sum_{\alpha \in \Phi^+} \alpha.
\end{equation*}
As a result, in a  representation $\pi$  of highest weight $\lambda$, 
$v$ is constant and takes the value $q^{-(\lambda,\lambda+2\rho)}$.

In the case where ${\mathfrak g}=sl(n+1)$, the fundamental representation 
(also called the vector representation) is denoted 
$\stackrel{f}{\pi}: U_q(\mathfrak{g})\rightarrow 
\mathrm{Mat}_{n+1}(\mathbb{C})$ 
and we have 
$\stackrel{f}{\pi}(h_{\alpha_i})=E_{i,i}-E_{i+1,i+1}$,
$\stackrel{f}{\pi}(e_{\alpha_i})=E_{i,i+1}$,
$\stackrel{f}{\pi}(f_{\alpha_i})=
E_{i+1,i},$ where $E_{i,j}$ is the basis of elementary matrices of 
$\mathrm{Mat}_{n+1}(\mathbb{C}).$
The explicit value of $v$ in the fundamental representation is given by 
$\stackrel{f}{\pi}(v)=q^{-\frac{n(n+2)}{n+1}} 1.$

\bigskip

Let us now review the explicit construction of the $R$-matrix using the 
quantum Weyl group.

In order to obtain simple formulas we need to introduce the $q$-exponential 
function.
The $q$-exponential is the meromorphic function
\begin{equation}
   e_q^{z}=\sum_{n=0}^{+\infty}\frac{z^n}{(n)_q!}
          =\frac{1}{((1-q^2)z;q^2)_{\infty}},
   \quad z\in \mathbb{C} \label{defqexp}
\end{equation}
with the notation $(z)_q=q^{z-1}[z]_q=\frac{1-q^{2z}}{1-q^2}.$ 
Elementary properties of the $q$-exponential function, useful to derive
combinatorially the properties of $R$-matrices and the quantum Weyl elements 
in the spirit of \cite{KT},  are now  recalled. 
We have, for any $z\in  {\mathbb C}$, 
\begin{align}
  &e_q^{z}\;e_{q^{-1}}^{-z}=1,\label{expinverse}\\
  &e_q^{q^2z}=e_q^{z}\;(1+(q^2-1)z),\label{expshift}
\end{align}
and for any elements $x,y$
\begin{align}
  &e_q^{x}\;y\;e_{q^{-1}}^{-x}=y+\sum_{k=1}^{+\infty}\frac{1}{(k)_q!}
                           [x,[\cdots,[x,y]]_{q^2}\cdots]_{q^{2k-2}},
\end{align}
where $[x,y]_{q^{2m}}=xy-q^{2m}yx.$

For any  $x,y$ such that  $xy=q^2yx,$
\begin{align}
  &e_q^{x+y}=e_q^{y}\;e_q^{x},\label{expsomme}\\
  &e_q^{x}\; e_q^y=e_q^y \; e_q^{(1-q^{-2})xy}\;e_q^{x}\label{expproduit}.
\end{align}

In the case where $\mathfrak{g}=sl(2)$, we have $R=K \hat{R}$ where 
$K=q^{\frac{h\otimes h}{2}}$ and $\hat{R}=e_{q^{-1}}^{(q-q^{-1})e\otimes f}$.
The quantum Weyl group of $U_q(sl(2))$ is formed by the element
$\hat{w}\in (U_q(sl(2)))^c$ defined as:
\begin{align}
  {\hat w}&=e_{q^{-1}}^{{f}}\, q^{-\frac{h^2}{4}}\,
            e_{q^{-1}}^{{-e}}\, q^{-\frac{h^2}{4}}\,
            e_{q^{-1}}^{{f}}\, q^{-\frac{h}{2}}\\
          &= 
            e_{q^{-1}}^{{-e}}\, q^{-\frac{h^2}{4}}\,
            e_{q^{-1}}^{{f}}\, q^{-\frac{h^2}{4}}\,
            e_{q^{-1}}^{{-e}}\, q^{-\frac{h}{2}}.\label{saitow}
\end{align}
The quantum Weyl group element $\hat{w}$ satisfies the two 
identities:
\begin{align}
   &\Delta(\hat{w})=\hat{R}^{-1}(\hat{w}\otimes\hat{w}),\\
   &\hat{w}^2=q^{\frac{h^2}{2}}\xi v,
\end{align}
where $\xi\in (U_q(sl(2)))^c $ is a central group element 
which value in each 
irreducible finite dimensional representation of dimension $k$ is $(-1)^{k-1}.$

In the general case, let $W$ be the Weyl group associated to 
the root system $\Phi$. For each root $\alpha$, let $s_{\alpha}\in W$ be
its associated reflection. 
For any two distinct nodes $i,j$ of the Dynkin diagram we define an integer 
$m_{ij}$ by  $m_{ij}=2,3,4,6$  respectively if $a_{ij}a_{ji}=0,1,2,3.$ 
The defining relations of $W$ are: 
$s_{\alpha_i}^2=1,\ (s_{\alpha_i}s_{\alpha_j})^{m_{ij}}=1.$

The braid group $\cal{B}$ associated to $W$ is the group generated by 
$\sigma_1,\ldots,\sigma_r$ and satisfying the braid relations
\begin{equation}
   [\sigma_i\sigma_j]^{m_{ij}/2}=[\sigma_j\sigma_i]^{m_{ij}/2}
  \label{braidrelations},
\end{equation}
where we have used the following notation:
if $a,b$ are two elements of a group we define $[ab]^{n/2}$ for $n\geq 0$ to  
be the element $(ab)^{n/2}$  if $n$ is even and $(ab)^{(n-1)/2}a$ if $n$ 
is odd.

To each simple root $\alpha$ we associate an element $\hat{w}_{\alpha}$ as 
follows:
if $U_q(sl(2))_{\alpha}$ is the Hopf subalgebra of $U_q(\mathfrak g)$ 
generated by $e_{\alpha},\ f_{\alpha},\ q^{h}\ (h\in \mathfrak{h}_{\alpha})$, 
we have $U_q(sl(2))_{\alpha}=U_{q_{\alpha}}(sl(2))$ as a Hopf algebra, 
and therefore we can construct the element $\hat{w}_{\alpha}$ of  
$(U_q(sl(2))_{\alpha})^c$ by the same procedure as (\ref{saitow}).
The elements $\hat{w}_{\alpha_i}$   satisfy the braid relations 
({\ref{braidrelations}}).
One therefore obtains a morphism from the braid group $\cal B$ to the group 
of invertible elements of $U_q(\mathfrak g)$ which is called the 
quantum Weyl group \cite{KR,LS}. 

Let $w=s_{\alpha_{i_1}}\ldots s_{\alpha_{i_k}}$ be a reduced expression of an 
element $w\in W$, then the element $\hat{w}=
\hat{w}_{\alpha_{i_1}}\ldots\hat{w}_{\alpha_{i_k}}$ does not depend on the 
choice of the reduced expression.
One can therefore associate to $w$ the automorphism $T_{w}$ of 
$U_q(\mathfrak{g})$ defined as 
$T_{w}(a)=\hat{w}a{\hat{w}}^{-1}.$

Let $w_0$ be the longest element of $W$ and let 
$w_0=s_{\alpha_{i_1}}\ldots s_{\alpha_{i_p}}$ be a reduced expression. 
We have $p=\vert \Phi^+ \vert.$
The set $\{\alpha_{i_1},s_{\alpha_{i_1}}(\alpha_{i_2}),\ldots , 
s_{\alpha_{i_1}}\ldots s_{\alpha_{i_{p-1}}}(\alpha_{i_p})\}$ contains every 
positive root exactly once and defines therefore an order on the set of 
positive root by:
\begin{equation}
   \alpha_{i_1}<s_{\alpha_{i_1}}(\alpha_{i_2})<\cdots <
   s_{\alpha_{i_1}}\ldots s_{\alpha_{i_{p-1}}}(\alpha_{i_p}).
\end{equation}
This ordering of positive roots is a normal ordering in the sense of 
V. Tolstoy \cite{KT}.

Using this ordering, one can now express the standard $R$-matrix of 
$U_q(\mathfrak{g})$ in a similar way as in the $U_q(sl(2))$ case.
Indeed, for $\alpha\in \Phi^+$, there exists a unique $k$ such that 
$\alpha=s_{\alpha_{i_1}}\ldots s_{\alpha_{i_{k-1}}}(\alpha_{i_k})$,
which enables us to define 
$e_{\alpha}=T_{\alpha_{i_1}}\ldots T_{\alpha_{i_{k-1}}}(e_{i_k})$ as well as  
$f_{\alpha}=T_{\alpha_{i_1}}\ldots T_{\alpha_{i_{k-1}}}(f_{i_k})$. 
The algebra generated by 
$e_{\alpha},\ f_{\alpha},\ q^h\ (h\in \mathfrak{h}_{\alpha})$ is 
$U_{q_{\alpha}}(sl(2))$. We therefore define
\begin{equation} 
  \hat{R}_{\alpha}=e_{q_{\alpha}^{-1}}^{(q_{\alpha}-q_{\alpha}^{-1})
                   e_{\alpha}\otimes f_{\alpha}}.
\end{equation}
The standard $R$-matrix of $U_q(\mathfrak{g})$ is expressed as:
\begin{equation}
   R=K{\widehat{R}} \ \mbox{ where }\ 
   K=\prod_{j=1}^r q^{h_{\alpha_j}\otimes\zeta^{\alpha_j}}\ \mbox{ and }\
   {\widehat{R}}=\prod_{\alpha\in \Phi^+}^{>}{\widehat{R}}_{\alpha},
\label{R=KR}
\end{equation}
Its associated classical r-matrix is
\begin{equation}
  r=\frac{1 }{2}\Omega_{\mathfrak{h}}+
    \sum_{\alpha\in \Phi^+} \frac{(\alpha,\alpha)}{2} e_{\alpha} \otimes f_{\alpha}.
\end {equation}

An important result states that the value of the standard $R$-matrix 
\eqref{R=KR} is independent of
the choice of the reduced expression of $w_{0}.$ 
Moreover, we have by construction 
\begin{equation}
  \Delta(\hat{w}_{0})={\widehat{R}}^{-1}(\hat{w}_{0}\otimes \hat{w}_{0}).
\end{equation}
If we define   $\omega=q^{-\frac{h_{\alpha_i}\zeta^{\alpha_i}}{2}}\hat{w}_{0}$,
this element satisfies:
\begin{align}
&\Delta(\omega)=R^{-1} \omega_1 \omega_2\label{deltaw2}.
\end{align}

Let us finally precise simpler  conventions in the case 
$\mathfrak{g}=sl(n+1).$ 
We will use the shortened 
notations $h_{(i)}=h_{\alpha_i},\ \zeta^{(i)}=\zeta^{\alpha_i},\ 
e_{(i)}= e_{\alpha_i},\ f_{(i)}=f_{\alpha_i}$,  as well as 
$w_{\alpha_i}=w_{(i)}.$ 
We will choose the following reduced expression of $w_{0}$:
\begin{equation}
   w_{0}= w_{(1)}(w_{(2)}w_{(1)})\dots (w_{(n)}\dots w_{(1)}),
\end{equation}
which implies the following ordering on roots:
\begin{multline}
   \alpha_1<\alpha_1+\alpha_2<\alpha_2<\cdots\\
   \cdots<\alpha_{n-1}< \alpha_1+\alpha_2+\cdots +\alpha_n
   <\alpha_2+\cdots+\alpha_n<\cdots<\alpha_{n-1}+\alpha_n<\alpha_n.
\end{multline}
If $1\leq i\leq  j\leq n$, we define the positive root 
$\alpha_{ij}=\sum_{k=i}^j\alpha_k.$
We can therefore construct a Cartan-Weyl basis in two different ways.
The first one has already been explained, we denote $e_{\alpha_{ij}},
f_{\alpha_{ij}}$ the corresponding elements.
The second one uses the inductive algorithm of \cite{KT}.
We denote by $e_{(ij)},f_{(ij)}$ the corresponding elements.
% which are recalled in the appendix.
We have $e_{\alpha_{ij}}=(-1)^{i-j}e_{(ij)}$, 
$f_{\alpha_{ij}}=(-1)^{i-j}f_{(ij)}.$

In the fundamental representation, we have 
$\stackrel{f}{\pi}(e_{(ij)})=E_{i,j+1}$, 
$\stackrel{f}{\pi}(f_{(ij)})=E_{j+1,i}$, and the explicit expression
$\mathbf{R}=(\stackrel{f}{\pi}\otimes \stackrel{f}{\pi})(R)$ of the standard
$R$-matrix is given by
\begin{equation}
   \mathbf{R}= q^{-\frac{1}{n+1}}
   \bigg\{q\sum_{i=1}^{n+1} E_{ii}\otimes E_{ii}+
          \sum_{1\leq i\not=j\leq n+1} \!\!\! E_{ii}\otimes E_{jj}
        +(q-q^{-1}) \!\!\! \sum_{1\leq i<j\leq n+1} \!\!\! E_{ij}\otimes E_{ji}\bigg\}.
 \end{equation}

%%%%%%%%%%%%%%%%%%%%%%%%%%%%%%%%%%%%%%%%%%%%%%%%%%%%%%%%%%%%%%%%%%%%%%%%%%%%%%%
\subsection{Belavin-Drinfeld Triples and   Cremmer-Gervais $R$-matrices}

We recall here the notion of Belavin-Drinfeld triple and study the 
$sl(n+1)$ case with a special attention to the shift and to the universal construction of the  Cremmer-Gervais's solution.

A Belavin-Drinfeld triple \cite{BD} for a simple Lie algebra $\mathfrak{g}$ is a triple $(\Gamma_1,\Gamma_2, T)$ where $\Gamma_1,\Gamma_2$ are subsets of the Dynkin diagram $\Gamma $  of $\mathfrak{g}$ and $T:\Gamma_1\rightarrow \Gamma_2$ is an isomorphism which preserves the inner product and satisfies the nilpotency condition: if
$\alpha\in \Gamma_1$ then there exists $k$ such that $T^{k-1}(\alpha)\in \Gamma_1$ but $T^{k}(\alpha)\notin \Gamma_1.$
We extend $T$ to a Lie algebra homomorphism $T:\mathfrak{n}_+ \rightarrow \mathfrak{n}_+$ by setting on simple root elements $T(e_{\alpha})=e_{T(\alpha)}$ for $\alpha\in \Gamma_1$, and zero otherwise.

Any solution $\bf{r}$ to the CYBE satisfying ${ \bf r}+{ \bf r}_{21}=\Omega$ is equivalent, under an automorphism of $\mathfrak{g}$, to a solution of the form
\begin{equation}
r_{T,s}=r-s+\sum_{\alpha\in \Phi^+}\sum_{l=1}^{+\infty}\frac{(\alpha,\alpha)}{2}
T^l(e_{\alpha})\wedge f_{\alpha}, 
\end{equation}
where $s\in \bigwedge^2 \mathfrak{h}$ is a solution of the affine equations:
\begin{equation}
2((\alpha-T\alpha)\otimes id)(s)=((\alpha+T\alpha)\otimes id)(\Omega_{\mathfrak{h}}).
\end{equation}
Given a Belavin-Drinfeld triple, the affine space of the $s$ satisfying the previous equations is of dimension $k(k-1)/2$ where 
$k=\vert \Gamma\setminus \Gamma_1\vert$.

In the present work, we are mainly concerned with the case where ${\mathfrak g}=sl(n+1)$, and with the following Belavin-Drinfeld triple, known as the shift:
\begin{align}
&\Gamma=\{\alpha_1,\ldots,\alpha_{n}\},
\quad 
\Gamma_1=\{\alpha_2,\ldots,\alpha_{n}\},
\quad
\Gamma_2=\{\alpha_1,\ldots,\alpha_{n-1}\},
              \nonumber\\
&\tau:\ \Gamma_1\rightarrow \Gamma_2,\  \alpha_i\mapsto\alpha_{i-1}.
              \label{belavintriple}
\end{align}
Since in this case $k=1$, this Belavin-Drinfeld triple selects a unique $s\in \bigwedge^2 \mathfrak{h}$ which is
\begin{equation}
s=-\frac{1}{2}\sum_{j=1}^{n-1}\zeta^{(j)}\wedge \zeta^{(j+1)}.\label{sCG}
\end{equation}

The quantization of the corresponding $r$-matrix in the fundamental representation is known as the  Cremmer-Gervais's solution \cite{CG,BG}.
A universal construction has been given in \cite{KM} in the case $\mathfrak{g}=sl(3)$, whereas the complete understanding of the explicit expression of the quantization of any solution of  Belavin-Drinfeld type has been given in \cite{ESS}.
We will not follow here this last result, and will instead construct directly the quantization of the $r$-matrix associated to the shift.
The expression that we will obtain is simpler than the one obtained by the method of \cite{ESS}.  

Let us first recall the notion of cocycle.
An invertible element $J\in (U_q(\mathfrak{g})^{\otimes 2})^c$ satisfying
 \begin{equation}
(\Delta \otimes id)(J)J_{12}=(id \otimes \Delta)(J)J_{23},\label{cocycle1}
\end{equation}
is called a cocycle.
If $J$ is a cocycle, $(U_q(\mathfrak{g}))^c$ can be endowed with a new 
quasitriangular Hopf algebra structure with twisted coproduct 
$\Delta^{J}(.)={ J}_{12}^{-1}\Delta(.){ J}_{12}$,  
and twisted $R$-matrix $R^{J}={ J}_{21}^{-1}R{ J}_{12}.$
We define $R^{ J(+)}=R^{J},$ $R^{J(-)}=(R_{21}^{J})^{-1}.$ 

We will now construct the universal twist associated to the shift triple.
We define two morphisms of algebras:
\begin{alignat}{2}
\tau :\ 
   &U_q({\mathfrak n}^+) \rightarrow U_q({\mathfrak n}^+) & &
                   \nonumber\\
   &e_{(i)}\mapsto e_{(i-1)},\ \forall i=2,\ldots n, &
	     \qquad &e_{(1)}\mapsto 0,	          
	            \label{deftau}\\
\tilde{\tau} :\ 
   &U_q({\mathfrak n}^-) \rightarrow U_q({\mathfrak n}^-) & &
                   \nonumber\\
   &f_{(i)}\mapsto f_{(i+1)},\ \forall i=1,\ldots n-1, &
             \qquad &f_{(n)}\mapsto 0.
	            \label{deftautilde}	         
\end{alignat}

\begin{rem}\label{rem-notations}
In order to simplify expressions we will use the notations $e_{(0)}=\zeta^{(0)}=f_{(n+1)}=\zeta^{(n+1)}=0$. 
\end{rem}
We have the following result:

\Theorem{\label{cocycleBCG}}{
For ${\mathfrak g}=sl(n+1)$, a solution $J$ of the cocycle equation
\begin{equation}
  (\Delta \otimes id)(J)\,J_{12}=(id \otimes \Delta)(J)\,J_{23},
  \label{cocycle}
\end{equation}
associated to the shift $\tau$  is given by:
\begin{equation}
J=\prod_{k=1}^{+\infty}J^{[k]} \qquad \text{with}\quad
          J^{[k]}=W^{[k]}\;\widehat{J}{}^{[k]},
\end{equation}
where, $\forall k \in {\mathbb N}^*$,
\begin{align}	  
&\widehat{J}{}^{[k]}=(\tau^{k} \otimes id)(\widehat{R})
    \quad\in 
    \left( 1\otimes1+(U^+_q({\mathfrak g})\otimes U^-_q({\mathfrak g}))^c\right),
    \\
&W^{[k]}_{12}=S^{[k]}(S^{[k+1]})^{-1},
    \qquad 
    S^{[k]}=q^{\sum_{i=1}^{n-k} \zeta^{(i)} \otimes
    \zeta^{(i+k)}}.
\end{align}
This solution will be called  {\em Cremmer--Gervais cocycle}
%\footnote{We have prefered the notation $J$ rather than $J^{\tau}$ because the choice of Belavin-Drinfeld triple %is unambiguous.},
 and $R^J=J^{-1}_{21}RJ_{12}$ will be called  {\em Cremmer--Gervais $R-$matrix}.
}

\begin{rem}
Due to the nilpotency of $\tau$ all the products are actually finite. More precisely, $\widehat{J}{}^{[k]}=S^{[k]}=1 \otimes 1$, $\forall k \geq n$. 
\end{rem}

\begin{rem}
For any $m$ such that $0\leq m\leq k$, one has also $\widehat{J}{}^{[k]}=(\tau^{k-m} \otimes \tilde{\tau}{}^{m})(\widehat{R})$, 
the resulting expression being independent of $m$.
\end{rem}

\proof
Let us give a direct proof of the cocycle identity for $J$ in this framework. 
First, let us remark the following elementary results deduced from the properties of the coproduct:
\begin{equation}
(\Delta \otimes id)(S^{[k]})=S_{13}^{[k]}\, S_{23}^{[k]},\qquad
 (id \otimes \Delta)(S^{[k]})=S_{13}^{[k]}\, S_{12}^{[k]},
 \label{prop11}
\end{equation}
and the properties of the $R$-matrix using $\tau$:
\begin{align}
&(\Delta \otimes id)(\widehat{J}{}^{[k]})
  =\widehat{J}{}^{[0,k]}_{1(2\mid 3}\,\widehat{J}{}^{[k]}_{23},
  \label{propJ1}\\
&(id \otimes \Delta)(\widehat{J}{}^{[k]})
  =\widehat{J}{}^{[k,0]}_{1\mid 2)3}\,\widehat{J}{}^{[k]}_{12},
  \label{propJ2}\\
&\widehat{J}{}^{[k]}_{12} \, \widehat{J}{}^{[k,m]}_{1 (2\mid 3}\, 
   \widehat{J}{}^{[m]}_{23}    =    \widehat{J}{}^{[m]}_{23} \,
   \widehat{J}{}^{[k,m]}_{1 \mid 2) 3}\, \widehat{J}{}^{[k]}_{12},
  \label{propJ3}
\end{align}
where we have defined
\begin{align}
&\widehat{J}{}^{[k,m]}_{1(2\mid 3}
  =(\tau^{k} \otimes id \otimes\tilde{\tau}{}^{m})
                 (K^{-1}_{23}\widehat{R}_{13}\,K_{23}),
           \label{defJ1}\\
&\widehat{J}{}^{[k,m]}_{1\mid 2) 3}
  =(\tau^{k} \otimes id \otimes \tilde{\tau}{}^{m})
                 (K^{-1}_{12}\widehat{R}_{13}\,K_{12}).
           \label{defJ2}
\end{align}
These elements are well defined because  $K^{-1}_{23}\widehat{R}_{13}K_{23}$ and 
$K^{-1}_{12}\widehat{R}_{13}K_{12}$ belong to 
$(U_q({\mathfrak n}^+)\otimes U_q({\mathfrak h})\otimes U_q({\mathfrak n}^-))^c.$
Using the explicit values of $W,J$ and the properties of $\tau$, one can also prove that:
\begin{align}
&[J^{[n-i]}_{23},J^{[j]}_{12}]=0,\quad\forall \;i,j\; / \;1 \leq i \leq j \leq n-1,
\label{prop1}\\
&[\widehat{J}{}^{[k]}_{12},W^{[k+m]}_{13}W^{[m]}_{23}]=[\widehat{J}{}^{[k]}_{23},W^{[k+m]}_{13}W^{[m]}_{12}]=0,
\;\forall k,m \in\{1,\ldots,n\!-\!1 \}\; /\; k\!+\!m \leq n\!-\!1,
\label{prop2}\\
&W^{[m+1]}_{23}\widehat{J}{}^{[l-m-1,m+1]}_{1(2\mid 3}(W^{[m+1]}_{23})^{-1}=W^{[l-m]}_{12}\widehat{J}{}^{[l-m,m]}_{1\mid 2) 3}(W^{[l-m]}_{12})^{-1},\;\forall 0\leq m < l \leq n\!-\!1.
\label{prop3}
\end{align}
\eqref{prop1} is an immediate consequence of the fact that $J^{[i]}_{12}\in A^+_{i} \otimes A^-_{i}$, where $A^+_{i}$ (resp. $A^-_{i}$) is the subalgebra of $U_q({\mathfrak g})$ generated by $q^{\zeta^{(k)}},e_{(k)},\;k=1,\ldots,n-i$ (resp. generated by $q^{\zeta^{(k)}},f_{(k)},\;k=i+1,\ldots,n$). 
In order to check \eqref{prop2} and \eqref{prop3}, we use the following notation: 
\begin{equation}
  W^{[k]}_{12}=\prod_{{i,j\in \{1,\cdots,n\}}} q^{\varepsilon^{k}_{i,j}\, 
            \zeta^{(i)} \otimes \zeta^{(j)}},
\end{equation}
with
\begin{equation}
\varepsilon^{(m)}_{i,j}=\left\{ \begin{array}{lll} 
1\qquad\! \text{if} \quad i,j,m\in \{1,\ldots,n-1\},\quad j=m+i,\\ 
-1\quad \text{if} \quad i,j,m\in \{1,\ldots,n-1\},\quad j=m+i+1,\\ 
0 \qquad\! \text{otherwise.}
\end{array}\right.
\end{equation}
\eqref{prop2} follows immediately from the fact that
\begin{align}
&\varepsilon^{m}_{i,j}=\varepsilon^{m-k}_{i+k,j},\;\; \forall i,j,k,m \in \{1,\ldots,n\!-\!1\}\;/\; i\!+\!k \leq n\!-\!1,1 \leq m\!-\!k,
\label{identiteepsilon1}\\
&\varepsilon^{m}_{i,j}=\varepsilon^{m-k}_{i,j-k},\;\;\forall i,j,k,m \in \{1,\ldots,n\!-\!1\}\;/\; 1 \leq m\!-\!k,1 \leq j\!-\!k,
\label{identiteepsilon2}
\end{align}
and \eqref{prop3} is equivalent to
\begin{equation}
  h_{(p-m-1)}+h_{(p-m)}=\varepsilon^{m+1}_{i,p}\,\zeta^{(i)}
                        +\varepsilon^{l-m}_{p-l,j}\,\zeta^{(j)},
  \quad\forall l,m,p\; /\; 
  1 \leq m\!+\!1 \!\leq\! l \!\leq \! p\!-\!1 \!\leq \!n\!-\!1,
\label{identiteepsilon3}
\end{equation}
which is satisfied due to the relation
\begin{equation*}
h_{(i)}=2\zeta^{(i)}-\zeta^{(i-1)}-\zeta^{(i+1)}.
\end{equation*}

We can now prove the cocycle identity for $J$ by recursion. Indeed, using 
properties \eqref{prop11}--\eqref{prop3}, we deduce easily the following 
recursion relation, proved in Appendix~\ref{sec-lemmas} 
(Lemma~\ref{appendlemma1}):
\begin{multline*}
  \prod_{k=p}^{n-1}\left\{ W^{[k]}_{13}\, W^{[k-p+1]}_{23}\,
    \widehat{J}{}^{[p-1,k-p+1]}_{1(2\mid 3}\,
    \widehat{J}^{[k-p+1]}_{23} \right\} 
  \Big\{ \prod_{k=p}^{n-1}J^{[k]}_{12} \Big\}\\
  =(id \otimes \Delta)(J^{[p]}) 
  \prod_{k=p+1}^{n-1}\left\{ W^{[k]}_{13}\, W^{[k-p]}_{23}\,  
     \widehat{J}{}^{[p,k-p]}_{1 (2 \mid 3}\,
     \widehat{J}^{[k-p]}_{23}\right\}
  \Big\{ \prod_{k=p+1}^{n-1} J^{[k]}_{12} \Big\}
  J^{[n-p]}_{23},
\end{multline*}
and having remarked that
\begin{align*}
  &(\Delta \otimes id)(J) \, J_{12}
   =\prod_{k=1}^{n-1}\left\{ W^{[k]}_{13} \, W^{[k]}_{23}
       \widehat{J}{}^{[0,k]}_{1(2\mid 3}\,
       \widehat{J}^{[k]}_{23} \right\} 
    \Big\{\prod_{k=1}^{n-1}J^{[k]}_{12}\Big\},\\
  &(id \otimes \Delta)(J)\, J_{23}
   =\prod_{k=1}^{n-1}(id \otimes \Delta)(J^{[k]})
 \;\prod_{k=1}^{n-1}J^{[k]}_{23},
\end{align*}
we conclude the proof of the cocycle identity verified by $J.$ 
\qed

In the following, we extend the morphisms of algebras $\tau$ and $\tilde\tau$ 
to $U_q(\mathfrak{b}^+)$ and $U_q(\mathfrak{b}^-)$ respectively 
as
\begin{alignat}{2}
   &\tau(\zeta^{(i)})=\zeta^{(i-1)},\ i=2,\ldots,n,
     & \qquad &\tau(\zeta^{(1)})=0,
      \label{deftau1}\\
   &\tilde\tau(\zeta^{(i)})=\zeta^{(i+1)},\ i=1,\ldots,n-1,
     &\qquad &\tilde\tau(\zeta^{(n)})=0.
      \label{deftautilde1}      
\end{alignat}
Note that, with this definition,  $\tau$ and $\tilde\tau$ are not 
morphisms of Hopf algebras and are different from the extension of   \cite{ESS}.
Indeed, their action on the coproduct is given as
\begin{alignat}{2}
   &(\tau\otimes\tau)\left(\Delta(a)\right) = q^{\zeta^{(n-1)}\otimes\zeta^{(n)}}
        \,\Delta\left(\tau(a)\right)\, q^{-\zeta^{(n-1)}\otimes\zeta^{(n)}},
	\quad
	& &\forall a \in U_q(\mathfrak{b}^+),
	  \label{deltatau}\\
   &(\tilde\tau\otimes\tilde\tau)\left(\Delta(a)\right) = 
   q^{\zeta^{(1)}\otimes\zeta^{(2)}}
        \,\Delta\left(\tilde\tau(a)\right)\, q^{-\zeta^{(1)}\otimes\zeta^{(2)}},
	& &\forall a \in U_q(\mathfrak{b}^-).
	  \label{deltatautilde}	  
\end{alignat}
Using this definition, we have
\begin{equation}\label{Stau}
  S^{[k]}=(\tau^k\otimes id)(S^{[0]}), \qquad\text{with}\quad
   S^{[0]}=q^{\sum_{i=1}^n \zeta^{(i)}\otimes\zeta^{(i)}}.
\end{equation}

\bigskip

In the fundamental representation,
we denote by $\mathbf{R}^J$ the explicit $(n+1)^2\times (n+1)^2$  
Cremmer-Gervais's solution of the Quantum Yang-Baxter Equation associated to 
the previous twist. It is  given by
\begin{equation}\label{RJ-fund}
{\mathbf R}^J= (D\otimes D)\; \widetilde{\mathbf{R}}^J\; (D\otimes D)^{-1},
\end{equation}
where
\begin{align}
&D=\sum_i\; q^{\frac{i^2-3i}{2(n+1)}}\; E_{i,i},\\
&\widetilde{\mathbf{R}}^J=q^{-\frac{1}{n+1}}\bigg\{\, q\, \sum_{r,s} \; q^{\frac{2(r-s)}{n+1}}  \; E_{r,r}\otimes E_{s,s}
  \nonumber\\
  &\hspace{4.5cm}
  +(q-q^{-1})\sum_{i,j,k}\; q^{\frac{2(i-k)}{n+1}}\; \eta(i,j,k)\;
                              E_{i,j+i-k}\otimes E_{j,k} \bigg\},
\end{align}
with
\begin{equation}
\eta(i,j,k)=\left\{ \begin{array}{lll} 
1\qquad\! \text{if} \quad i\le k<j,\\ 
-1\quad \text{if} \quad j\le k<i,\\ 
0 \qquad\! \text{otherwise.}
\end{array}\right.
\end{equation}

\begin{rem}
$R^J$ is not $\mathfrak{h}$-invariant, but is $t_\rho$-invariant. 
This implies that $\mathbf{R}^J$ is homogeneous in the sense that 
$(\mathbf{R}^J)_{ij}^{kl}\ne 0$ only if $i+j=k+l$.
\end{rem}

%%%%%%%%%%%%%%%%%%%%%%%%%%%%%%%%%%%%%%%%%%%%%%%%%%%%%%%%%%%%%%%%%%%%%%%%%%%%%%
\section{Results on  Dynamical Quantum Groups}
\label{sectionDQG}
%%%%%%%%%%%%%%%%%%%%%%%%%%%%%%%%%%%%%%%%%%%%%%%%%%%%%%%%%%%%%%%%%%%%%%%%%%%%%%

%%%%%%%%%%%%%%%%%%%%%%%%%%%%%%%%%%%%%%%%%%%%%%%%%%%%%%%%%%%%%%%%%%%%%%%%%%%%%%
\subsection{Quantum Dynamical Yang-Baxter Equation}

We first begin with a formulation of the dynamical Yang--Baxter equation which 
is not sufficiently mathematically precise for our future purposes.

\Definition{Quantum Dynamical Yang-Baxter Equation (Formal)}{A universal
solution of  the {\em Quantum Dynamical Yang--Baxter
Equation (QDYBE)}, also known as
{\em Gervais--Neveu--Felder 
equation}, is a  map $R:\mathbb{C}^r\rightarrow 
U_q(\mathfrak{g})^{\otimes 2}$ such that
$R(x)$ is $\mathfrak{h}$-invariant 
and
\begin{equation}
  R_{12}(x)R_{13}(xq^{h_2})R_{23}(x)=
  R_{23}(xq^{h_1})R_{13}(x)R_{12}(xq^{h_3}) \label{sYB1}
\end{equation}
where we have denoted 
$xq^{h}=(x_1q^{h_{\alpha_1}},\ldots,x_r
q^{h_{\alpha_r}})$.}

This is sufficient for formal manipulations but it is not enough precise in 
the sense that the standard universal solution of the dynamical Yang-Baxter 
equation is such that $R(x)$ is an infinite  series in $U_q(\mathfrak{g})$ with
coefficients being rational function of $x_1,\ldots,x_r$ with coefficients in 
$U_q(\mathfrak{h}).$ 

As a result we can extend the construction of 
$\End({\bf U})$ and $ \End({\bf U}^{\otimes 2})$ as follows.
Let $A$ be a unital  algebra over the complex field, we define the
 functor $A\otimes {\bf U}:{\bf Rep}_{U_q(\mathfrak{g})}\rightarrow 
A\otimes{\bf Vect}$, where $A\otimes{\bf Vect}$ is the category which 
objects are $A\otimes V$ where $V$ is a vector space and the maps are 
$\id_A\otimes\phi$ where $\phi$ is a linear map between vector spaces.
We can define the functor $A\otimes {\bf U}$ which associates to each finite dimensional 
 module $V$ the vector space $A\otimes V$ and to each interwiner $\phi$ the map $\id_A\otimes \phi.$
We define 
$\End(A\otimes {\bf U})$ to be the set of additive natural transformations 
between the functors 
$A\otimes {\bf U}$ and $A\otimes {\bf U}.$ An element of
 $\End(A\otimes {\bf U})$ is a family 
$a_V\in A\otimes \End(V)$ such that it satisfies the naturality and 
additivity condition. 
We can define similarly  $\End(A\otimes {\bf U}^{\otimes n}).$

  $\End(A\otimes {\bf U}^{\otimes n})$ is  
naturally endowed with a structure of algebra.
We will denote $(A\otimes U_q(\mathfrak g)^{\otimes n})^c=\End(A\otimes {\bf U}^{\otimes n}).$

We could similarly have defined 
$( U_q(\mathfrak g)^{\otimes n}\otimes A)^c.$

\Definition{Quantum Dynamical Yang-Baxter Equation (Precise)}{A universal
solution of  the {\em Quantum Dynamical Yang--Baxter
Equation (QDYBE)}, also known as
{\em Gervais--Neveu--Felder 
equation}, is an $\mathfrak{h}$-invariant element $R(x)$ of 
$(\mathbb{C} (x_1,\ldots,x_r)\otimes U_q(\mathfrak g)^{\otimes 2})^c$ satisfying 
\begin{equation}
  R_{UV}(x)R_{UW}(xq^{h_V})R_{VW}(x)=
  R_{VW}(xq^{h_U})R_{UW}(x)R_{UV}(xq^{h_W}) \label{sYB}
\end{equation}
where we have denoted 
$xq^{h}=(x_1q^{h_{\alpha_1}},\ldots,x_r q^{h_{\alpha_r}})$ and where
$U,V,W$ are any finite dimensional $U_q(\mathfrak g)$-modules.}

If moreover $R(x)$ belongs to 
$(\mathbb{C} (x_1^2,\ldots,x_r^2)\otimes  U_q(\mathfrak g)^{\otimes 2})^c$ 
it is called $x^2$-rational. 

\bigskip

We will often study the particular case where $\mathfrak{g}=sl(n+1).$

\Definition{}{
Let $\stackrel{f}{\pi}$ be the fundamental representation of 
${\mathfrak g}=U_q(sl(n+1)).$ Once a $n$-uplet  $(x_1,\ldots,x_{n})$ of 
complex numbers is given, we will denote 
\begin{equation}
\stackrel{f}{\pi}\bigg(\prod_{i=1}^n x_i^{2\zeta_{\alpha_i}}\bigg)
=\mathrm{diag}(\nu_1,\ldots,\nu_{n+1}).\label{deflambda}
\end{equation}
 As a result we have   $\prod_{i=1}^{n+1} \nu_i=1$ and $x^2_i=\nu_i\nu_{i+1}^{-1}.$
}

Therefore, in the $sl(n+1)$ case, we can define a notion of regularity as 
follows:

Let  $\nu_1,\ldots,\nu_n$ be $n$ indeterminates and define $\nu_{n+1}$ by  
$\prod_{i=1}^{n+1} \nu_i=1$ and  $x^2_i=\nu_i\nu_{i+1}^{-1}$, 
for $i=1,\ldots,n.$
An element $a\in 
(\mathbb{C} (\nu_1,\ldots,\nu_n)\otimes  U_q(\mathfrak g)^{\otimes n})^c$ will be called 
{\em $\nu$-rational}.
An element  $a\in (\mathbb{C} (x_1^2,\ldots,x_n^2)\otimes U_q(\mathfrak g)^{\otimes n})^c $ will be 
called {\em $x^2$-rational}.
Note that obviously an $x^2$-rational element   is also $\nu$-rational.

This distinction is important because in the $U_q(sl(n+1))$ case the standard 
solution of QDYBE   is $x^2$-rational, whereas the universal solution of the 
dynamical coboundary equation is (almost) $\nu$-rational.

%%%%%%%%%%%%%%%%%%%%%%%%%%%%%%%%%%%%%%%%%%%%%%%%%%%%%%%%%%%%%%%%%%%%%%%%%%%%%%%
\subsection{Quantum Dynamical coCycles}

The first to understand universal aspects of the dynamical Yang-Baxter equation
was O. Babelon in his work on quantum Liouville on a lattice \cite{Bab}. 
There he  introduced the notion of {\em Quantum Dynamical coCycle} 
$F(x)\in (\mathbb{C}(x^2)\otimes U_q(sl(2))^{\otimes 2})^c$ and  gave  an exact formula 
for $F(x)$ expressed as a series. 

More generally, in $U_q(\mathfrak{g})^{\otimes 2}$, 
a universal solution of QDYBE equation can be obtained from a 
solution of the {\em Quantum Dynamical coCycle Equation (QDCE)} 
(\ref{eq:s-cocycle}) in the following sense: 

\Theorem{Quantum Dynamical Cocycle Equation\label{TheoremQDCE}}{
If $F(x)\in (\mathbb{C}(x_1^2,\ldots,x_r^2)\otimes U_q(\mathfrak{g})^{\otimes 2})^c $ is an  
$x^2$-rational map such that
\begin{enumerate}
\item $F(x)$ is invertible
\item  $F(x)$  is $\mathfrak{h}$-invariant , i.e.
\begin{equation}
  [F_{12}(x),h\otimes 1+1\otimes h]=0,\ \forall h\in \mathfrak{h},
  \label{hinvariance}
\end{equation}
\item $F(x)$ satisfies the  {\em Quantum Dynamical coCycle Equation (QDCE)},
\begin{equation}
  (\Delta\otimes \id)(F(x))\; F_{12}(xq^{h_3})=
  (\id\otimes \Delta)(F(x))\; F_{23}(x),
  \label{eq:s-cocycle}
\end{equation}
\end{enumerate}
then
\begin{equation}
   R(x)=F_{21}(x)^{-1}R_{12}\,F_{12}(x)
   \label{R=FRF}
\end{equation} 
satisfies the universal QDYBE \eqref{sYB} where 
($R$ is the standard universal $R-$matrix in $(U_q(\mathfrak{g})^{\otimes 2})^c$ defined by (\ref{R=KR}) ).
}

Although simpler than the QDYBE equation, (\ref{eq:s-cocycle}) is difficult 
to solve directly.
However, it is possible to determine its general solutions through 
an auxiliary linear equation, the 
{\em Arnaudon-Buffenoir-Ragoucy-Roche Equation (ABRR)}
(see \cite{ABRR}, but it was first remarked in \cite{BR1}):

\Theorem{ABRR Equation\label{TheoremABRR}}{ 
We define $B(x)\in ({\mathbb A}_r(x_1,\ldots,x_r)\otimes U_q(\mathfrak{h}))^c$ by 
\begin{equation}\label{Bx}
  B(x)=\prod_{j=1}^r 
      x_{j}^{2\zeta^{\alpha_j}}q^{h_{\alpha_j}\zeta^{\alpha_j}},
\end{equation}
where  ${\mathbb A}_r(x_1,\ldots,x_r)=
\otimes_{i=1}^r{\mathbb A}[x_i]$ with ${\mathfrak A}[x]$ being the algebra 
generated by $x^{\alpha},\ \alpha\in {\mathbb R}$, with the relations 
$x^{\alpha+\alpha'}=x^{\alpha}x^{\alpha'}.$\\
An element $F(x)\in ({\mathbb C}(x_1^2,\ldots,x_r^2)\otimes U_q(\mathfrak{g})^{\otimes 2})^c$ satisfies the  {\em ABRR Equation}
if and only if 
\begin{equation}
  F_{12}(x)\, B_{2}(x)={\widehat{R}}_{12}^{-1}\,B_{2}(x)\,F_{12}(x).
  \label{lineareq}
\end{equation}
The importance of the ABRR equation comes from the following theorem \cite{ABRR}:\\
Under the hypothesis 
\begin{align*}
  &
(F(x) - 1\otimes 1)\; \in (\mathbb{C}(x_1^2,\ldots,x_r^2) \otimes U_q ({\mathfrak g})\otimes U^-_q ({\mathfrak g}))^c, 
 \end{align*}
there exists a unique solution of Equation (\ref{lineareq}).
This solution is invertible, $\mathfrak{h}-$ invariant and satisfies the {\em QDCE} (\ref{eq:s-cocycle}). It is called the standard solution of the {\em QDCE}. \\
Let $V$ be a finite dimensional $U_q(\mathfrak g)$-module, there exists a positive number $c_V$ such that, if  $x_1,\ldots,x_r\in {\mathbb C}$ with $\vert x_i\vert< c_V$,  the infinite product  $\prod_{k=0}^{+\infty}\big( B_{2}^{-k-1}(x)\, {\widehat{R}}_{12}\,
                 B_{2}^{k+1}(x)\big)$ is convergent when represented on $V$.
It satisfies moreover the $ABRR$ Equation and belongs to $1\otimes 1+(\mathbb{C}(x_1^2,\ldots,x_r^2) \otimes U_q^+ ({\mathfrak g})\otimes U^-_q ({\mathfrak g}))^c$. As a result, the standard solution $F(x)$ of the {\em QDCE} is given by:
\begin{equation}
 F(x)=\prod_{k=0}^{+\infty}\Big( B_{2}^{-k-1}(x)\, {\widehat{R}}_{12}\,
                 B_{2}^{k+1}(x)\Big).
  \label{Fprod}
\end{equation} 
}
\begin{rem} 
$B(x)$ satisfies the following useful relations:
\begin{align}
   &R_{12}(x)\,B_{2}(x)=B_{2}(x)\,K_{12}^2 \,R_{21}(x)^{-1},
        \label{Rlinear1}\\
   &\Delta(B(x))=B_1(x)\,B_2(x)\,K^2=B_1(xq^{h_2})\,B_2(x)
        =B_1(x)\,B_2(xq^{h_1}).
 \label{DeltaB}
\end{align}
\end{rem}

\bigskip

In the fundamental representation $\stackrel{f}{\pi}$ of $U_q(sl(n+1))$, 
let us denote by $\mathbf{F}(x)$ the explicit expression of the quantum 
dynamical cocycle \eqref{Fprod}. It is given by
\begin{equation}\label{Fx-fund}
  \mathbf{F}(x)=1\otimes 1 - (q-q^{-1}) \sum_{i<j} 
  \Big(1-\frac{\nu_j}{\nu_i}\Big)^{-1} E_{i,j}\otimes E_{j,i}.
\end{equation}
The corresponding expression \eqref{R=FRF} of the standard dynamical 
$R$-matrix for $sl(n+1)$ in the fundamental representation, denoted 
${\mathbf R}(x)$,
is then
\begin{multline}\label{Rx-fund}
  {\mathbf R}(x)=q^{-\frac{1}{n+1}} \bigg\{
     q\sum_{i}E_{i,i}\otimes E_{i,i}
     + \sum_{i \not= j}E_{i,i}\otimes E_{j,j}
     +(q-q^{-1})\sum_{i\not= j} \Big(1-\frac{\nu_i}{\nu_j}\Big)^{-1}
                         E_{i,j}\otimes E_{j,i}\\
     - (q-q^{-1})^2\sum_{i>j} 
             \frac{\nu_i}{\nu_j}\Big(1-\frac{\nu_i}{\nu_j}\Big)^{-2}
                         E_{i,i}\otimes E_{j,j}\bigg\}.
\end{multline}

\bigskip

The construction of Theorem \ref{TheoremABRR} has been used by P. Etingof, 
T. Schedler and O. Schiffmann \cite{ESS} to build the quantization 
of $r$-matrices associated to any Belavin-Drinfeld triple. 
Indeed, for any Belavin-Drinfeld triple $T$, they have constructed 
explicitely a twist $J^{T}\in (U_q({\mathfrak g})^{\otimes 2})^c$ 
which satisfies
\begin{equation}
  (\Delta \otimes id)(J^{T})\,J^{T}_{12}
  =(id \otimes \Delta)(J^{T})\,J^{T}_{23},
  \label{cocycle2}
\end{equation}
such that $(J^{T}_{21})^{-1}RJ^{T}_{12}$ is a solution of the Yang-Baxter 
equation and is a quantization of the classical $r^{{}_{T}}$ matrix 
associated to $T$.
The general expression for $J^{T}$ was obtained through a nice use of dynamical
quantum groups and of a modification of the ABRR equation.

In general, for any finite dimensional simple Lie algebra ${\mathfrak g}$ and 
any Belavin-Drinfeld triple $T$, $J{}^{T}$ is expressed as a finite product of 
explicit invertible elements as
\begin{align}
  J{}^{T}=S{}^{T}\widehat{J}{}^{T}\qquad\text{with}\quad
  S^{T}\in (U_q({\mathfrak h})^{\otimes 2})^c,\quad
  \widehat{J}{}^{T}\in  1\otimes1
                    + \left(U^+_q({\mathfrak g})\otimes U^-_q({\mathfrak g})\right)^c.
\end{align} 
The reader is invited to read \cite{ESS} for a general construction of 
$S{}^{T}$ and $\widehat{J}{}^{T}$ using the ABRR identity. 
The result of \cite{ESS} is completely general. However, for reasons explained 
below, we will only be interested in the case where 
${\mathfrak g}=sl(n+1),$ $n\in {\mathbb N}^*$, and where $T=\tau$ is the 
 Cremmer-Gervais triple. 
The construction of $J^{\tau}$ described in the last section is therefore enough for our 
purpose.

%%%%%%%%%%%%%%%%%%%%%%%%%%%%%%%%%%%%%%%%%%%%%%%%%%%%%%%%%%%%%%%%%%%%%%%%%%%%%%%
\subsection{Quantum Dynamical coBoundary Problem}

%Since the beginning of the study of dynamical quantum structures, it has been suggested to write the dynamical cocycle $F(x)$ as a dynamical coboundary. This leads to a new equation, 
%%
%\begin{equation}
%F(x)=\Delta(M(x)){ J} M_2(x)^{-1}(M_1(xq^{h_2}))^{-1},
%\end{equation}
%%
%for a certain  map $M$ from ${\mathbb C}^r$ to $U_q({\mathfrak g})$, associated to a certain solution $F$ of the Dynamical Cocycle Equation (\ref{eq:s-cocycle}), and to a cocycle $J$ in $U_q({\mathfrak g})^{\otimes 2}$ obeying (\ref{cocycle}). Inspired by its explicit value in the fundamental representation we will search $M$ with the property that $\exists M^{(0)}$ mapping from ${\mathbb C}^r$ to $U_q({\mathfrak h})$ such that $M^{(0)\;-1}M$ is $\nu$-rational.

We have seen in the last section that Theorem \ref{TheoremABRR} provides a way
to solve the QDCE \eqref{eq:s-cocycle}.  
However in \cite{Bab}, for the $U_q(sl(2))$ case, O. Babelon  used a 
different approach: he noticed that $F(x)$ is a quantum dynamical 
coboundary, i.e. that there exists an explicit invertible element 
$M(x)\in (\mathbb{C}(x)\otimes U_q(sl(2)))^c$ solution of the following Dynamical coBoundary 
Equation
:
\begin{equation}
   F(x)=\Delta(M(x))\; M_2(x)^{-1}\,(M_1(xq^{h_2}))^{-1}.
\end{equation}
Note that in the  $U_q(sl(2))$ this is a particular case of \eqref{coboundaryM} since ${J}$ is merely equal
 to $1.$

%

%Quite remarkably, O. Babelon obtained an explicit  formula for the solution $F(x)$ in the $U_q(sl(2))$ case, equivalent to previous one but expressed as an infinite sum, by noting that $F(x)$ is a quantum dynamical coboundary, i.e there exists an explicit invertible element $M(x)\in U_q(sl(2))$ solution of dynamical coboundary equation (\ref{coboundaryM}) (where ${ J}$ is simply equal to $1$ in the $U_q(sl(2))$ case). Indeed, Babelon's proof was inspired by the following property:

More generally, we have the following property which is formal at this stage because we do not 
precise the analytic property of $M(x)$: 

\Theorem{Quantum Dynamical coBoundary equation (QDBE) (Formal)\label{coboundcocycl}}{
Let $J$ be a given cocycle in $U_q({\mathfrak g})^{\otimes 2}$, 
if  there exists a map $M$ from ${\mathbb C}^r$ to $U_q({\mathfrak g}),$ 
such that $M(x)$ is invertible for any $x$ where it exists, 
and such that the map $F$ defined by 
\begin{equation}
  F(x)=\Delta(M(x))\;{ J}\; M_2(x)^{-1}\,(M_1(xq^{h_2}))^{-1}
  \label{coboundaryM}
\end{equation}
 verifies
\begin{equation}\label{hinv}
  [F_{12}(x),h\otimes 1+1\otimes h]=0,\quad\forall h\in \mathfrak{h},
\end{equation}
then $F$ is a solution of the Quantum Dynamical coCycle Equation 
(\ref{eq:s-cocycle}) (see \cite{EN} for details). 
As a result, the corresponding quantum dynamical $R$-matrix \eqref{R=FRF}
can be obtained in terms of $M$ as
\begin{equation}
  R(x)= M_2(xq^{h_1})\, M_1(x)\;R^{J}\;M_2(x)^{-1}\,M_1(xq^{h_2})^{-1}.
  \label{RMM=MMR}
\end{equation}
}

\negvspace

\proof
If $F$ satisfies (\ref{coboundaryM}) and (\ref{hinv}), we have, 
using the cocycle equation \eqref{cocycle}:
\begin{align*}
  &(\Delta \otimes  id)(F(x)^{-1})\,(id \otimes \Delta)(F(x))\\
  &\qquad =\Delta_{12}(M(xq^{h_3}))\,M_3(x)\,
           (\Delta \otimes id)(J^{-1})\,(id \otimes \Delta)(J)\,
           \Delta_{23}(M(x)^{-1})\,M_1(xq^{h_2+h_3})^{-1}\\
  &\qquad =\Delta_{12}(M(xq^{h_3}))\,M_3(x)\,J_{12}\,J_{23}^{-1}\,
           \Delta_{23}(M(x)^{-1})\,M_1(xq^{h_2+h_3})^{-1}\\
  &\qquad =\Delta_{12}(M(xq^{h_3}))\,J_{12}\,M_{2}(xq^{h_3})^{-1}\,
           M_{1}(xq^{h_2+h_3})^{-1}\, M_{1}(xq^{h_2+h_3})\\
  &\hspace{4.9cm}\times M_{2}(xq^{h_3})\,M_3(x)\,J_{23}^{-1}\,
           \Delta_{23}(M(x)^{-1})\,M_1(xq^{h_2+h_3})^{-1}\\
  &\qquad =F_{12}(xq^{h_3})\, M_{1}(xq^{h_2+h_3})\,F_{23}(x)^{-1}\,
           M_1(xq^{h_2+h_3})^{-1}\\
  &\qquad=F_{12}(xq^{h_3})\, F_{23}(x)^{-1},
\end{align*}
and thus $F$ is a solution of \eqref{eq:s-cocycle}. 
\eqref{RMM=MMR} follows directly from \eqref{R=FRF} and \eqref{coboundaryM}
using \eqref{quasitriangularity2}.
\qed

Since the article of O. Babelon on the $sl(2)$ case, 
the theory of dynamical quantum groups 
has been the subject of numerous works, and is now well 
understood for $U_q({\mathfrak g})$ where ${\mathfrak g}$ is a Kac-Moody 
algebra of affine or finite type. 
However, quite surprisingly, only little progress has been made 
concerning the dynamical coboundary equation (see however 
the articles \cite{BBB},\cite{BR2},\cite{EN},\cite{St}).

The first results in this subject concerning higher rank cases 
have been obtained by  Bilal-Gervais in \cite{BG}, where they have 
found the expression for $R(x)$ in the fundamental representation of
$sl(n+1)$. Cremmer-Gervais \cite{CG} have then shown that, 
in the fundamental representation of $sl(n+1)$, 
it is possible to absorb the dynamical dependence of $R(x)$ through a 
dynamical gauge transformation. 
More precisely, we have the following theorem: 

\Theorem{\cite{CG},\cite{BDF},\cite{Hodges}\label{lemme1}}{
In the fundamental representation of $U_q({sl(n+1)})$, 
the expression $\mathbf{R}(x)$ \eqref{Rx-fund} of the standard dynamical 
$R$-matrix is related to the expression $\mathbf{R}^J$ \eqref{RJ-fund} 
of the  Cremmer-Gervais $R$-matrix as,
\begin{equation}\label{RMM-fund}
  \mathbf{R}(x)= {\mathbf M}_2(xq^{h_1}) {\mathbf M}_1(x)\;
      {\mathbf R}^{J}\;{\mathbf M}_2(x)^{-1}{\mathbf M}_1(xq^{h_2})^{-1},
\end{equation}
where the $(n+1)\times (n+1)$  matrix ${\mathbf M}(x)$ is given by
\begin{equation}\label{Minv-fund}
   {\mathbf M}(x)^{-1}=D\;\mathcal{V}(x)\;{\mathcal U}(x)D^{-1}
\end{equation}
in term of the Vandermonde matrix
\begin{equation*}
   \mathcal{V}(x)=\sum_{i,j}\;\nu_j^{i-1}\; E_{i,j},
\end{equation*}
and of the following diagonal matrices 
\begin{align*}
  &D=\sum_i\, q^{\frac{i^2-3i}{2(n+1)}}\, E_{i,i},\\
  &\mathcal{U}(x)=\sum_i \frac{\nu_i^{-i+1}\prod_{k=2}^{n+1}\nu_k^{-\frac{1}{2}(\delta_{i\leq k-1}-\frac{k-1}{n+1})}q^{\frac{1}{2}(\delta_{i\leq k-1}-\frac{k-1}{n+1})^2}}{\prod_{r=i+1}^{n+1}(1-\frac{\nu_i}{\nu_r})}E_{i,i}.
\end{align*}
}

\negvspace

\proof
{}From \eqref{Rx-fund}, it is easy to check that
\begin{equation}
  \mathbf{R}(x)=\mathcal{U}_2(x\, q^{h_{1}})^{-1}\,\mathcal{U}_1(x)^{-1}\;
  \widetilde{\mathbf{R}}(x)\;\mathcal{U}_2(x) \, \mathcal{U}_1(x\, q^{h_{2}}),
\end{equation}
where
\begin{multline}
 \widetilde{\mathbf{R}}(x)=q^{-\frac{1}{n+1}} \bigg\{
   q\sum_{i}E_{i,i}\otimes E_{i,i}
  +\sum_{i \not= j} \Big(q-q^{-1}\frac{\nu_i}{\nu_j}\Big)
                     \Big(1-\frac{\nu_i}{\nu_j}\Big)^{-1}  
                                       E_{i,i}\otimes E_{j,j}\\
  +(q-q^{-1})\sum_{i\not= j} \Big(1-\frac{\nu_i}{\nu_j}\Big)^{-1}
                                       E_{i,j}\otimes E_{j,i}\bigg\}.
\end{multline}
Note that this expression is the one that can be found in \cite{BG,ABB}

Therefore, it remains to show that
\begin{equation}
\widetilde{\mathbf{R}}(x)= \mathcal{V}_2(xq^{h_1})^{-1} \, \mathcal{V}_1(x)^{-1}\;\widetilde{\mathbf R}^{J}\; \mathcal{V}_2(x)\, \mathcal{V}_1(xq^{h_2}),
\end{equation}%
which is proved in \cite{CG,BDF,Hodges}.
\qed

The following theorem \cite{BDF} shows that it is only in the case where $\mathfrak{g}=sl(n+1)$ that the coboundary equation can be eventually solved.

\Theorem{\cite{BDF}\label{lemme2}}{
Let $\mathfrak{g}$ be a finite dimensional simple Lie algebra, and let 
$R:\mathbb{C}^r\rightarrow U_q(\mathfrak{g})^{\otimes 2}$ be the standard 
solution \eqref{R=FRF} of the QDYBE. 
If  $\mathfrak{g}$ is not of $A$-type, it is not possible to find any pair 
$(R^J,M)$, with $M:\mathbb{C}^r\rightarrow U_q(\mathfrak{g})^{\otimes 2}$ and 
$R^J$ a solution of the (non-dynamical) QYBE, such that \eqref{RMM=MMR} is 
satisfied. If  $\mathfrak{g}$ is of $A$-type and if such a pair exists, then 
$R^J$ can be expanded as $R^J=1+\hbar\, r_J+o(\hbar)$ where $r_J$  coincides, up to an automorphism of the Lie algebra, with $r_{\tau,s}$ associated to the shift. 
} 

\proof
We refer the reader to \cite{BDF} for the proof of this theorem. 
For completeness we have given an explanation of this proof in Appendix~\ref{sec-BDF}.
\qed

Therefore, quite surprisingly, except for $n=1,$ ${\mathbf R}^{J}$  is not the Drinfeld's solution of the Yang-Baxter equation and is not of
zero weight, but is instead the  Cremmer-Gervais $R$-matrix.
We will therefore assume in the rest of this section that $\mathfrak{g}=sl(n+1).$

\bigskip

{}From the expression of ${\bf M}(x)$ in the fundamental representation one sees, because of the expression of ${\cal U}$, that ${\bf M}(x)$ lies in
 ${\mathbb C}(\tilde{\nu_1},\ldots,\tilde{\nu_n})\otimes \End(\stackrel{f}{V})$, where $\nu_i=\tilde{\nu_i}^{2(n+1)}.$
The next section and the explicit form of ${\bf M}(x)$ motivates the following definition:
  
An element $a\in ({\mathbb C}(\tilde{\nu_1},\ldots,\tilde{\nu_n})\otimes U_q(\mathfrak{g}))^c$ is said to be {\em almost $\nu$-rational} if there exists an invertible element $b\in ({\mathbb C}(\tilde{\nu_1},\ldots,\tilde{\nu_n})\otimes U_q(\mathfrak{h}))^c$ such that 
\begin{enumerate}
\item $b^{-1}a$ is $\nu$-rational,
\item $\Delta(b)\, b_2(x)^{-1}b_1(xq^{h_2})^{-1}$ is $\nu$-rational,
\end{enumerate}
where $\nu_i=\tilde{\nu_i}^{2(n+1)},\ i=1,\ldots,n+1$, and 
$\tilde{\nu}_1\ldots\tilde{\nu}_{n+1}=1.$

For ${\mathfrak g}=sl(n+1)$ and for the  Cremmer-Gervais cocycle $J\in (U_q(\mathfrak{g})^{\otimes 2})^c$, we are now ready to address the following precise problems which will be solved in this paper:

\Problem{Weak Quantum Dynamical coBoundary Problem (WQDBP)\label{WP}}{
Find an almost $\nu$-rational invertible element   
${\cal M}\in ({\mathbb C}(\tilde{\nu_1},\ldots,\tilde{\nu_n})\otimes  U_q({\mathfrak g}))^c$ such that 
\begin{equation}\label{eq1-pb1}
  \mathcal{F}(x)=\Delta(\mathcal{M}(x))\, J \, \mathcal{M}_2(x)^{-1}
           \mathcal{M}_1(xq^{h_2})^{-1}
\end{equation}
is a  zero-weight solution of the QDCE \eqref{eq:s-cocycle} 
and such that 
\begin{equation}\label{eq2-pb1}
   {\cal R}(x)={\cal F}_{21}(x)^{-1} R_{12} \, {\cal F}_{12}(x) 
\end{equation}
satisfies  the following linear equation
\begin{equation}
  {\cal R}_{12}(x)\, B_2(x)=B_2(x)\, K_{12}^2\, {\cal R}_{21}(x)^{-1},
\label{Rlinear}
\end{equation} 
where  $B(x)$  is defined by
Eq.~\eqref{Bx}.
}

Note that, since ${\cal M}$  is almost $\nu$-rational and  $ \mathcal{F}(x)$ is $\mathfrak{h}$-invariant, $ \mathcal{F}(x)$ is $\nu$-rational.

\Problem{Strong Quantum Dynamical coBoundary Problem (SQDBP)\label{SP}}{
Find an almost $\nu$-rational invertible element $M(x)\in ({\mathbb C}(\tilde{\nu_1},\ldots,\tilde{\nu_n})\otimes  U_q({\mathfrak g}))^c$ solution of the Weak Quantum Dynamical coBoundary Problem such that 
${\cal F}(x)$ is equal to the standard solution $F(x)$ of the QDCE.
}

\begin{rem} 
It would have been better to denote $M(\tilde{\nu})$ such an element, 
but for notational reasons we prefer to call it $M(x)$. This should cause no confusion.
\end{rem}

%%%%%%%%%%%%%%%%%%%%%%%%%%%%%%%%%%%%%%%%%%%%%%%%%%%%%%%%%%%%%%%%%%%%%%%%%%%%%%%
\section{Solving the Quantum Dynamical coBoundary problem}
\label{sectionDBP}
%%%%%%%%%%%%%%%%%%%%%%%%%%%%%%%%%%%%%%%%%%%%%%%%%%%%%%%%%%%%%%%%%%%%%%%%%%%%%%%

In this section we assume that   ${\mathfrak g}=sl(n+1)$, and
$J$ will denote the Cremmer-Gervais cocycle defined in 
Theorem~\ref{cocycleBCG}. $F(x)$ and $R(x)$
will respectively denote the standard solutions of zero weight 
(\ref{Fprod}) and (\ref{R=FRF}) of the equations
(\ref{lineareq},\ref{eq:s-cocycle}) and  (\ref{sYB}) respectively.

We present here some general results concerning the weak and strong QDBP. 
Our aim is to identify and construct  elementary objects 
obeying simple algebraic rules  which will be the building blocks of the solutions 
of these problems.

We first propose, in Section \ref{sec-prim-loop}, a procedure to solve 
the WQDBP using the notion of {\em  primitive loop}. 
The study of these primitive loops does not however
enable us to solve the SQDBP.
Therefore, in Section \ref{sec-GaussQDC}, we introduce a different approach, 
based on the Gauss decomposition of $\mathcal{M}(x)$,
which leads to the solution of the SQDBP. 

%%%%%%%%%%%%%%%%%%%%%%%%%%%%%%%%%%%%%%%%%%%%%%%%%%%%%%%%%%%%%%%%%%%%%%%%%%%%%%%
\subsection{Primitive loops}\label{sec-prim-loop}

In order to study the properties of the solutions of the WQDBP, let us 
first introduce the notion of {\em primitive loop}, defined as follows:

\Definition{Primitive loop }{
For  any solution ${\cal M}(x)$ of the WQDBP, we define an element ${\cal P}(x) \in ({\mathbb C}(\nu_1,\ldots,\nu_n)\otimes  U_q({\mathfrak g}))^c $ by:
\begin{equation}
 {\cal P}(x)=v\,{\cal M}(x)^{-1}B(x)\,{\cal M}(x).\label{PdeXsl2}
\end{equation}
${\cal P}$ will be called the {\em primitive loop} associated to ${\cal M}(x).$
}

Note that, because ${\cal M}(x)$ is almost $\nu$-rational, ${\cal P}(x)$ is $\nu$-rational.

%Quite surprisingly this {\em primitive loop} is related to reflection algebras described in the last section, indeed we have the following proposition

A primitive loop satisfies various  properties and   verifies a {\em reflection equation} 
(see \eqref{echangeP} below) which is 
related to the notion of {\em reflection algebra} that will be introduced in
Section~\ref{sectionQRA}. More precisely, 

\Proposition{Properties of the primitive loop}{
Let ${\cal P}(x)$ be a primitive loop associated to a solution ${\cal M}(x)$ 
of the WQDBP. We have the following relations:
\begin{align}
  & R^J_{12}\,{\cal P}_{2}(x)\,R^J_{21}={\cal M}_1(x)^{-1}
         {\cal P}_{2}(xq^{h_1})\,{\cal M}_1(x), \label{linearPdeX}\\
  & \Delta^J({\cal P}(x))=(R^J_{12})^{-1}{\cal P}_1(x)\,R^J_{12}
            \,{\cal P}_{2}(x). \label{deltaPdeX}
\end{align}
As a consequence, ${\cal P}(x)$  satisfy the reflection equation:
\begin{align}
 &R^J_{21}\,{\cal P}_1(x)\,R^J_{12}\,{\cal P}_{2}(x)
  ={\cal P}_{2}(x)\,R^J_{21}\,{\cal P}_1(x)\,R^J_{12}.
 \label{echangeP}
\end{align}
}

\negvspace

\Proof{Eq.\eqref{linearPdeX} follows immediatly from the linear 
equation (\ref{Rlinear}) and from the relation (\ref{RMM=MMR}), 
itself inherited from \eqref{eq1-pb1}-\eqref{eq2-pb1} using 
\eqref{quasitriangularity2}.
Then, using successively (\ref{PdeXsl2}), (\ref{eq1-pb1}) and 
\eqref{propertyv}, the zero-weight property and \eqref{DeltaB}, and 
\eqref{linearPdeX}, 
we have:
\begin{align*}
  \Delta^J({\cal P}(x))
  &=J^{-1}\,\Delta(v)\,\Delta({\cal M}(x)^{-1})\,\Delta(B(x))\,
            \Delta({\cal M}(x))\,J\\
  &=(R^J_{21}R^J_{12})^{-1} v_1v_2\;{\cal M}_2(x)^{-1}
    {\cal M}_1(xq^{h_2})^{-1}{\cal F}_{12}(x)^{-1}\nonumber\\
  &\hspace{5.8cm}\times\Delta(B(x))\;{\cal F}_{12}(x)\,{\cal M}_1(xq^{h_2})\,
           {\cal M}_2(x)\\
  &=(R^J_{21}R^J_{12})^{-1} v_1v_2\;{\cal M}_2(x)^{-1}
    {\cal M}_1(xq^{h_2})^{-1}B_1(xq^{h_2})\,B_2(x)\,
    {\cal M}_1(xq^{h_2})\,{\cal M}_2(x)\\
  &=(R^J_{21}R^J_{12})^{-1} {\cal M}_2(x)^{-1}{\cal P}_1(xq^{h_2})\,
    {\cal M}_2(x)\,{\cal P}_2(x)\\
  &=(R^J_{21}R^J_{12})^{-1}R^J_{21}\,{\cal P}_{1}(x)\,R^J_{12}\, 
    {\cal P}_2(x)\\
  &=(R^J_{12})^{-1}{\cal P}_1(x)\,R^J_{12}\,{\cal P}_{2}(x), 
\end{align*}
which concludes the proof of \eqref{deltaPdeX}.
The relation (\ref{echangeP}) 
is a direct consequence of (\ref{deltaPdeX}) 
using (\ref{quasitriangularity2}).}

\begin{rem}
If ${\cal M}(x)$ is a solution of the WQDBP and if  ${\cal P}(x)$ is the associated primitive loop, then, in the fundamental representation of $U_q(sl(n+1))$, 
one has
 $\mathrm{tr}(\stackrel{f}{\pi}({\cal P}(x))=
\mathrm{tr}(\stackrel{f}{\pi}(vB(x))=q^{-n}(\nu_1+\cdots+\nu_{n+1}).$
\end{rem}

\begin{rem}
If $M(x)$ is a solution of the SQDBP such that its expression (\ref{Minv-fund})  in the fundamental 
representation of $U_q(sl(n+1))$ is given by ${\mathbf M}(x)$, 
the associated primitive  loop ${ P}(x)$  can also be computed in the fundamental representation and its  explicit expression is given by
\begin{equation} 
  {\bf P}(x)=
   \stackrel{f}{\pi}({ P}(x))=D\; q^{-n}\Big\{ \sum_{j=1}^n E_{j,j+1}
   +\sum_{k=0}^{n}(-1)^{n-k}{\cal S}_{n+1-k}(x)\; E_{n+1,k+1}\Big\} \;D^{-1},
   \label{Pinfund}
\end{equation}
where, for $1\leq m\leq n+1$, ${\cal S}_{m}(x)$ denotes the symmetric 
polynomial in $\nu_1,\ldots,\nu_{n+1}$  defined as
${\cal S}_{m}(x)=\sum_{1\leq i_1 < \cdots <i_m \leq n+1}\prod_{k=1}^{m}\nu_{i_k}.$
Note that ${\cal S}_{n+1}(x)=\nu_1\ldots\nu_{n+1}=1.$
\end{rem}

\bigskip

We now prove a sufficient condition for a given ${\cal M}$ to be a solution of the  WQDBP.

\Proposition{\label{theoreme}}{
Let ${\cal M}$ be an almost $\nu$-rational element of $(\mathbb{C}(\tilde{\nu}_1,\ldots,\tilde{\nu}_n)\otimes U_q(\mathfrak g))^c $, and let us define 
${\cal P}(x)\in (\mathbb{C}(\nu_1,\ldots,\nu_n)\otimes U_q(\mathfrak g))^c $  as 
\begin{align}
   &{\cal P}(x)=v{\cal M}(x)^{-1}B(x)\,{\cal M}(x).
   \label{diagonalisation}
\end{align}
If ${\cal P}$ satisfies the property 
\begin{align}
%   & \Delta^J({\cal P}(x))=(R^J_{12})^{-1}{\cal P}_1(x)\,R^J_{12}\,
%               {\cal P}_{2}(x), 
%\label{deltaPdeX2}\\
   & R^J_{12}\,{\cal P}_{2}(x)\,R^J_{21}={\cal M}_1(x)^{-1}
   {\cal P}_{2}(xq^{h_1})\,{\cal M}_1(x), 
\label{linearPdeX2}
\end{align}
then ${\cal F}(x)=\Delta({\cal M}(x))\,J\,{\cal M}_2(x)^{-1}
{\cal M}_1(xq^{h_2})^{-1}$ 
is of weight zero and is $\nu$-rational. \\
As a result,  ${\cal F}$ is a solution of 
the QDCE (\ref{eq:s-cocycle}) and 
${\cal R}(x)={\cal F}_{21}(x)^{-1}R_{12}\,{\cal F}_{12}(x)$ 
obeys the QDYBE (\ref{sYB}).
${\cal R}(x)$ satisfies the linear equation  (\ref{Rlinear}) and therefore  
 ${\cal M}(x)$ is a solution of the WQDBP.
}

\proof
The zero-weight property of ${\cal F}(x)$ is obtained from the 
following computation, using \eqref{diagonalisation}, \eqref{linearPdeX2}
and the quasitriangularity of $( U_q(\mathfrak{g}),\Delta^J,R^J)$:
\begin{align*}
  {\cal P}_{3}(xq^{h_1+h_2}&) \left\{ 
     \Delta_{12}({\cal M}(x))\,J_{12}\,{\cal M}_2(x)^{-1}
     {\cal M}_1(xq^{h_2})^{-1} \right\}
    \\
  &=J_{12}\,
   (\Delta^J \otimes id)({\cal M}_1(x)\,R^J_{12}\,{\cal P}_{2}(x)\,R^J_{21})\;
   {\cal M}_2(x)^{-1}{\cal M}_1(xq^{h_2})^{-1}
    \\
  &=\Delta_{12}({\cal M}(x))\,J_{12}\,R^J_{13}\,R^J_{23}\,{\cal P}_{3}(x)\,
    R^J_{32}\,R^J_{31}\,{\cal M}_2(x)^{-1}{\cal M}_1(xq^{h_2})^{-1}
    \\
  &=\Delta_{12}({\cal M}(x))\,J_{12}\,R^J_{13}\, {\cal M}_2(x)^{-1}
    {\cal P}_{3}(xq^{h_2})\,R^J_{31}\,{\cal M}_1(xq^{h_2})^{-1}
    \\
  &=\left\{\Delta_{12}({\cal M}(x))\,J_{12}\, {\cal M}_2(x)^{-1}
    {\cal M}_1(xq^{h_2})^{-1}\right\}  {\cal P}_{3}(xq^{h_2+h_1}). 
\end{align*}
By representing the third space on the fundamental representation and by tracing on it, we obtain that ${\cal F}$ commutes with ${\cal S}_1(xq^{h_1+h_2}).$
We now prove that this condition implies that   ${\cal F}$ is of zero weight.
Let $V,W$ be two finite dimensional $U_q(\mathfrak g)$-modules, 
${\cal F}_{VW}\in {\mathbb C}(\tilde{\nu}_1,\ldots,\tilde{\nu}_n)\otimes \End(V\otimes W).$
We can decompose
 $V\otimes W=\oplus_{\lambda\in \mathfrak{h}^*}(V\otimes W)[\lambda].$
Let us denote by $P_{\lambda}$ the projection on  $(V\otimes W)[\lambda]$  and 
let ${\cal F}_{\lambda,\mu}=P_{\lambda} {\cal F}_{VW}P_{\mu}.$
The previous condition implies that ${\cal F}_{\lambda,\mu}(a_\lambda-a_\mu)=0$, where 
\begin{equation}
a_{\lambda}={\cal S}_1(xq^{h(\lambda)})=
\sum_{i=1}^{n+1}\nu_i(x) q^{2(\zeta^{(i)}-\zeta^{(i-1)})(\lambda)}.
\end{equation}
When $\lambda\not=\mu,\ a_\lambda-a_\mu\not=0,$ and therefore, since ${\cal F}_{\lambda,\mu}\in {\mathbb C}(\tilde{\nu}_1,\ldots,\tilde{\nu}_n)\otimes \End(V\otimes W),$ we obtain that ${\cal F}_{\lambda,\mu}=0$ which is equivalent to the zero-weight condition.

As a result,
${\cal F}$ is a solution of 
the QDCE (\ref{eq:s-cocycle}) due to Theorem \ref{coboundcocycl},
and thus $\mathcal{R}(x)$ satisfies the QDYBE \eqref{sYB} due to 
Theorem \ref{TheoremQDCE}. It  obviously verifies 
\eqref{RMM=MMR}, and
the fact that it satisfies (\ref{Rlinear}) follows immediately from   
(\ref{diagonalisation},\ref{linearPdeX2}).
\qed

\begin{rem} 
The previous propositions show that the primitive loop element is of fundamental importance for solving the WQDBP and that the analog of the ABRR equation for ${\cal M}(x)$ is Eq.~(\ref{linearPdeX}), a (major) difference being that ${\cal P}(x)$ is also not known.
\end{rem}

\begin{rem} 
The previous proposition leads naturally to an algorithmic approach to the 
computation of the universal  $M(x).$ 
Indeed, the explicit value of ${\bf M}(x)$ being known in the fundamental 
representation, we can view the system of $(n+1)^2$ equations
\begin{equation}
   \big(\stackrel{f}{\pi} \otimes id\big)
   \big( R^J_{12}{\cal P}_{2}(x)R^J_{21}\big)
    ={\bf M}_1(x)^{-1}{\cal P}_{2}(xq^{h_1}){\bf M}_1(x),
\end{equation}
with moreover $\stackrel{f}{\pi}({\cal P})$ given by 
${\bf P}(x)=v_f{\bf M}(x)^{-1} {\bf B}(x)\,{\bf M}(x),$
as a system of universal equations fixing ${\cal P}(x)$ up to a 
central element which can  be
determined using the relation
\begin{equation}
   \Delta^J({\cal P}(x))=(R^J_{12})^{-1}{\cal P}_1(x)\,R^J_{12}\,
               {\cal P}_{2}(x). 
\label{deltaPdeX2}
\end{equation}
The universal expression of ${\cal P}(x)$ being known, the equation 
\begin{equation}
  {\cal M}(x){\cal P}(x)=vB(x){\cal M}(x)
\end{equation} 
is then a universal linear relation 
fixing ${\cal M}(x)$ up to a left-multiplication by an element of 
$({\mathbb C}(\nu_1,\ldots,\nu_n)\otimes {C}_q({\mathfrak h}))^c.$ 
In order to obtain a generic solution of the WQDBP,
it is sufficient  to show that this expression verifies universally 
(\ref{linearPdeX2}).
This is the path that we had initially followed in order to obtain explicit universal expressions for  ${\cal P}(x)$ and    ${\cal M}(x)$ 
such as those presented in the next section.

However, although the properties of the primitive loop ${\cal P}(x)$  
associated to ${\cal M}(x)$  are powerful tools to solve very explicitely 
the WQDBP, these  relations are not sufficient to ensure that 
${\cal M}(x)$ is a solution of the SQDBP. 
More precisely, let  ${\cal M}^{(1)}(x)$  be a solution of 
(\ref{diagonalisation}), (\ref{linearPdeX2}) and (\ref{deltaPdeX2}), and let
$u(x)\in ({\mathbb C}(\nu_1,\ldots,\nu_n)\otimes {C}_q({\mathfrak h}))^c$, 
we define  ${\cal M}^{(2)}(x)=u(x){\cal M}^{(1)}(x).$ 
Let us denote ${\cal P}^{(1)}(x),\ {\cal P}^{(2)}(x)$ the primitive loops 
corresponding respectively to ${\cal M}^{(1)}(x),\ {\cal M}^{(2)}(x)$. 
We have ${\cal P}^{(1)}(x)={\cal P}^{(2)}(x),$ and 
${\cal M}^{(2)}(x)$ is also a solution of (\ref{diagonalisation}),
(\ref{linearPdeX2}) and (\ref{deltaPdeX2}). 
Nevertheless, in general,  the corresponding 
${\cal F}^{(1)}(x)$ and ${\cal F}^{(2)}(x)$ are different. This shows that  
solutions of the WQDBP are in general  not solutions of 
the SQDBP.

In the next section we will solve this problem and 
 obtain sufficient conditions on 
${\cal M}(x)$  solution of the WQDBP to ensure that it is  also a solution to  
the SQDBP.
\end{rem}
 
%%%%%%%%%%%%%%%%%%%%%%%%%%%%%%%%%%%%%%%%%%%%%%%%%%%%%%%%%%%%%%%%%%%%%%%%%%%%%% 
\subsection{Gauss Decomposition of Quantum Dynamical coBoundary}
\label{sec-GaussQDC}

We propose here a new approach to construct the solutions of the SQDBP, 
based on the study of some fundamental building blocks  entering the Gauss 
decomposition of $M(x)$. We will prove in this section the
following theorem:

\Theorem{\label{maintheorem}}{ 
Let ${ \cal M}^{(0)}\in 
({\mathbb C}(\tilde{\nu}_1,\ldots,\tilde{\nu}_n)\otimes  U_q({\mathfrak h}))^c$ and ${\mathfrak C}^{[\pm]}\in 1\oplus \left({\mathbb C}(\nu_1,\ldots,\nu_n)
\otimes U^{\pm}_q({\mathfrak g})\right)^c$.
Let  ${\cal  M}^{(\pm)}\in 1\oplus \left({\mathbb C}(\nu_1,\ldots,\nu_n)
\otimes U^{\pm}_q({\mathfrak g})\right)^c$ 
be given by
\begin{align}
  {\cal M}^{(\pm)}(x)=\prod_{k=1}^{+\infty} {\mathfrak C}^{[\pm k]}(x)^{\pm 1},
  \quad \text{with}\quad
  &{\mathfrak C}^{[+ k]}(x)= \tau^{k-1}\left({\mathfrak C}^{[+]}(x)\right),
          \nonumber\\
  &{\mathfrak C}^{[- k]}(x)=B(x)^{-k}\,{\mathfrak C}^{[-]}(x)\, B(x)^{k} .
   \label{MdeXprod2}
\end{align}
Because $\tau$ is nilpotent, the  product defining  ${\cal  M}^{(+)}(x)$ is finite.
We define  ${ \cal M}\in \left({\mathbb C}(\tilde{\nu}_1,\ldots,\tilde{\nu}_n)\otimes U_q({\mathfrak g})\right)^c$ as 
\begin{equation}\label{MdeXprod}
  { \cal M}(x)= { \cal M}^{(0)}(x)\, {\cal M}^{(-)}(x)^{-1} 
               { \cal M}^{(+)}(x).
\end{equation}
The following algebraic relations on ${\cal M}^{(0)}$ and 
${\mathfrak C}^{[\pm]} $ are sufficient conditions to ensure that 
 ${\cal M}(x)$ is a solution of the  SQDBP:
\begin{align}
   &\Delta({\cal M}^{(0)}(x))\,S^{[1]}_{12}\,{\cal M}^{(0)}_2(x)^{-1} 
     {\cal M}^{(0)}_1(xq^{h_2})^{-1}=1 \otimes 1,
    \label{axiomABRR0}\\
%   &\Delta({\mathfrak C}^{[\pm]}(x))=V^{[\pm]}_{12}(x)\,K^{\mp 1}_{12}\,
%     V^{[\pm]}_{21}(x)\,K_{12}^{\pm 1},\ \text{with}\
%     V^{[\pm]}_{12}(x)=S^{[1]}_{21}\,K_{12}\, {\mathfrak C}^{[\pm]}_{1}(x)\, 
%       K_{12}^{-1}(S^{[1]}_{21})^{-1},
   &K_{12}^{-1}\,\Delta({\mathfrak C}^{[\pm]}(x))\,K_{12}=
    \big\{S^{[1]}_{21}\,{\mathfrak C}^{[\pm]}_{1}(x)\, (S^{[1]}_{21})^{-1}\big\}\,
    K^{\mp 1}_{12}\,
    \big\{S^{[1]}_{12}\,{\mathfrak C}^{[\pm]}_{2}(x)\, (S^{[1]}_{12})^{-1}\big\}\,
    K^{\pm 1}_{12},
     \label{axiomABRR1d}\\
   &{\mathfrak C}^{[\pm]}_{1}(xq^{h_2})=
     \big\{(S^{[1]}_{12})^{-1}S^{[1]}_{21}\,K_{12}\big\}\,
     {\mathfrak C}^{[\pm]}_{1}(x)\,
     \big\{K_{12}^{-1}(S^{[1]}_{21})^{-1}S^{[1]}_{12}\big\},
     \label{axiomABRR1}\\
   &{\mathfrak C}^{[-]}_2(x)\, {\mathfrak C}^{[+]}_1(xq^{h_2}) \,
     \big\{  B_2(x)\, (S^{[2]}_{12})^{-1}\widehat{J}{}^{[1]}_{12}\,
             S^{[2]}_{12} B_2(x)^{-1}\big\}\nonumber\\ 
   &\hspace{6.5cm}=\big\{ (S^{[1]}_{12})^{-1} \widehat{R}_{12}\, 
                                              S^{[1]}_{12}\big\}\,  
     {\mathfrak C}^{[+]}_1(xq^{h_2})\, {\mathfrak C}^{[-]}_2(x).
     \label{axiomABRR2prime}
\end{align}
}

\begin{rem}
In the next section we will find explicit solutions to these
 sufficient algebraic equations.
\end{rem}

The proof of this theorem decomposes in three lemmas.
The first lemma contains elementary results on $\mathcal{M}^{(+)}$ 
and $\mathcal{M}^{(-)}$:

\Lemma{\label{firstresults}}{
The infinite products (\ref{MdeXprod2}) define elements  $\mathcal{M}^{(+)}$ and $\mathcal{M}^{(-)}$ belonging to 
$1\oplus ({\mathbb C}(\nu_1,\ldots,\nu_n)\otimes U^{\pm}_q(\mathfrak{g}))^c.$
If Equations  \eqref{axiomABRR1d}, \eqref{axiomABRR1} and
\eqref{axiomABRR2prime} are  satisfied, the element
 $U(x)\in ({\mathbb C}(\nu_1,\ldots,\nu_n)\otimes  
U_q({\mathfrak g})^{\otimes 2})^c$ defined as 
\begin{equation}
U(x)= \Delta({\cal M}^{(+)}(x)) J {\cal M}^{(+)}_2(x)^{-1},
\end{equation}
can also be written as
\begin{equation}
   U_{12}(x)=S^{[1]}_{12}
      \prod_{k=1}^{n}\left( 
      {\mathfrak C}^{[+k]}_1(xq^{h_2})\, (S^{[k+1]}_{12})^{-1}
      \widehat{J}{}^{[k]}_{12}\,S^{[k+1]}_{12} \right) ,\label{Uprod}
\end{equation}
and therefore belongs to $({\mathbb C}(\nu_1,\ldots,\nu_n)\otimes U_q(\mathfrak{b}^+)\otimes U_q(\mathfrak{b}^-))^c.$
It satisfies the properties
\begin{align}
   & 
    (id \otimes \iota_{-})(U_{12}(x))
                  = S^{[1]}_{12}\,{\cal M}_1^{(+)}(xq^{h_2}),
   \label{axiomABRR3}\\
   &{\mathfrak C}^{[-]}_2(x)\,  B_2(x)\, (S^{[1]}_{12})^{-1}U_{12}(x)
    =  (S^{[1]}_{12})^{-1}\widehat{R}_{12}\, U_{12}(x)\, 
     {\mathfrak C}^{[-]}_2(x)\, B_2(x) .
    \label{axiomABRR2bis}
\end{align}
Moreover, $\mathcal{M}^{(-)}$ satisfies the following relations:
\begin{align}
   &{\mathfrak C}^{[-]}(x)\, B(x)\, {\cal M}^{(-)}(x)
   =  {\cal M}^{(-)}(x)\, B(x),
   \label{BM-T}\\
   &B_2(x)\,K_{12}\,\Delta( {\cal M}^{(-)}(x)^{-1})
   =\Delta'( {\cal M}^{(-)}(x)^{-1})\,S_{12}^{[1]}\, K_{12}\, 
    {\mathfrak C}^{[-]}_2(x)\, B_2(x)\,(S_{12}^{[1]})^{-1},
   \label{BdeltaMV}\\
   &(id \otimes \iota_{-})(\Delta({\cal M}^{(-)}(x)^{-1}))
                  =S^{[1]}_{12}\,{\cal M}_1^{(-)}(xq^{h_2})^{-1}\,
        (S^{[1]}_{12})^{-1}.
   \label{axiomABRR4}
\end{align}
}

\negvspace

\proof
\eqref{Uprod} is shown in Appendix~\ref{sec-lemmas} 
(see Lemma \ref{appendlemma2}).
\eqref{axiomABRR3} can then be derived from (\ref{Uprod})  using the fact that 
$\widehat{J}{}^{[k]}$ belongs to 
$\left( 1\otimes 1+(U^+_q({\mathfrak g})\otimes U^-_q({\mathfrak g}))^c\right)$.
\eqref{axiomABRR2bis} 
can be proved as follows:
first, acting by $\tau^{k-1}$ on the first space of (\ref{axiomABRR2prime}), 
we obtain
\begin{multline*}
    {\mathfrak C}^{[-]}_2(x)\, {\mathfrak C}^{[+ k]}_1(xq^{h_2})\, 
    \big\{ B_2(x)\, (S^{[k+1]}_{12})^{-1}\widehat{J}{}^{[k]}_{12}\,
    S^{[k+1]}_{12}\, B_2(x)^{-1}\big\} \\
   = \big\{ (S^{[k]}_{12})^{-1} \widehat{J}{}^{[k-1]}_{12}\, S^{[k]}_{12}
          \big\}\, 
     {\mathfrak C}^{[+ k]}_1(xq^{h_2})\, {\mathfrak C}^{[-]}_2(x),
%   \label{axiomABRR2ter}
\end{multline*}
then, using (\ref{Uprod}), we conclude the proof of (\ref{axiomABRR2bis}) 
by recursion.
 
\eqref{BM-T}  follows directly from the definition of
${\cal M}^{(-)}(x)$, and
 \eqref{BdeltaMV} follows directly from the definition of
${\cal M}^{(-)}(x)$ and from condition \eqref{axiomABRR1d}.
Finally, from \eqref{axiomABRR1d}, \eqref{axiomABRR1} and \eqref{DeltaB}, 
we have
\begin{multline*}
  \Delta({\cal M}^{(-)}(x)^{-1})
  =\prod_{k=+\infty}^{1}\big\{ B_1(xq^{h_2})^{-k} B_2(x)^{-k}
   S^{[1]}_{12}\,{\mathfrak C}_1^{[-]}(xq^{h_2})\,
   K_{12}^2\,\\
   \times{\mathfrak C}_2^{[-]}(x)\,(S^{[1]}_{12})^{-1}
   K_{12}^{-2} B_2(x)^{k}\,B_1(xq^{h_2})^{k}\big\},
\end{multline*}
which directly implies \eqref{axiomABRR4}.		
\qed

\Lemma{\label{lemme-ABRR}}{ 
Let us suppose that ${\cal M}$ is defined  as in Theorem~\ref{maintheorem} and assume that the hypothesis of Theorem~\ref{maintheorem} are satisfied.
Then ${\cal F}(x)$ defined as
${\cal F}(x)=\Delta({\cal M}(x))\,J\,{\cal M}_2(x)^{-1}
{\cal M}_1(xq^{h_2})^{-1}$ belongs to 
$ 1 \otimes 1 + \left({\mathbb C}(\tilde{\nu}_1,\ldots,\tilde{\nu}_n)\otimes U_q({\mathfrak g})\otimes U^-_q({\mathfrak g})\right)^c$ and satisfies 
 the ABRR identity \eqref{lineareq}. 
}

\proof
Since, from Lemma~\ref{firstresults}, $U(x)$ belongs to   
$({\mathbb C}(\nu_1,\ldots,\nu_n)\otimes U_q(\mathfrak{g})\otimes U_q(\mathfrak{b}^-))^c$, 
we have ${\cal F}(x)\in  ({\mathbb C}(\nu_1,\ldots,\nu_n)\otimes U_q(\mathfrak{g})\otimes U_q(\mathfrak{b}^-))^c.$
We can therefore define $(id\otimes \iota_-)({\cal F}(x))$. 
{}From the relations \eqref{axiomABRR3}, \eqref{axiomABRR4}
and the identity \eqref{axiomABRR0},  we deduce that
$(id \otimes \iota_-)({\cal F}(x))=1 \otimes 1$,
which means that ${\cal F}(x)\in 
1 \otimes 1 + \big({\mathbb C}(\tilde{\nu}_1,\ldots,\tilde{\nu}_n)\otimes U_q({\mathfrak g})\otimes U^-_q({\mathfrak g})\big)^c.$

The fact that $\mathcal{F}$ satisfies the ABRR relation 
can be proved as follows. 
Using respectively (\ref{MdeXprod}), (\ref{BdeltaMV}), the quasitriangularity
property \eqref{quasitriangularity2}, (\ref{axiomABRR2bis}) and
(\ref{BM-T}), we have
\begin{align*}
   {\widehat{R}}{}^{-1}B_2(x)\,{\cal F}_{12}(x)
   &=R^{-1} K_{12}\, B_2(x)\, \Delta({\cal M}^{(0)}(x))\,
     \Delta( {\cal M}^{(-)}(x)^{-1})\, U_{12}(x)\nonumber\\
    &\hspace{3cm} \times
     {\cal M}_2^{(-)}(x)\, {\cal M}_2^{(0)}(x)^{-1}{\cal M}_1(xq^{h_2})^{-1}\\ 
   &= \Delta({\cal M}^{(0)}(x))\, R^{-1} \Delta'( {\cal M}^{(-)}(x)^{-1})\,
     S_{12}^{[1]}\, K_{12}\,\nonumber\\ 
    &\hspace{3cm} \times
    \big\{ {\mathfrak C}^{[-]}_2(x)\, B_2(x)\,(S_{12}^{[1]})^{-1} U_{12}(x)\big\}
    \nonumber\\ 
    &\hspace{3cm} \times
    {\cal M}_2^{(-)}(x)\, {\cal M}_2^{(0)}(x)^{-1}{\cal M}_1(xq^{h_2})^{-1}\\
   &= \Delta({\cal M}^{(0)}(x))\,\Delta( {\cal M}^{(-)}(x)^{-1})\,U_{12}(x)
                        \nonumber\\ 
    &\hspace{3cm} \times
    \big\{{\mathfrak C}^{[-]}_2(x)\, B_2(x)\, 
          {\cal M}_2^{(-)}(x)\big\} 
%	  \nonumber\\ 
%    &\hspace{3cm} \times
    {\cal M}_2^{(0)}(x)^{-1}{\cal M}_1(xq^{h_2})^{-1}\\
   &=\Delta({\cal M}(x))\,J{\cal M}_2(x)^{-1} B_2(x)\, {\cal M}_1(xq^{h_2})^{-1}
     \\
   &={\cal F}_{12}(x)\, B_2(x),
\end{align*}  
%
%Then, using (\ref{MdeXprod}), (\ref{BdeltaMV}), (\ref{BM-T}) and 
%(\ref{axiomABRR2}), we have
%
%\begin{align*}
%   &{\widehat{R}}{}^{-1}B_2(x)\,{\cal F}_{12}(x)\,B_2(x)^{-1}\\ 
%   &\quad=R^{-1}\Delta({\cal M}^{(0)}(x))\,\Delta'( {\cal M}^{(-)}(x)^{-1})\, 
%            V^{[-]}_{21}(x)\,B_2(x)\,K_{12}\,\Delta({\cal M}^{(+)}(x))\,J_{12}
%    \\
%   &\hspace{3.5cm}
%   \times {\cal M}^{(+)}_2(x)^{-1} B_2(x)^{-1} {\mathfrak C}^{[-]}_2(x)^{-1}
%          {\cal M}^{(+)}_2(x)\, {\cal M}_2(x)^{-1}{\cal M}_1(xq^{h_2})^{-1},
%     \\
%   &\quad=\Delta({\cal M}(x))\, \big\{
%          \Delta({\cal M}^{(+)}(x)^{-1})\, {{R}}{}^{-1} V^{[-]}_{21}(x)\,
%          B_2(x)\,K_{12}\,U_{12}(x)\,           
%     \\
%   &\hspace{3.5cm}
%   \times B_2(x)^{-1}{\mathfrak C}^{[-]}_2(x)^{-1} {\cal M}^{(+)}_2(x) \big\}\,
%     {\cal M}_2(x)^{-1}{\cal M}_1(xq^{h_2})^{-1},
%     \\
%   &\quad=\Delta({\cal M}(x))\,J{\cal M}_2(x)^{-1}{\cal M}_1(xq^{h_2})^{-1},
%     \\
%   &\quad={\cal F}_{12}(x).
%\end{align*}
%
which ends the proof of the ABRR identity for ${\cal F}(x).$
\qed 

Finally, we have the following uniqueness lemma: 

\Lemma{\label{uniqueness}}{ 
Let ${\cal F}\in\left( 1 \otimes 1 
+({\mathbb C}(\tilde{\nu}_1,\ldots,\tilde{\nu}_n)\otimes U_q({\mathfrak g})\otimes U^-_q({\mathfrak g}))^c\right)$ and assume that  ${\cal F}(x)$ is a solution of the ABRR equation, then  ${\cal F}(x)$ is equal to the standard 
solution $F(x)$
 \eqref{Fprod} of the QDCE. 
}

\proof
${\cal F}(x)$ and $F(x)$ being  both in 
$\left(1 \otimes 1 + ({\mathbb C}(\tilde{\nu}_1,\ldots,\tilde{\nu}_n)\otimes U_q({\mathfrak g})\otimes U^-_q({\mathfrak g}))^c\right)$, 
we define $Y=F^{-1}(x){\cal F}(x)-1\otimes 1$. 
Then $Y \in ({\mathbb C}(\tilde{\nu}_1,\ldots,\tilde{\nu}_n)\otimes U_q({\mathfrak g})\otimes U^-_q({\mathfrak g}))^c.$
 ${\cal F}(x)$ and $F(x)$ being  both solutions of the ABRR identity, 
we also have $[Y(x),B_2(x)]=0.$ 

Let $V,W$ be finite dimensional $U_q(\mathfrak{g})$-modules. We can decompose 
$W=\oplus_{\lambda\in \mathfrak{h}^*}W[\lambda]$,
and consider $P_\lambda$ the associated projection on $W[\lambda]$. 
We define 
$Y_{\lambda,\mu}=(\id\otimes P_{\lambda})Y_{V,W}(\id\otimes P_{\mu}).$
Then, the fact that $[Y(x),B_2(x)]=0$ implies that 
$(b(\lambda)-b(\mu)) Y_{\lambda,\mu}=0$ with 
$B(x)_{\vert V[\lambda]}=b(\lambda)\id_{V[\lambda]}.$ 
Since  $Y_{V,W}$ is strictly lower  triangular on $W$, the only possible nonzero $Y_{\lambda,\mu}$ are associated to 
$\lambda<\mu.$
In this case the rational fraction $(b(\lambda)-b(\mu))\not=0$ and therefore 
$Y_{\lambda,\mu}=0.$ 

As a result, $Y=0.$
\qed

Theorem~\ref{maintheorem} is a direct consequence of Lemma~\ref{lemme-ABRR}
and Lemma~\ref{uniqueness}.

\bigskip

We now give a direct derivation of a weaker result but which proof is interesting in itself, namely that, under the hypothesis of Theorem~\ref{maintheorem},  $\mathcal{M}(x)$ is a solution of the WQDBP.

\Lemma{\label{lemme-Mweak}}{
If ${\cal M}(x)$ is defined as in 
Theorem~\ref{maintheorem},
then ${\cal F}(x)=\Delta({\cal M}(x))\,J\,{\cal M}_2(x)^{-1}
{\cal M}_1(xq^{h_2})^{-1}$ is a $\nu$-rational  zero-weight solution of 
the QDCE.
As a result, ${\cal M}(x)$ is also a solution of the WQDBP.
}

\proof
Let us consider ${\cal P}(x)=v {\cal M}(x)^{-1} B(x)\,{\cal M}(x)$ and show
that the hypothesis of Proposition~\ref{theoreme} are satisfied.

First, from (\ref{MdeXprod}) and (\ref{BM-T}), we have
\begin{align}
  {\cal P}(x)=v{\cal M}^{(+)}(x)^{-1} {\mathfrak C}^{[-]}(x)\, B(x)\, 
           {\cal M}^{(+)}(x).\label{PM+M-}
\end{align}
%
%Hence, we have 
%\begin{align*}
%&J(R^{J}_{12})^{-1}{\cal P}_1(x)R^{J}{\cal %P}_2(x)J^{-1}=v_1v_2(R_{21}R_{12})^{-1}{\cal X}_{12}(x){\cal X}_{21}(x)
%\end{align*}
%with
%\begin{align*}
%{\cal X}_{12}(x)&=R_{21} J_{21}{\cal M}_1^{(+)}(x)^{-1}{\mathfrak %C}_1^{[-]}(x)B_1(x){\cal M}_1^{(+)}(x)J^{-1}_{21}\\
%&=\Delta({\cal M}^{(+)}(x)^{-1}) R_{21}\left(  U_{21}(x) \right) 
%\left({\mathfrak C}^{[-]}_1(x)B_1(x)\right)  \left(
%U_{21}(x)^{-1}
%\right)\Delta'({\cal M}^{(+)}(x))\\
%&=\Delta({\cal M}^{(+)}(x)^{-1}) V^{[-]}_{12}(x)KB_1(x)
%\Delta'({\cal M}^{(+)}(x)).
%\end{align*}
%These relations are obtained  by 
%using (\ref{axiomABRR3})(\ref{axiomABRR2}). We then conclude the proof of (\ref{deltaPdeX}) using (\ref{axiomABRR1d}) and properties of $v$ recalled in the appendix.

 In order to show that $\mathcal{P}$ satisfies the linear 
 equation (\ref{linearPdeX2}), 
 we consider the quantity
\begin{equation*} 
  {\cal X}_{12}(x)={\cal M}^{(+)}_1(x)\, R^J_{12}\,{\cal P}_{2}(x)\,
  R^J_{21} \, {\cal M}^{(+)}_1(x)^{-1}
\end{equation*}
 and observe, using successively (\ref{PM+M-}), the definition of $U(x)$
 and Eq.~(\ref{axiomABRR2bis}) of Lemma~\ref{firstresults}, that 
 \begin{align*}
 {\cal X}_{12}(x)
   &=v_2 \,U_{21}(x)^{-1} R_{12}\, U_{12}(x)\,
     {\mathfrak C}_2^{[-]}(x)\, B_2(x)\, U_{12}(x)^{-1} R_{21}\, U_{21}(x)\\
   &=v_2\, U_{21}(x)^{-1} K_{12}\, S^{[1]}_{12}\,\mathfrak{C}_2^{[-]}(x)\,B_2(x)\,
     (S^{[1]}_{12})^{-1} R_{21}\,U_{21}(x).
\end{align*}
Note that, from \eqref{axiomABRR3}, 
${\cal X}_{12}(x)\in (U_q({\mathfrak b}^-)\otimes U_q({\mathfrak g}))^c$.
Then, using successively (\ref{axiomABRR2bis}), (\ref{axiomABRR1d}) and
\eqref{DeltaB}, the quasitriangularity property (\ref{quasitriangularity2}),
(\ref{axiomABRR1d}) and \eqref{DeltaB} again, and finally
(\ref{axiomABRR2bis}), we have
\begin{align*}
  \big\{ {\mathfrak C}^{[-]}_1(x)\, B_1(x) \big\}\,  {\cal X}_{12}(x)
  &=v_2 \, \big\{ {\mathfrak C}^{[-]}_1(x)\, B_1(x)\, U_{21}(x)^{-1} \big\}
  \\
  &\hspace{2cm}\times 
      K_{12}\,S^{[1]}_{12}\,\mathfrak{C}_2^{[-]}(x)\,B_2(x) \,  
     (S^{[1]}_{12})^{-1} R_{21}\,U_{21}(x)\\
  &=v_2 \, U_{21}(x)^{-1} \widehat{R}^{-1}_{21} \\
  &\hspace{2cm}\times\big\{
    S^{[1]}_{21}\,\mathfrak{C}_1^{[-]}(x)\,(S^{[1]}_{21})^{-1} K_{12}\,
    S^{[1]}_{12}\,\mathfrak{C}_2^{[-]}(x)\,(S^{[1]}_{12})^{-1} K_{12}^{-1}
    \big\}\\
  &\hspace{2cm}\times  \big\{B_1(x)\, B_2(x)\, K^2\}\,\widehat{R}_{21}\,U_{21}(x)
  \\
  &=v_2\, U_{21}(x)^{-1} R^{-1}_{21} \Delta({\mathfrak C}^{[-]}(x))\, 
    \Delta(B(x))\, R_{21}\, U_{21}(x)
    \\
  &=v_2\, U_{21}(x)^{-1} \Delta'({\mathfrak C}^{[-]}(x))\, \Delta(B(x))\,U_{21}(x)
  \\
  &=v_2\, U_{21}(x)^{-1}\, K_{12}\,S^{[1]}_{12}\,\mathfrak{C}_2^{[-]}(x)\, 
  (S^{[1]}_{12})^{-1} B_2(x)\, K_{12}\,\\ 
  &\hspace{2cm}\times  \big\{
  S^{[1]}_{21}\,\mathfrak{C}_1^{[-]}(x)\,B_1(x)\,(S^{[1]}_{21})^{-1} 
   U_{21}(x)\big\}\\
 &=\big\{v_2\, U_{21}(x)^{-1} K_{12}\, S^{[1]}_{12}\,\mathfrak{C}_2^{[-]}(x)\,
     B_2(x)\,
     (S^{[1]}_{12})^{-1} R_{21}\,U_{21}(x)\big\}\,\\ 
  &\hspace{2cm}\times
     \big\{{\mathfrak C}^{[-]}_1(x)\, B_1(x)\big\}\\
 &={\cal X}_{12}(x)\, \big\{ {\mathfrak C}^{[-]}_1(x)\, B_1(x)\big\} .
 \end{align*}
 As a consequence, if we denote 
 \begin{align*}
   {\cal M}^{(-)\;[N]}(x)&=\prod_{k=1}^{N}
   \Big\{ B(x)^{-k}{\mathfrak C}^{[-]}(x)^{-1}B(x)^{k} \Big\}
 \end{align*}
 we obtain, using (\ref{axiomABRR1}), 
 \begin{align*}
  {\cal M}^{(-)\;[N]}_1(x)^{-1}  
%{\cal M}^{(+)}_1(x) R^J_{12}{\cal P}_{2}(x)R^J_{21}  {\cal M}^{(+)}_1(x)^{-1}
  {\cal X}_{12}(x)\,{\cal M}^{(-)\;[N]}_1(x)
  &=B_1(x)^{-N}{\cal X}_{12}(x)\,  B_1(x)^{N}\\
  &=v_2\, \big\{ B_1(x)^{-N} 
               U_{21}(x)\, B_1(x)^{N} \big\}^{-1} S^{[1]}_{21}\,
    {\mathfrak C}_2^{[-]}(xq^{h_1})\,\\
  &\hspace{1.8cm}\times 
    B_2(xq^{h_1})\,(S_{21}^{[1]})^{-1}
        \big\{ B_1(x)^{-N}\widehat{R}_{21}\, B_1(x)^{N} \big\}\,\\
  &\hspace{1.8cm}\times
      \big\{ B_1(x)^{-N} U_{21}(x)\,B_1(x)^{N} \big\}.
 \end{align*}
 Using the fact that 
\begin{align*}
  \forall \xi \in U_q({\mathfrak b}^-),\ 
  \text{lim}_{N \rightarrow +\infty}\big\{B(x)^{-N} \xi\, B(x)^{N}\big\}
  =\iota_-(\xi),
\end{align*}
(which is shown in each finite dimensional $U_q(\mathfrak g)$ module),
as well as the fact that $\widehat{R}_{12}\in 
\big(1\otimes1 +U^+_q({\mathfrak g})\otimes U^-_q({\mathfrak g})\big)$ 
and the property (\ref{axiomABRR3}), we obtain
\begin{align*}
 {\cal M}_1(x)\, R^J_{12}\,{\cal P}_{2}(x)\,R^J_{21}\,{\cal M}_1(x)^{-1}
  &=v_2\, {\cal M}^{(+)}_2(xq^{h_1})^{-1} {\mathfrak C}_2^{[-]}(xq^{h_1})\,
      B_2(xq^{h_1})\,{\cal M}^{(+)}_2(xq^{h_1})\\
  &={\cal P}_2(xq^{h_1}),
\end{align*}
which concludes the proof of (\ref{linearPdeX2}).

${\cal M}(x)$ being  almost $\nu$-rational, the hypothesis of Proposition~\ref{theoreme} are satisfied and  $\mathcal{F}$ is 
a $\nu$-rational zero-weight solution of the QDCE. Therefore $\mathcal{M}$ is
a solution of the WQDBP.
\qed

\begin{rem} Trivial Gauge transformations.
\par\noindent
It is important to remark that the solutions to the SQDBP and WQDBP we found 
along the previous construction are by no means unique.
Indeed, let $M(x)$ be a solution of the SQDBP associated to the standard 
solution $F(x)$ of the QDCE and to a given cocycle $J.$ 
Let $u$ be a  dynamical group like element  
$u \in ({\mathbb C}(\tilde{\nu}_1,\ldots,\tilde{\nu}_n)\otimes U_q({\mathfrak h}))^c$, i.e. verifying  $\Delta(u(x))=u_1(xq^{h_2})u_2(x),$ and let  
$y\in (U_q({\mathfrak g}))^c,$ then the element $u(x)M(x)y^{-1}$ is a solution of the QDBE 
associated to the standard solution $F(x)$ of the QDCE and to 
the cocycle $\Delta(y)Jy_2^{-1}y_1^{-1}.$
It is  an interesting problem, not adressed here,  to  find the entire set of solutions of the SQDBP up to these transformations.
\end{rem}

%%%%%%%%%%%%%%%%%%%%%%%%%%%%%%%%%%%%%%%%%%%%%%%%%%%%%%%%%%%%%%%%%%%%%%%%%%%%%%%
\section{Explicit construction of Quantum Dynamical Coboundaries for 
$U_q(sl(n+1))$}
\label{sectionConstr}
%%%%%%%%%%%%%%%%%%%%%%%%%%%%%%%%%%%%%%%%%%%%%%%%%%%%%%%%%%%%%%%%%%%%%%%%%%%%%%%

%%%%%%%%%%%%%%%%%%%%%%%%%%%%%%%%%%%%%%%%%%%%%%%%%%%%%%%%%%%%%%%%%%%%%%%%%%%%%%%
\subsection{The $U_q(sl(2))$ case}

In this case the following result holds, which gives a new derivation of the result of \cite{BR2}. 

\Theorem{\label{th-case1}}{
In the $U_q(sl(2))$ case ,  a  solution ${\cal M}(x)$ to the SQDBP is given by  (\ref{MdeXprod}) and 
(\ref{MdeXprod2})  with
\begin{align}
   &{\cal M}^{(0)}(x)=1,\qquad
   {\mathfrak C}^{[+]}(x)=e_{q^{-1}}^{- x e},\qquad
   {\mathfrak C}^{[-]}(x)^{-1}= e_{q^{-1}}^{(xq^{h+1})^{-1}f}.
  \label{MdeXsl2}
\end{align}
}

\negvspace

\proof
Equations (\ref{axiomABRR1d},\ref{axiomABRR1},\ref{axiomABRR2prime}) are obtained from the properties (\ref{expsomme},\ref{expproduit}) of the $q$-exponential.
\qed

%%%%%%%%%%%%%%%%%%%%%%%%%%%%%%%%%%%%%%%%%%%%%%%%%%%%%%%%%%%%%%%%%%%%%%%%%%%%%%%
\subsection{General solution of the SQDBP for $U_q(sl(n+1))$}

\Theorem{\label{theo-Mexpr}}{
A solution ${\cal M}(x)$ of the SQDBP  for $U_q(sl(n+1))$ ($n\geq 1$) 
is  given by the infinite product 
\begin{align}
  &{\cal M}(x)= {\cal M}^{(0)}(x)\ 
  \prod_{k=+\infty}^{1}\Big(B(x)^{-k}{\mathfrak C}^{[-]}(x)\,B(x)^{k}\Big)\
  \prod_{k=0}^{+\infty} \tau^{k}\big({\mathfrak C}^{[+]}(x)\big), 
\end{align}
with
%\footnote{We use here the notations of Remark~\ref{rem-notations}, i.e.
%$\zeta^{(0)}=0$.}
%
\begin{align}
  &{\cal M}^{(0)}(x)= \prod_{k=1}^{n} \nu_{k+1}^{\frac{1}{2}\zeta^{(k)}} 
                      q^{-\frac{1}{2}(\zeta^{(k)})^2},
     \nonumber\\
  &{\mathfrak C}^{[+]}(x)= %e_{q^{-1}}^{-\nu^{-1}_{2}e_{(1)}} 
    \prod_{k=1}^{n}e_{q^{-1}}^{-\nu^{-1}_{k+1}\,q^{\zeta^{(k-1)}}e_{(k)}},
     \nonumber\\
  &{\mathfrak C}^{[-]}(x)^{-1}=%e_{q^{-1}}^{\nu_{2}q^{-h_{(1)}-1}f_{(1)}} 
    \prod_{k=1}^{n}e_{q^{-1}}^{\nu_{k+1}\,q^{-\zeta^{(k-1)}-h_{(k)}-1}f_{(k)}}.
 \label{MdeXsln}
\end{align}
}

%\negvspace

\begin{rem}
One may wonder why the expression of ${\cal M}^{(0)}(x)$ that one obtains here in the $n=1$ case is not equal to $1$ as in Theorem~\ref{th-case1}. Actually,
since in this case $J=1,$ one can also choose the simpler solution ${\cal M}^{(0)}(x)=1$, which gives the result of the previous section. 
\end{rem}

We first begin by proving a lemma interesting in itself.

\Lemma{}{Let ${\mathfrak C}^{[\pm]}$ solutions to the equations 
(\ref{axiomABRR1d},\ref{axiomABRR1})  we define 
 \begin{align*}
   &{\cal W}_{12}(x)= {\mathfrak C}^{[+]}_1(xq^{h_2})\,
   \{B_{2}(x)\,(S^{[2]}_{12})^{-1}\widehat{J}^{[1]}_{12}\,
     S^{[2]}_{12}\,B_{2}(x)^{-1}\}\,
   {\mathfrak C}_2^{[-]}(x)^{-1}, \\
   &\widetilde{\cal W}_{12}(x)=
    {\mathfrak C}^{[-]}_2(x)^{-1}
   \{(S^{[1]}_{12})^{-1}\widehat{R}_{12}\,S^{[1]}_{12}\}\,
   {\mathfrak C}_1^{[+]}(xq^{h_2}).
\end{align*}
These elements satisfy the relations:
\begin{align*}
   &(id \otimes \Delta)({\cal W}_{12} (x))
        = K_{23}\, S^{[1]}_{32}\, {\cal W}_{13}(x q^{h_2})\,
	  K^{-1}_{23}\,(S^{[1]}_{32})^{-1}\,
	  \\
	  &\hspace{5cm}\times{\mathfrak C}^{[+]}_1(xq^{h_2+h_3})^{-1}
	  S^{[1]}_{23}\,
	  {\cal W}_{12}(x q^{h_3})\,(S^{[1]}_{23})^{-1},\\
%       = K_{32}\, S^{[1]}_{32}\,
%         {\mathfrak C}^{[+]}_1(xq^{h_2+h_3})\,\big\{B_{3}(xq^{h_2})\,
%         (S^{[2]}_{13})^{-1} \widehat{J}^{[1]}_{13} S^{[2]}_{13}\,
%         B_{3}(xq^{h_2})^{-1}\big\}\\
%   &\hspace{3cm}\times{\mathfrak C}_3^{[-]}(xq^{h_2})^{-1} 
%         \, K^{-1}_{32}\,(S^{[1]}_{32})^{-1}\,S^{[1]}_{23}\\
%   &\hspace{3cm}\times 
%         \big\{B_{2}(xq^{h_3})\,(S^{[2]}_{12})^{-1} \widehat{J}^{[1]}_{12}
%         S^{[2]}_{12}\,B_{2}(xq^{h_3})^{-1}\big\}\, 
%         {\mathfrak C}_2^{[-]}(xq^{h_3})^{-1} (S^{[1]}_{23})^{-1},\\
   &(id \otimes \Delta)(\widetilde{\cal W}_{12}(x))
       = K_{23}\, S^{[1]}_{32}\, \widetilde{\cal W}_{13}(x q^{h_2})\,
	  K^{-1}_{23}\,(S^{[1]}_{32})^{-1}\,
	  \\
	  &\hspace{5cm}\times{\mathfrak C}^{[+]}_1(xq^{h_2+h_3})^{-1}
	  S^{[1]}_{23}\,
	  \widetilde{\cal W}_{12}(x q^{h_3})\,(S^{[1]}_{23})^{-1},
%	  K_{32}\,S^{[1]}_{32}\,{\mathfrak C}_3^{[-]}(xq^{h_2})^{-1}
%         {\mathfrak C}_2^{[-]}(x)^{-1} (S^{[1]}_{13})^{-1} \widehat{R}_{13}  
%         S^{[1]}_{13} \, S^{[1]}_{23}\,\\
%   &\hspace{3cm}\times   K^{-1}_{23} (S^{[1]}_{32})^{-1}
%        ( S^{[1]}_{12})^{-1}
%        \widehat{R}_{12}\, S^{[1]}_{12}\, (S^{[1]}_{23})^{-1} 
%        {\mathfrak C}^{[+]}_1(xq^{h_2+h_3}).
 \end{align*}
as well as
\begin{align*}
  &(\Delta  \otimes id)({\cal W}_{12} (x))=
        K_{12}\,S^{[1]}_{21}\,{\cal W}_{13} (x)\, K^{-1}_{12} (S^{[1]}_{21})^{-1}
	{\mathfrak C}_3^{[-]}(x)\, S^{[1]}_{12}\,{\cal W}_{23} (x)\, 
	(S^{[1]}_{12})^{-1},\\
%     = K_{12}\,S^{[1]}_{21}\,\mathfrak{C}^{[+]}_1(xq^{h_3})
%       \big\{B_{3}(x)\,(S^{[2]}_{13})^{-1} \widehat{J}^{[1]}_{13}\,
%       S^{[2]}_{13}\, B_{3}(x)^{-1}\big\}\, K^{-1}_{12} (S^{[1]}_{21})^{-1} \\
%   &\hspace{3cm}\times S^{[1]}_{12}\,
%       {\mathfrak C}_2^{[+]}(xq^{h_3})\,
%       \big\{ B_{3}(x)\,(S^{[2]}_{23})^{-1} \widehat{J}^{[1]}_{23}\,
%       S^{[2]}_{23}\, B_{3}(x)^{-1}\big\}\, 
%       {\mathfrak C}_3^{[-]}(x)^{-1} (S^{[1]}_{12})^{-1},\\
  &( \Delta  \otimes id)(\widetilde{\cal W}_{12}(x))=
        K_{12}\,S^{[1]}_{21}\,\widetilde{\cal W}_{13} (x)\, 
	K^{-1}_{12} (S^{[1]}_{21})^{-1}
	{\mathfrak C}_3^{[-]}(x)\, S^{[1]}_{12}\,\widetilde{\cal W}_{23} (x)\, 
	(S^{[1]}_{12})^{-1},  
%      = K_{12}\,S^{[1]}_{21}\,{\mathfrak C}_3^{[-]}(x)^{-1} (S^{[1]}_{13})^{-1}
%        \widehat{R}_{13}\, S^{[1]}_{13}\,    
%       {\mathfrak C}_1^{[+]}(xq^{h_3})\, S^{[1]}_{12}\, \\
%   &\hspace{3cm}\times  K^{-1}_{12}
%       (S^{[1]}_{21})^{-1} (S^{[1]}_{23})^{-1}
%       \widehat{R}_{23}\, S^{[1]}_{23}\, 
%       {\mathfrak C}^{[+]}_2(xq^{h_3})\, (S^{[1]}_{12})^{-1}.
 \end{align*}
}

\proof
This is easy to show using only  only (\ref{axiomABRR1d}), (\ref{axiomABRR1}), the definition of the
cocycle in Theorem~\ref{cocycleBCG} and the quasitriangularity properties of
the $R$-matrix. 

We now give the proof of the theorem.

\proof
As $\mathcal{M}$ is of the form \eqref{MdeXprod}, \eqref{MdeXprod2}, with
${\mathfrak C}^{[\pm]} $ being $\nu$-rational, it is enough to show that
the hypothesis of Theorem~\ref{maintheorem} are verified, i.e. that
\eqref{axiomABRR0}, \eqref{axiomABRR1d}, \eqref{axiomABRR1} and 
\eqref{axiomABRR2prime} are satisfied.
 
Equations (\ref{axiomABRR0}) and (\ref{axiomABRR1}) are trivial to check 
using the elementary property
\begin{align}
   \nu_k(xq^{h})=\nu_k(x)\, 
         q^{2(\zeta^{(k)}-\zeta^{(k-1)})},
   \quad \forall k=1,\ldots,n+1.
\end{align}
Equation (\ref{axiomABRR1d}) can be easily deduced from property 
(\ref{expsomme}) and from the following commutation properties, 
verified for all $i,j$ satisfying  $2 \leq j+1 \leq i \leq n$,
\begin{align*}
  & [ q^{-\zeta^{(i-1)}-2h_{(i)}} \otimes q^{-\zeta^{(i-1)}-h_{(i)}}f_{(i)}
    \ ,\  
       q^{-\zeta^{(j-1)}-h_{(j)}}f_{(j)} \otimes 
       q^{-\zeta^{(j-1)}-h_{(j)}} ]=0,
   \\
  &[ q^{\zeta^{(i-1)}}e_{(i)} \otimes q^{\zeta^{(i-1)}+h_{(i)}}
   \ ,\  
   q^{\zeta^{(j-1)}} \otimes q^{\zeta^{(j-1)}}e_{(j)} ]
   =0.
\end{align*}

Equation (\ref{axiomABRR2prime}) is slightly more difficult to show. 

Although our proof is a bit unsatisfactory, we have chosen a  method preventing us 
 to enter too deeply in the combinatorics of  $q$-exponentials. 
Our proof consists in two steps: first, we show that this relation holds in 
the fundamental representation, and then that it can be obtained in other 
representations by fusion from the fundamental representation. 

The fact that \eqref{axiomABRR2prime} is satisfied in the fundamental 
representation, i.e. that
\begin{equation*}
   (\stackrel{f}{\pi}\otimes \stackrel{f}{\pi})
       ({\cal W}_{12}-\widetilde{\cal W}_{12})(x)=0,
\end{equation*} 
is proved by an explicit computation in  Appendix~\ref{sec-lemmas} 
(Lemma \ref{appendlemma3}). 

Let us now study the fusion properties of this relation. This is the consequence of the previous lemma and from it we obtain:
if  
$(\stackrel{\Lambda_i}{\pi},\stackrel{\Lambda_i}{V})$, $i=1,2,3$, are 
 representations of $U_q({\mathfrak g})$,  %
\begin{align*}
 ( \stackrel{\Lambda_1}{\pi}\otimes \stackrel{\Lambda_i}{\pi})
 ({\cal W}_{12}-\widetilde{\cal W}_{12})(x)=0,\ i=2,3,\qquad
 \text{implies}\qquad
 ( \stackrel{\Lambda_1}{\pi}\otimes \stackrel{\Lambda}{\pi})
 ({\cal W}_{12}-\widetilde{\cal W}_{12})(x)=0,
\end{align*}
for any submodule 
$\stackrel{\Lambda}{V}$ of $\stackrel{\Lambda_2}{V}\otimes \stackrel{\Lambda_3}{V}.$
From this lemma we obtain that:
\begin{align*}
   (\stackrel{\Lambda_i}{\pi} \otimes \stackrel{\Lambda_1}{\pi})
   ({\cal W}_{12}-\widetilde{\cal W}_{12})(x)=0,\ i=2,3,\qquad
 \text{implies}\qquad
   (\stackrel{\Lambda}{\pi} \otimes \stackrel{\Lambda_1}{\pi})
   ({\cal W}_{12}-\widetilde{\cal W}_{12})(x)=0
 \end{align*}
for any  submodule 
$\stackrel{\Lambda}{V}$ of $\stackrel{\Lambda_2}{V}\otimes \stackrel{\Lambda_3}{V}.$

It is a basic result that any irreducible  module of  
$U_q(sl(n+1))$ is obtained as a submodule of   some tensor power of the fundamental 
representation. Therefore, 
$(\stackrel{\Lambda_1}{\pi}\otimes \stackrel{\Lambda_2}{\pi})
({\cal W}_{12}-\widetilde{\cal W}_{12})(x)=0$ for any couple 
$(\stackrel{\Lambda_i}{\pi},\stackrel{\Lambda_i}{V})$, $i=1,2$, 
of irreducible representations of $U_q( sl(n+1)).$

This concludes the proof.
\qed

\begin{rem}
Note that Lemma~\ref{lemme-Mweak}, Lemma~\ref{lemme-ABRR}, the uniqueness
condition of Lemma~\ref{uniqueness} and the theorem   give a new  proof of the fact that
$F(x)$ is a zero-weight solution of the QDCE in the $U_q(sl(n+1))$ case.
\end{rem}

\begin{rem}
It is easy to compute the expression of ${\cal  M}^{(\pm)}(x),{\cal  M}^{(0)}(x),
{\cal  M}(x)$ in the fundamental representation from the explicit universal expression given by theorem  \ref{theo-Mexpr}. One obtains: 
\begin{eqnarray*}
&&\stackrel{f}{\pi}({\cal  M}^{(+)}(x))=1+\sum_{1\leq i<j\leq n+1}a_{ij}(-1)^{i-j}
q^{-\frac{(j-i)(j+i-3)}{2(n+1)}}E_{i,j} \\
&&\stackrel{f}{\pi}(({\cal  M}^{(+)}(x))^{-1})=1+\sum_{1\leq i<j\leq n+1}b_{ij}
q^{-\frac{(j-i)(j+i-3)}{2(n+1)}}E_{i,j} \\
&&\text{with}\;\; a_{ij}=\sum_{i<a_1<...<a_{j-i}\leq n+1} 
\nu_{a_1}^{-1}\cdots \nu_{a_{j-i}}^{-1}\;\;,\; b_{ij}=\sum_{j<a_1<...<a_{j-i}
\leq n+1} \nu_{a_1}^{-1}\cdots \nu_{a_{j-i}}^{-1}.
\end{eqnarray*}
One also has:
\begin{eqnarray*}
&&\stackrel{f}{\pi}({\cal  M}^{(-)}(x))=1+\sum_{1\leq i<j\leq n+1}c_{ij}
q^{\frac{(j-i)(j+i-3)}{2(n+1)}}E_{j,i} \\
&&\stackrel{f}{\pi}(({\cal  M}^{(-)}(x))^{-1})=1+\sum_{1\leq i>j\leq n+1}d_{ij}
q^{\frac{(i-j)(j+i-3)}{2(n+1)}}E_{i,j} \\
&&\text{with}\;\; c_{ij}=\frac{\nu_i^{j-i}}{\prod_{r=i+1}^j (1-\nu_i \nu_r^{-1})}\;\; ,\;
 d_{ij}=
\frac{\nu_i^{j-i}}{\prod_{r=j}^{i-1} (1-\nu_i \nu_r^{-1})}. 
\end{eqnarray*}

Using these results an  explicit computation shows that:
\begin{equation*}
\stackrel{f}{\pi}(({\cal  M}(x)^{-1})=D\;\mathcal{V}(x)\;{\mathcal U}(x)D^{-1}
\end{equation*}
with $D,\mathcal{V}(x), {\mathcal U}(x)$ defined in Theorem \ref{lemme1}.

\end{rem}

%%%%%%%%%%%%%%%%%%%%%%%%%%%%%%%%%%%%%%%%%%%%%%%%%%%%%%%%%%%%%%%%%%%%%%%%%%%%%%%
\section{Quantum Dynamical coBoundaries and Quantum Weyl Group}
\label{sectionQWG}
%%%%%%%%%%%%%%%%%%%%%%%%%%%%%%%%%%%%%%%%%%%%%%%%%%%%%%%%%%%%%%%%%%%%%%%%%%%%%%%

%%%%%%%%%%%%%%%%%%%%%%%%%%%%%%%%%%%%%%%%%%%%%%%%%%%%%%%%%%%%%%%%%%%%%%%%%%%%%%%

\subsection{ Primitive Loops and Quantum Coxeter Element }

Let $\mathfrak{g}$ be a finite dimensional simple Lie algebra, and $W$ its 
associated Weyl group. 
We define the shifted action of $W$ on $\mathfrak{h}$ by 
$w.\lambda= w(\lambda+\rho)-\rho.$ 
We will denote $w.x$ the corresponding action on the variable 
$x_i,\  i=1,\ldots,r.$
In the case where $\mathfrak{g}=sl(n+1)$, we can identify the Weyl group with 
the permutation group $S_{n+1}.$
Its shifted action on $\nu_1,\ldots,\nu_{n+1}$ is given by 
$w.(\nu_1,\ldots,\nu_{n+1})=(\nu_{w(1)},\ldots,\nu_{w({n+1})}).$

We will assume in this subsection that $\mathfrak{g}=sl(n+1)$ and that $M(x)$ is the solution to the SQDBP defined in Theorem \ref{theo-Mexpr}.
We will study the behaviour of  ${ P}(x),\ M(x)$ under the shifted action of the Weyl group.

\Proposition{}{
The primitive loop is invariant under the shifted action of the Weyl group, 
i.e.
${ P}(w.x)={ P}(x), \forall w\in W.$}

\proof
Because of the fusion property (\ref{deltaPdeX}), it is sufficient to prove that ${ P}(x)$ is invariant in the fundamental representation. We have computed ${ P}(x)$ in the fundamental representation using the explicit form of $M(x)$, it is given by (\ref{Pinfund}) and depends on $\nu_i$ only through the symmetric polynomials.
As a result, ${ P}(x)$ is invariant under the shifted action of $W.$
\qed

An interesting question is what is the explicit expression of  
${ P}(x)$? From equation (\ref{PM+M-},\ref{MdeXsln}), $P(x)$ can be expressed as a finite product of $q$-exponential.
However, on this expression, the invariance of $P(x)$ under the shifted action of the Weyl group is not explicit.
We can simplify the expression of $P(x)$ by using the quantum Weyl group.

Indeed,  in the $sl(2)$ case, we have the property:

\Proposition{Expression of $P(x)$ for $sl(2).$}{
\vspace{-0.5cm}
\begin{align}
 {P}(x)&= v \; e_{q}^{xe}\; e_{q}^{-(xq^{h+1})^{-1}f}\; B(x)\;e_{q^{-1}}^{-xe}
         \label{explicitPdeX2}\\
  &=\omega\;e_{q^{-1}}^{-xe}\;e_{q^{-1}}^{-x^{-1}e}.
         \label{explicitPdeX1}
\end{align}}

\negvspace

\proof Formula (\ref{explicitPdeX1})
is deduced from (\ref{explicitPdeX2}) using (\ref{saitow}). This last formula is explicitely symmetric under the exchange of $x$ and $x^{-1}.$\qed

This result can be  generalized as follows:

\Theorem{}{${P}(x)$ satisfies 
\begin{align}
{P}(x)&= v \hat{w}_C^{-1}\; Q(x),
\end{align}
where  $Q\in (\mathbb{C}(\nu_1,\ldots,\nu_{n})\otimes U_q({\mathfrak b}^+))^c$ 
is invariant under the shifted action of the Weyl group and $\hat{w}_C$ is a 
quantum analog of the Coxeter element 
$\hat{w}_C= \prod_{i=1}^{n}\hat{w}_{(i)}.$}

\proof 
{}From the  expression (\ref{PM+M-}), we obtain  that
\begin{equation*}
P(x)=v\tau(M^{(+)-1}(x))   (\prod_{k=n}^1{ E}^{[k]}(x))
(\prod_{k=n}^1{ F}^{[k]}(x)) B(x)\ M^{(+)}(x),
\end{equation*}
with 
\begin{equation*}
{ E}^{[k]}(x)=e_{q}^{\nu^{-1}_{k+1}\,q^{\zeta^{(k-1)}}e_{(k)}},\qquad
{ F}^{[k]}(x)=e_{q}^{-\nu_{k+1}\,q^{-\zeta^{(k-1)}-h_{(k)}-1}f_{(k)}}.
\end{equation*}
As a result,
\begin{equation*}
P(x)=v\tau(M^{(+)-1}(x))\prod_{k=n}^1({ E}^{[k]}(x){ F}^{[k]}(x))B(x) M^{(+)}(x).
\end{equation*}
Using Saito's formula \eqref{saitow}, we have
\begin{equation*}
 q^{-\frac{h}{2}}\hat{w}^{-1}e_{q^{-1}}^{-e}q^{-\frac{h^2}{2}}
 =e_q^{e}\cdot e_q^{(-q^{-1-h}f)},
\end{equation*}
which implies that
\begin{equation*}
{ E}^{[k]}(x)\; { F}^{[k]}(x)=G_k\; \hat{w}_{(k)}^{-1}\; 
e_{q^{-1}}^{-e_{(k)}}\; H_k
\end{equation*}
with $G_k=q^{-\frac{h_{(k)}}{2}}\nu_{k+1}^{-\frac{h_{(k)}}{2}}
q^{\frac{h_{(k)}\zeta^{(k-1)}}{2}}$ and 
$ H_k^{-1}=G_k q^{\frac{h_{(k)}^2}{2}}q^{\frac{h_{(k)}}{2}}.$
We therefore have proven that
\begin{equation*}
P(x)=v\tau(M^{(+)-1})\prod_{k=n}^1(G_k \hat{w}_{(k)}^{-1} 
e_{q^{-1}}^{-e_{(k)}} H_k)B(x) M^{(+)}(x),
\end{equation*}
and it remains to move all the $ \hat{w}_{(k)}^{-1}$ on the left.
Using the fact that $s_{\alpha_1}\ldots s_{\alpha_k}(\alpha_{k+1})$ 
is a positive root, 
we obtain that
 $(\prod_{p=1}^k\hat{w}_{(p)})e_{(k)} (\prod_{p=1}^k\hat{w}_{(p)})^{-1} \in
 ( U_q(\mathfrak{b}_+))^c.$
Because $\hat{w}_C e_{k}\hat{w}_C^{-1}\in (U_q(\mathfrak{b}^+))^c,$ for $1\leq k\leq n-1,$
we obtain  $\hat{w}_C \tau (M^{(+)-1})\hat{w}_C^{-1}\in 
(U_q(\mathfrak{b}^+))^c.$ This finishes  the proof.
\qed

\begin{rem} 
An interesting question would be to find  the simplest  exact form of $Q(x)$ expressed as a finite product of $q$-exponentials when $n\geq 2,$ exhibiting invariance under the shifted action of the Weyl group.
\end{rem}
 
%%%%%%%%%%%%%%%%%%%%%%%%%%%%%%%%%%%%%%%%%%%%%%%%%%%%%%%%%%%%%%%%%%%%%%%%%%%%%%%
\subsection{Dynamical coBoundary and Dynamical Quantum Weyl Group}

We show in this section that the shifted action of $W$ on $M(x)$ is controlled via a dynamical quantum Weyl group \cite{EV}.
Although the notion of dynamical Weyl group can be defined for any simple Lie algebra, we will assume here that $\mathfrak{g}=sl(n+1)$.

We first recall that the dynamical Weyl group controls the shifted action of 
$W$ on $F(x)$, where $F(x)$  denotes the standard solution to the QDCE. 
Etingof and Varchenko have constructed  a map 
$A:W\rightarrow  ({\mathbb C}(\nu_1,\ldots,\nu_n)\otimes U_q(\mathfrak g))^c,\ 
w\mapsto A_w(x),$ 
satisfying the following properties:
\begin{align}
& \hat{w} A_w(x)  \;\;\text{is a zero-weight element}, \label{dynWeyl1}\\
& \Delta A_{w}(x)\, F(x)=F(w.x)\, (A_w)_{2}(x)\,(A_w)_{1}(xq^{h_2}),
      \label{dynWeyl2}\\
& A_{ww'}(x)=A_{w}(w'.x)\, A_{w'}(x),\  \forall w,w' \in W.\label{dynWeyl3}
\end{align}
The third property can be reformulated as $A_w(x)$ being a one $W$-cocycle 
taking values in  
$({\mathbb C}(\nu_1,\ldots,\nu_n)\otimes U_q(\mathfrak{g}))^c.$

We have the following proposition:

\Proposition{}{Let $M$ denote the solution of the SQDBP defined in 
Theorem ~\ref{theo-Mexpr}. 
We can define a map 
$\tilde{A}:W\rightarrow  ({\mathbb C}(\nu_1,\ldots,\nu_n)\otimes 
U_q(sl(n+1)))^c,\  w\mapsto\tilde{A}_w(x)$ by
\begin{equation}
  \tilde{A}_w(x)=M(w.x)\,M(x)^{-1}.\label{A=M}
\end{equation}
This map satisfies all the properties (\ref{dynWeyl1},\ref{dynWeyl2},\ref{dynWeyl3}).
}

\proof
 We first show that $\hat{w}
 \tilde{A}_w(x)$ is a zero-weight element. We have:
\begin{align*}
\tilde{A}_w(x)_1\, P_2(xq^{h_1})
  &=M_1(w.x)\, M_1(x)^{-1}{ P_2 }(xq^{h_1}) \\
  &=M_1(w.x)\, R_{12}^J\, { P}_2(x)\, R_{21}^J\, M_1(x)^{-1} \\
  &=M_1(w.x)\, R_{12}^J\, { P}_2(w.x)\, R_{21}^J\, M_1(x)^{-1} \\
  &=M_1(w.x)\, M_1(w.x)^{-1}{ P}_2((w.x)q^{h_1})\, M_1(w.x)\, M_1(x)^{-1} \\
  &={ P}_2((w.x)q^{h_1})\, \tilde{A}_w(x)_{1} \\
  &={ P}_2(xq^{w(h)_1})\, \tilde{A}_w(x)_{1}. 
\end{align*}
This shows that $(\hat{w}
 \tilde{A}_w(x))_1$ commutes with ${ P}_2(xq^{h_1})$. Therefore,  using the 
same argument as in Proposition \ref{theoreme}, we obtain that $\hat{w}
 \tilde{A}_w(x)$ is a zero-weight element.

We now prove that ${\tilde A}_w(x)$ satisfies the property (\ref{dynWeyl2}).
The dynamical coboundary equation implies that
\begin{equation*}
F(w.x)=\Delta M(w.x)\, J\,  M_2(w.x)^{-1} M_1(w.x q^{h_2})^{-1}.
\end{equation*}
Therefore,
\begin{align*}
\Delta M(x)^{-1}\, F(x)&=J\, M_2(x)^{-1} M_1(x q^{h_2})^{-1 }\\
&=\Delta M(w.x)^{-1}\, F(w.x)\, M_1(w.x q^{h_2})\, M_2(w.x)\, M_2(x)^{-1}
   M_1(x q^{h_2})\\
&=\Delta M(w.x)^{-1}\, F(w.x)\, M_2(w.x)\,  M_2(x)^{-1} M_1((w.x) q^{wh_2})\, 
   M_1(x q^{h_2})\\
&=\Delta M(w.x)^{-1}\, F(w.x)\, \tilde{A}_w(x)_{2}\, \tilde{A}_w(xq^{h_2})_{1},
\end{align*}
which ends the proof of the proposition.
\qed

One can precise the relation between $A$ and $\tilde{A}.$

\Proposition{}{
If $A,A'$ are two maps $W\rightarrow  ({\mathbb C}(\nu_1,\ldots,\nu_n)\otimes U_q(\mathfrak g))^c$ satisfying the axioms $(\ref{dynWeyl1},
\ref{dynWeyl2}, \ref{dynWeyl3}),$ 
they are related as $A'_w=A_w Y_w$, where $Y$ is a map
 $Y: W\rightarrow  ({\mathbb C}(\nu_1,\ldots,\nu_n)\otimes U_q(\mathfrak h))^c$ which is a one-cocycle, i.e.
\begin{equation}
  Y_{ww'}(x)= \hat{w}'(Y_{w}(w'.x))\hat{w}'{}^{-1}\  Y_{w'}(x),
\end{equation}
and such that each $Y_w(x)$ is a dynamical group like element, i.e.
\begin{equation}
\Delta(Y_{w}(x))=Y_w(x)_2 Y_w(xq^{h_2})_1.
\end{equation}
}

\vspace{-0.7cm}
\proof
We define $Y_{w}=A_{w}^{-1}A_{w}'=(\hat{w}A_{w})^{-1}(\hat{w}A_{w}')$ which is 
zero-weight.
We have,
\begin{equation}
F(x)^{-1}\Delta(Y_{w}(x))\, F_{12}(x)=A_{w}(xq^{h_2})_1^{-1}\,
Y_{w}(x)_2 \,
A_{w}'(xq^{h_2})_1=Y_w(x)_2\, Y_w(xq^{h_2})_1.
\end{equation}
We now apply the Lemma 2.15 of Etingof-Varchenko \cite{EV2} and we obtain that
 $Y_{w}(x)$ lies in $({\mathbb C}(\nu_1,\ldots,\nu_n)\otimes U_q(\mathfrak h))^c.$
This ends the proof.
\qed

\begin{rem} It would be interesting to precise the value $Y_w(x)$ relating  $A_w$  the dynamical quantum Weyl group of Etingof-Varchenko to  $\tilde{A}_w(x)$ defined by (\ref{A=M}).
\end{rem} 

\begin{rem}The solutions of (\ref{dynWeyl1},\ref{dynWeyl2},\ref{dynWeyl3}) 
exist for all simple Lie algebra $\mathfrak{g}$, however we give a solution of these equations as a coboundary of the form (\ref{A=M}) only in the $sl(n+1)$ case.
It would be an interesting problem to study if the dynamical Weyl group of Etingof-Varchenko can be written as $A_w(x)=M(w.x)M(x)^{-1}Y_w(x)$ with $Y_w(x)$ satisfying properties of the previous proposition and $M:\mathbb {C}^r\rightarrow U_q(\mathfrak{g}).$ In the negative this would provide an alternative proof of the Balog-D\c abrowski-Feh\'er Theorem.
\end{rem}

%%%%%%%%%%%%%%%%%%%%%%%%%%%%%%%%%%%%%%%%%%%%%%%%%%%%%%%%%%%%%%%%%%%%%%%%%%%%%%%
\section{Quantum reflection algebras}
\label{sectionQRA}
%%%%%%%%%%%%%%%%%%%%%%%%%%%%%%%%%%%%%%%%%%%%%%%%%%%%%%%%%%%%%%%%%%%%%%%%%%%%%%%

%%%%%%%%%%%%%%%%%%%%%%%%%%%%%%%%%%%%%%%%%%%%%%%%%%%%%%%%%%%%%%%%%%%%%%%%%%%%%%%
\subsection{Definitions, properties and primitive representations}

In this section we define and study the quantum reflection algebra 
${\cal L}_q({\mathfrak g},J)$ associated to any cocycle 
$J\in (U_q(\mathfrak{g}))^{\otimes 2}$. 
Although there exists a morphism of algebra 
 ${\cal L}_q({\mathfrak g},J)\rightarrow U_q(\mathfrak g)$, this morphism of algebra is not an isomorphism. 
In particular, in the case where $\mathfrak{g}=sl(n+1)$ and  $J=J^{\tau}$,  
${\cal L}_q( sl(n+1),J)$ admits one-dimensional representations which do not 
extend to   $U_q(sl(n+1))$-representations. 
We call them primitive representations and study their relation with the primitive loop.

\Definition{Quantum Reflection Algebra}{
For any simple Lie algebra $\mathfrak{g}$ and any cocycle $J\in  (U_q(\mathfrak{g}))^{\otimes 2}$, we define  the algebra ${\cal L}_q({\mathfrak g},J)$ as being the associative algebra generated by the components ${\mathfrak P}_ {V}$ for all finite dimensional module $V$ of the element  ${\mathfrak P} \in ( U_q({\mathfrak g})\otimes 
{\cal L}_q({\mathfrak g},J))^c $ with the relations:
\begin{equation}
(\Delta^{J} \otimes id)({\mathfrak P})=(R^{J}_{12})^{-1} {\mathfrak P}_1 R^{J}_{12} {\mathfrak P}_{2},
\;\;\;\;\;\;\;\;(\varepsilon \otimes id)({\mathfrak P})=1.\label{quantummoduli}
\end{equation}
}
Let $(V,\pi)$ be a finite dimensional $U_q(\mathfrak g)$-module, we denote
 $\stackrel{\pi}{U} =(\pi\otimes \id)(\mathfrak P)\in \mathrm{End}(V)\otimes {\cal L}_q({\mathfrak g},J).$
These matrices satisfy the following reflection equations for all finite dimensional  $U_q(\mathfrak g)$ modules $V,V':$ 
\begin{equation}
{ R}^{J}_{21} \stackrel{\pi}{U}_1 {R}^{J}_{12} \stackrel{\pi'}{U}_2
=\stackrel{\pi'}{U}_2 { R}^{J}_{21}  \stackrel{\pi}{U}_1 { R}^{J}_{12}.\label{echangeloop}
\end{equation}

Reflection  algebras satisfy the following property:

\Proposition{}{The  map $\kappa:{\cal L}_q({\mathfrak g},J)\rightarrow U_q(\mathfrak{g})$ defined by
\begin{equation}
(\id\otimes \kappa)(\mathfrak{P})=R^{J(-)-1}R^{J}
\end{equation}
is a morphism of algebra.
}

\proof It follows from the quasitriangularity of $R^J.$
\qed
This map is usually thought as being an isomorphism of algebra but this is not true, it is true only if we localize certain elements of
 ${\cal L}_q({\mathfrak g},J).$  This aspect will be central in the rest of this section.

${\cal L}_q({\mathfrak g},J)$ is not a Hopf algebra but it is naturally endowed with a structure of left $U_q({\mathfrak g})$-comodule algebra.

\Proposition{}{The map $\sigma:{\cal L}_q({\mathfrak g},J) \rightarrow 
         U_q({\mathfrak g}) \otimes {\cal L}_q({\mathfrak g},J)$ defined by 
\begin{equation}
 (\id\otimes \sigma)(\mathfrak{P})=R^{J(-)-1}_{12}\mathfrak{P}_{13}R^{J}_{12}
\end{equation}
is a morphism of algebra and is a left $U_q(\mathfrak g)$-coaction with the  coproduct being  $\Delta^J.$}

\proof Trivial.
\qed

Although this map does not define a structure of Hopf algebra on ${\cal L}_q({\mathfrak g},J)$ we still have the following commutative diagram:
\begin{equation}
(\kappa\otimes \id)\sigma=\Delta^J\kappa.
 \end{equation}
 A natural problem  is  the study of the representation theory of 
${\cal L}_q({\mathfrak g},J).$
Because of the two previous propositions we have the following general properties:
\begin{enumerate}
 \item a representation $\pi$ of $U_q(\mathfrak g)$ defines a representation $\pi\circ \kappa$ of $ {\cal L}_q({\mathfrak g},J).$
\item if $(V,\pi)$ is a representation of $U_q(\mathfrak{g})$ and
 $(W,\omega)$ is a representation of  $ {\cal L}_q({\mathfrak g},J),$  one can define the tensor product $\pi\hat{\otimes}\omega=(\pi\otimes \omega)\sigma$ representation of $  {\cal L}_q({\mathfrak g},J)$ acting on $V\otimes W.$
\end{enumerate}

As an  example we first classify the irreducible finite dimensional representations of ${\cal L}_q(sl(2))$ with $J=1.$

We first define a simpler presentation of ${\cal L}_q(sl(2))$ using the fundamental representation of $U_q(sl(2)).$
 ${\cal L}_q(sl(2))$ is generated 
 by the matrix elements 
$a=\stackrel{f}{U}{}^1_1,b=\stackrel{f}{U}{}^1_2,c=\stackrel{f}{U}{}^2_1,d=\stackrel{f}{U}{}^2_2$  with relations (\ref{quantummoduli}) 
which can be explicited respectively as (\ref{echangeloop}) with an additional relation:
\begin{align*}
&ac=q^2 ca \;\;\;\;\;\;\;\;\;\;\;\;ba=q^2 ab\;\;\; \;\;\;\;\;\;\;\;\;bc-cb=(1-q^{-2})a(d-a)\\
&cd-dc=(1-q^{-2})ca \;\;\;\;\;\;\;\;\;\;\;\;db-bd= (1-q^{-2})ab\;\;\;\;\;\;\;\;\;\;\;\;ad=da\\
&\mbox{\rm and}\;\;\;\;\;\;\;\;\;\;\;\;ad-q^2 cb=1.
\end{align*}
The center of ${\cal L}_q(sl(2))$ is  the  polynomial algebra  generated by 
$z=q^{-1}a+qd.$ \\

We denote  $\hat{\cal L}_q(sl(2))$ the algebra ${\cal L}_q(sl(2))$ localized in $a$ and  we denote $a^{-1}$ the inverse of $a$. 
$\hat{\cal L}_q(sl(2))$ is now isomorphic to $U_q(sl(2))$ through the explicit homomorphism
\begin{equation}
\rho: U_q(sl(2)) \longrightarrow {\hat{\cal L}}_q(sl(2)),
\;\;\;\;\;\;\;\;\rho(q^h)=a,\;\;\;\rho(e)=\frac{c}{1-q^{-2}}, \;\;\;\rho(f)=\frac{a^{-1}b}{q-q^{-1}},
\end{equation}
which satisfies 
 $\rho \circ \kappa =id_{\mid {\cal L}_q(sl(2))}.$\\

It is easy to classify the irreducible representations of ${\cal L}_q(sl(2)).$

\Proposition{}{The irreducible representations of ${\cal L}_q(sl(2))$ consists  in: \\
-The representation $\pi\circ \kappa$ where $\pi$ is an irreducible 
 representation of $U_q(sl(2)).$ These representations extends to  
 $\hat{\cal L}_q(sl(2)).$\\
-The one-dimensional representations ${\cal E}_{x,\alpha}$ with $x,\alpha\in \mathbb C^*,$ defined as:
 \begin{equation}
  {\cal E}_{x,\alpha}(a)=0\;\;\;\;\;{\cal E}_{x,\alpha}(d)=q^{-1}(x+x^{-1})\;\;\;\;\;
{\cal E}_{x,\alpha}(b)=q^{-1}\alpha\;\;\;\;\;{\cal E}_{x,\alpha}(c)=-q^{-1}\alpha^{-1}.
\;\;\;\;\;
\end{equation} Note that ${\cal E}_{x,\alpha}(z)=(x+x^{-1})$ and ${\cal E}_{x,\alpha}(\stackrel{f}{U})=D_{\alpha}{\bf P}(x)D_{\alpha}^{-1}$ with $D_{\alpha}^2=
diag (\alpha,\alpha^{-1}).$
 }

\proof
Let $(\Pi, V)$ be a finite dimensional  irreducible  ${\cal L}_q(sl(2))$ module. $z$ being central is represented by a complex number 
$z\in {\mathbb C}.$
If $\Pi(a)$ is invertible one obtains a representation of  $\hat{\cal L}_q(sl(2))\simeq U_q(sl(2))$ which from the classification of irreducible representations of $U_q(sl(2))$ when $q$ is not a root of unit is of the form $\Pi=\pi\circ \kappa.$ 
If not, we define $W=ker (\Pi(a))$, it is a submodule and the restriction of $a,b,c,d$ to $W$ satisfies:
\begin{equation}
a=0,\quad bc=cb,\quad -q^2bc=1,\quad d=q^{-1}z.
\end{equation}
This is an abelian algebra and therefore the irreducible finite dimesional representations are of dimension one. Therefore $W=V$ and is one dimensional. The representation is therefore of the type ${\cal E}_{x,\alpha}.$ 
\qed

We will call the  family of one-dimensional representations of  
${\cal L}_q(sl(2))$ for which  $a=0$ {\it primitive representations.} 

\begin{rem} 
The study of the representation theory of the algebra ${\cal L}_q(\mathfrak g, J)$ is an interesting problem that we will only look at in the case where $\mathfrak{g}=sl(n+1)$ and $J=J^\tau.$ 
\end{rem}

The primitive representation in the $sl(2)$ case appears first in \cite{KSS}

The classification of all characters, i.e. all one dimensional representations, of 
${\cal L}_q(sl(n+1), J)$ for $J=1$ has been obtained in \cite{M}.

\begin{rem} 
Finite dimensional representations of reflection algebras are not completely reducible. An example is given for ${\cal L}_q(sl(2))$ by \cite{Sa}. 
\end{rem}

In the next section we will study the generalisation of these primitive
 representations to the case of ${\cal L}_q(sl(n+1),J^{\tau}).$ 
We will study the decomposition of the tensor product of an irreducible representation $\pi$ of $U_q(sl(n+1))$ with an irreducible representation $\omega$ of 
${\cal L}_q(sl(n+1),J^{\tau}).$
When the irreducible representation $\omega=\pi'\circ\kappa$, where $\pi'$ is a representation of $U_q(sl(n+1)),$ the intertwining map is governed by ordinary Clebsch-Gordan map of $U_q(sl(n+1))$ and there is nothing new, but 
when $\omega$ is a primitive representation the intertwining map is entirely governed by the coboundary element evaluated in the representation $\pi.$ 
We can also invert this process and define the coboundary in 
the representation $\pi$ as being the intertwining map.

We now   assume that  $\mathfrak{g}=sl(n+1)$ and $J=J^\tau.$ 
In the first part of this section  we assume that $M(x)$ is the solution of the SQDBP  with associated primitive loop ${ P}(x).$

\Definition{Primitive Representations }{
The primitive representation ${\cal E}$ of ${\cal L}_q(sl(n+1),J^{\tau})$ is the representation of ${\cal L}_q(sl(n+1),J^{\tau})$ with values in ${\mathbb C}[\nu_1,\nu_1^{-1},\ldots,\nu_r,\nu_r^{-1}] $
defined by 
\begin{equation}
(id \otimes {\cal E}_x)({\mathfrak P})={P}(x)
\end{equation}
We will abusively denote  ${\cal E}_x$ this representation.}

\proof
Trivial. 
\qed

Of course if we fix $\nu_1,\ldots,\nu_r$  to non zero complex numbers, we obtain characters of  ${\cal L}_q(sl(n+1),J^{\tau}).$

We have the following proposition

\Proposition{Intertwining map and Coboundary }{\label{PropIntMapCobound}
Let $(V,\pi)$ be a finite dimensional representation of $U_q(\mathfrak g)$, the following decomposition property holds:
\begin{equation}
\pi\hat{\otimes} {\cal E}_x=\bigoplus_{\lambda\in {\mathfrak h}^* }  
{\cal E}_{xq^\lambda}^{\oplus m_\lambda}\;\;\; \text{where}\;\;
 m_\lambda=dim V[\lambda],
\end{equation}
and $M(x)_V$ is an intertwining map between these representations. 
}

\proof
${\cal E}_x$ acts on the module $\mathbb C$ generated by $1,$ we will show that if $v\in V[\lambda]$ then the action of $\pi\hat{\otimes} {\cal E}_x $ on $w=M(x)_V^{-1}v\otimes 1$ is 
${\cal E}_{xq^\lambda}.$
Indeed
\begin{align*} 
 (\pi\hat{\otimes} {\cal E}_x)(\stackrel{\pi'}{U}{}^a_b)(w\otimes 1)
 &=\pi((\stackrel{\pi'\;.}{R}{}^{J(-)-1}){}^{a}_m
       (\stackrel{\pi'\;.}{R}{}^{J(+)}){}^{n}_b)w\otimes 
       \stackrel{\pi'}{P}(x){}^m_n\;1\\
 &=(\pi'\otimes id)(R^{J(-)-1}_{12}P_1(x)R^{J(+)}_{12})^a_b.(w\otimes 1)
          \\
 &=(M_2(x)^{-1}P_1(xq^{h_2})M_2(x))^{a}_b M_2(x)^{-1} (v\otimes 1)\\ 
 &=(M_2(x)^{-1}P_1(xq^{\lambda}))^{a}_b(v\otimes 1)
  =P(xq^{\lambda})^{a}_b (w\otimes 1).
\end{align*}
The proposition follows.
\qed

%%%%%%%%%%%%%%%%%%%%%%%%%%%%%%%%%%%%%%%%%%%%%%%%%%%%%%%%%%%%%%%%%%%%%%%%%%%%
\subsection{A representation framework for weak solution of dynamical coBoundary elements }
In this subsection we give hints towards a purely   representation theoretical approach of QDBE. Because we already have obtained an explicit universal solution of this problem in its strongest formulation, the aim of this subsection  is not to give completely rigorous reconstruction of $M(x).$ We would like to emphasize the puzzling fact that $M(x)$ and $F(x)$ can be constructed solely with interwining operators involving primitive representations of loop algebras and finite dimensional representations of $U_q(\mathfrak g)$ whereas in the work of \cite{EV} $F(x)$ is built using intertwining operators between Verma modules and finite dimensional representations of  $U_q(\mathfrak g).$

${\cal L}_q(sl(n+1),J^{\tau})$ can be presented in term of  matrix
 elements of $\stackrel{f}{U}.$
We use the results of \cite{GPS} and denote ${\check{\bf R}}^{J}=
 {\bf R}^{J}_{12}P_{12},$ where $P:V_f^{\otimes 2}\rightarrow V_f^{\otimes 2}$ is the permutation operator.

Because $\check{{\bf R}}$ satisfies the Hecke symmetry it is also true for ${\check{\bf R}}^{J}.$
As a result all the constructions of \cite{GPS} concerning the $Tr_q$ and $det_q$  can be applied.
We denote $A^{(k)}$ the antisymmetrizer associated to ${\check{\bf R}}$ acting on $V_f^{\otimes k},$ it is denoted $P_{-}^k$  in \cite{GPS}.
Similarly we can define $A_J^{(k)}$ the antisymmetrizer associated to ${\check{\bf R}}^J$ acting on $V_f^{\otimes k}.$ 
 ${\check{\bf R}}$ and ${\check{\bf R}}^J$ are both Hecke symmetry of rank $n+1$, i.e 
 $A^{(n+2)}=A_J^{(n+2)}=0.$

${\cal L}_q(sl(n+1),J^{\tau})$ is isomorphic to the algebra 
$A(sl(n+1),J)$ generated by the matrix elements of 
$U\in End(V_f)\otimes A(sl(n+1),J)$ with relations 
\begin{align}
   &{ \bf R}^{J}_{21} U_1 {\bf R}^{J}_{12} U_2
    =U_2 {\bf  R}^{J}_{21} U_1 {\bf  R}^{J}_{12}\\
   &det_q(U)=det_q({\bf P}(x))
\end{align}
where 
$det_q(U)$ is the central element defined by
\begin{equation}
 det_q(U)=tr_{1\cdots n+1}(A_{J}^{(n+1)}(U_1
{\check{\bf R}}^J_{12}{\check{\bf R}}^J_{23}... {\check{\bf R}}^J_{n,n+1})^{n+1}).
\end{equation}
The isomorphism is obtained by identifying $\stackrel{f}{U}$ and $U.$ 
\\
Note that $det_q({\bf P}(x))=q^{-n(n+1)}.$

Let ${\bf M}(x)$  be a matrix such that (\ref{RMM-fund}) is satisfied and define 
${\bf  P}(x)=v{\bf M}(x)^{-1}{\bf B}(x){\bf M}(x).$
By construction we have 
\begin{align}
{ \bf R}^{J}_{21} {\bf P}(x)_1 {\bf R}^{J}_{12}  {\bf P}(x)_2
&=  {\bf P}(x)_2{\bf  R}^{J}_{21}  {\bf P}(x)_1 {\bf  R}^{J}_{12},\\
{ \bf R}^{J}_{21} {\bf P}(x)_1 {\bf R}^{J}_{12}
&=  {\bf M}(x)^{-1}_2 {\bf P}(xq^{h_2})_1 {\bf M}(x)_2.\label{PM}
\end{align}

As a result we still define the primitive representation ${\cal E}$ of $ A(sl(n+1),J)$ as the representation of $ A(sl(n+1),J)$ with values in ${\mathbb C}[\nu_1,\nu_1^{-1},\ldots,\nu_r,\nu_r^{-1}] $
defined by 
\begin{equation}
(id \otimes {\cal E}_x)(U)={\bf P}(x).
\end{equation}

Let $(V,\pi)$ be any finite dimensional representation of $U_q(sl(n+1))$, we want to show the following decomposition property:
\begin{equation}
\pi\hat{\otimes} {\cal E}_x=\bigoplus_{\lambda\in {\mathfrak h}^* }  
{\cal E}_{xq^\lambda}^{\oplus m_\lambda}\;\;\; \text{where}\;\;
 m_\lambda=dim V[\lambda]\label{reduc}
\end{equation}
and  will denote ${\cal M}_V(x)$ any intertwiner between these representations.

As a consequence of (\ref{PM}), the same proof as Proposition \ref{PropIntMapCobound} shows that
\begin{equation}
\pi_f\hat{\otimes} {\cal E}_x=\bigoplus_{\lambda\in {\mathfrak h}^* }  
{\cal E}_{xq^\lambda}^{\oplus m_\lambda}\;\;\; \text{where}\;\;
 m_\lambda=dim V_f[\lambda].
\end{equation}
Hence (\ref{reduc}) holds also for $\pi_f^{\otimes p} \circ ( \Delta^J)^{(p)},$ i.e

\begin{equation}(\pi_f)^{\otimes_{_{J}} p}\hat{\otimes} {\cal E}_x=\bigoplus_{\lambda\in {\mathfrak h}^* }  
{\cal E}_{xq^\lambda}^{\oplus m_\lambda}\text{where}\;\;
 m_\lambda=dim V_f^{\otimes_{_{J}} p}[\lambda].
\end{equation}

Let $H_p(q)$ be the $A_p$-Hecke algebra generated by $\sigma_1,...,\sigma_{p-1}.$
If $W$ is an irreducible submodule of $\pi_f^{\otimes p}$ there exists an idempotent element 
${\cal Y}_W={\cal Y}_W(\sigma_1,...,\sigma_{p-1})\in H_p(q)$ such that 
${\cal Y}_W(\check{\bf R}_{12},...,\check{\bf R}_{p-1,p})$ is the projector on the submodule $W.$
As a result if $W$ is the submodule of $\pi_f^{\otimes_{_{J}} p}$ the associated projector is 
$J_{1,...,p}^{-1}{\cal Y}_WJ_{1,...,p}=
{\cal Y}_W(\check{\bf R}^J_{12},...,\check{\bf R}^J_{p-1,p})={\cal Y}_{W}^J$ where 
$(\Delta^J)^{(p)}(a)=J_{1,...,p}^{-1}\Delta^{(p)}(a)J_{1,...,p}.$
It remains to show that

 \begin{equation}({\cal Y}_{W}^J (\pi_f)^{\otimes_{_{J}} p})\hat{\otimes} {\cal E}_x=\bigoplus_{\lambda\in {\mathfrak h}^* }  
{\cal E}_{xq^\lambda}^{\oplus m_\lambda}\text{where}\;\;
 m_\lambda=dim({\cal Y}_{W}^J V_f^{\otimes_{_{J}} p})[\lambda].
\end{equation}

It is straightforward to verify, using (\ref{RMM-fund}), the following properties:
\begin{eqnarray}
{\cal Y}_{W}^J {\bf M}_{1,..,p}(x)^{-1}&=& {\bf M}_{1,..,p}(x)^{-1}{\cal Y}_{W}(\check{\bf R}_{12}(x),...,\check{\bf R}_{p-1,p}(x))\\
&=& {\bf M}_{1,..,p}(x)^{-1}F_{1,..,p}(x)^{-1}
{\cal Y}_W(\check{\bf R}_{12},...,\check{\bf R}_{p-1,p})F_{1,..,p}(x)
\end{eqnarray}
 where \begin{eqnarray*} {\bf M}_{1,..,p}(x)&=& {\bf M}_1(xq^{h_1+\cdots +h_{p}})\cdots  {\bf M}_{p-1}(xq^{h_p})
{\bf M}_p(x)\\
F_{1,..,p}(x)&=&{\bf F}_{(1...p-1,p)}(x){\bf F}_{(1...p-2,p-1)}(xq^{h_p})...{\bf F}_{1,2}(xq^{h_p+...+h_3}),
\end{eqnarray*}
with ${F}_{(1...k,k+1)}=(\Delta^{(k)}\otimes id)(F(x)),$
which concludes the proof. \\

Having defined  ${\cal M}_V(x)$ for $V$ simple we straightforwardly define ${\cal M}_W(x)$ for $W$ semi-simple. Let $(\pi,V),(\pi',W)$ be two representations of $U_q(sl(n+1)).$ Using decomposition property \eqref{reduc} we obtain the property $$J^{-1}_{V,W}{\cal M}_{V\otimes W}(x)^{-1}\;(V\otimes_{{}_{J}}W)[\lambda]=({\cal M}_{V}(xq^{h_W}){\cal M}_{W}(x))^{-1}\;(V\otimes_{{}_{J}}W)[\lambda].$$ As a result, ${\cal F}_{V,W}(x)$ defined by $${\cal F}_{V,W}(x)={\cal M}_{V\otimes W}(x)J_{V,W}({\cal M}_{V}(xq^{h_W}){\cal M}_{W}(x))^{-1}$$ is of zero weight. The family of intertwing maps ${\cal M}_V(x)$ defines therefore a solution of the WQDBP.

The  previous framework  shows that  one can obtain a definition of 
${\cal M}_V(x)$ in a purely representation theoretical setting and that  
${\cal M}$ is a solution of the WQDBP.
We think that this method can be further pursued to obtain a purely representation theoretical approach to the SQDBP. As an example, we give in Appendix~\ref{sec-constrM} some remarks concerning this point in the rank one case, the higher rank case being still unclear for us.
%%%%%%%%%%%%%%%

%%%%%%%%%%%%%%%%%%%%%%%%%%%%%%%%%%%%%%%%%%%%%%%%%%%%%%%%%%%%%%%%%%%%%%%%%%%%%%
\section{Conclusion}
\label{sectionConcl}
%%%%%%%%%%%%%%%%%%%%%%%%%%%%%%%%%%%%%%%%%%%%%%%%%%%%%%%%%%%%%%%%%%%%%%%%%%%%%%

In this work we have given a universal explicit solution of the Quantum Dynamical coBoundary Equation. This was obtained through the use of the primitive loop,  which study led us to this solution.

However many points are still  unclear to us.

The first one concerns the   Balog-D\c abrowski-Feh\'er Theorem. Although the 
result is unquestionable, the proof seems unnatural. The fact that it selects precisely $sl(n+1)$ and the 
 Cremmer-Gervais 's $r$-matrix  still remains  unclear.

The second one comes from the possible  various generalizations. We have studied here the QDBE  in the case where $F(x)$ is the standard solution and $\mathfrak{h}$ is the Cartan sub-algebra, but we could imagine considering also the case when $F(x)$ is associated to a generalized Belavin-Drinfeld triple of the type considered by O. Schiffmann \cite{Sc}, or generalizing this equation to the non abelian case. These problems are still completely open.
   
As a third point, one may also wonder whether one can generalize the straightforward proof of  
Appendix~\ref{sec-constrM}, presented here  in the $U_q(sl(2))$ case, to higher rank.  

Finally, we would like to mention that the coboundary equation originates in the IRF-Vertex transform \cite{Bax}, and all the tools are now present for the construction of a universal IRF-Vertex transform in the quantum affine case. This universal coboundary element will relate  the face type twistors  and the vertex type twistors of elliptic quantum algebras of \cite{JKOS}.

\newpage

%%%%%%%%%%%%%%%%%%%%%%%%%%%%%%%%%%%%%%%%%%%%%%%%%%%%%%%%%%%%%%%%%%%%%%%%%%%%%%
%%%%%%%%%%%%%%%%%%%%%%%%%%%%%%%%%%%%%%%%%%%%%%%%%%%%%%%%%%%%%%%%%%%%%%%%%%%%%%
\section{Appendix}
%%%%%%%%%%%%%%%%%%%%%%%%%%%%%%%%%%%%%%%%%%%%%%%%%%%%%%%%%%%%%%%%%%%%%%%%%%%%%%

%%%%%%%%%%%%%%%%%%%%%%%%%%%%%%%%%%%%%%%%%%%%%%%%%%%%%%%%%%%%%%%%%%%%%%%%%%%%%%
\subsection{The Balog-D\c abrowski-Feh\'er Theorem}
\label{sec-BDF}

We give here elements of the proof of Theorem \ref{lemme2} following the 
arguments of \cite{BDF}.

\medskip

For a finite dimensional simple Lie algebra $\mathfrak{g}$, 
let $R(x)$ be the universal standard solution of the QDYBE
and assume that there exist $M(x)$, $R^J$ such that the equation 
(\ref{RMM=MMR}) holds.
We fix $x$ and expand each factor of \eqref{RMM=MMR} in term of $\hbar$, 
with $q=e^{\frac{\hbar}{2}}$.
We have
$R(x)=1+\hbar\, r(x)+o(\hbar)$, where  
\begin{equation}
 r(x)=\frac{1}{2}\Omega_{\mathfrak h}+\sum_{\alpha\in \Phi}
 r_{\alpha}(x)\, e_\alpha\otimes f_{\alpha}
 \quad\text{with}\quad 
 r_{\alpha}(x)=\frac{(\alpha,\alpha)}{2}
           \Big(1-\prod_j x_j^{2\alpha(\zeta^{\alpha_j})}\Big)^{-1}.
\end{equation}
As a result, $R^J=1+\hbar\, r_J+o(\hbar)$, and the linear term $r_J$ is given 
through \eqref{RMM=MMR} in terms of $r(x)$ as
\begin{equation}
r_{J}=M_1(x)^{-1}M_2(x)^{-1}\Big(r(x)+\frac{1}{2}\sum_{j=1}^r A_j(x)\wedge h_{\alpha_j}\Big)\,M_1(x)\, M_2(x),
\end{equation}
where $A=A_i dx^i$ is a flat connection defined as 
$A_i=x_i (\partial_{i}M)M^{-1}\in \mathfrak{g}$. 
The condition $\partial_{i}r_{J}=0$  can then be expressed only in term of 
$r(x)$ and of the connection $A$, and reads:
\begin{equation}\label{cond-rJ}
  x_i\partial_i\Big(r(x)+\frac{1}{2}\sum_{j=1}^r A_j\wedge h_{\alpha_j}\Big)+
  \Big[r(x)
  +\frac{1}{2}\sum_{j=1}^r A_j\wedge h_{\alpha_j}, A_i\otimes 1+1\otimes A_i
  \Big]=0.
\end{equation}

\bigskip

Balog, D\c abrowski and Feh\'er have shown that the set of flat connections 
satisfying this equation is empty when $\mathfrak{g}$ does not belong to the 
$A_n$ series.
In order to prove this result, we decompose $A$ on the root subspaces as 
\begin{equation}\label{Aroot}
  A_j=\sum_{i=1}^r A_j^i \, h_{\alpha_i}
    +\sum_{\alpha\in \Phi}A_j^{\alpha} \, e_\alpha.
\end{equation}
The differential equations satisfied by $A$ are the flatness condition 
$D_A A=0$ and the equation $D_A(r(x)+\frac{1}{2}\sum_{j=1}^r A_j\wedge h_j)=0$
(Eq. \eqref{cond-rJ}),
which give respectively, when projected on the root subspaces,
\begin{align}
 &x_i\partial_i A_j^m -x_j\partial_j A_i^m 
   -\frac{2}{(\alpha,\alpha)}\sum_{\alpha\in \Phi}
    A_i^\alpha A_j^{-\alpha} \alpha(\zeta^{\alpha_m})=0,
       \label{flatness1}\displaybreak[0]\\
 &x_i\partial_i A_j^\alpha-x_j\partial_j A_i^\alpha
   -\sum_{\substack{\beta,\gamma \\ \beta+\gamma=\alpha}}
    N_{\beta\gamma}^{\alpha}A_i^\beta A_j^\gamma
   -A_j^\alpha\sum_{n}A_i^n \alpha(h_{\alpha_n})
   +A_i^\alpha\sum_n A_j^n \alpha(h_{\alpha_n})=0,
       \label{flatness2}\displaybreak[1]\\
% &\sum_{m,n}x_i\partial_i A_n^m h_m\wedge h_n
%   +\sum_{\alpha,j}A_j^\alpha A_i^{-\alpha}h_\alpha\wedge h_j=0,
 &x_i\partial_i ( A_n^m - A_m^n ) 
  + \frac{2}{(\alpha,\alpha)}\sum_{\alpha\in \Phi} \big[
     A_n^\alpha A_i^{-\alpha} \alpha(\zeta^{\alpha_m})
     -A_m^\alpha A_i^{-\alpha} \alpha(\zeta^{\alpha_n}) \big] =0,
       \label{DAR1}\displaybreak[0]\\
% &\frac{1}{2}\sum_{n}x_i\partial_i A_n^\alpha e_\alpha\otimes h_n
%   +\frac{1}{2}\sum_{i} A_i^\alpha e_\alpha\otimes t_\alpha-r_\alpha(x)
%    A_i^\alpha e_\alpha\otimes h_\alpha+
%       \label{DAR2}\\
% &+\sum_{j,\beta,\gamma, \alpha=\beta+\gamma}
%  N_{\alpha\beta}^{\gamma}A_i^\gamma A_j^\beta e_\alpha\otimes h_j
%  +\frac{1}{2}\sum_{j,m} (A_j^m A_i^\alpha -A_j^\alpha A_i^m -A_m^j A_i^\alpha)
%   \alpha(h_m)e_\alpha\otimes h_j=0,
%       \nonumber
 &\frac{1}{2} x_i \partial_i A_n^\alpha 
  +\Big(\frac{1}{2} -\frac{2}{(\alpha,\alpha)} 
           r_\alpha(x)\Big) A_i^\alpha  \alpha(\zeta^{\alpha_n})
                    \nonumber\\
 &\hspace{2cm}
  +\sum_{\substack{\beta,\gamma \\ \beta+\gamma=\alpha}} 
   N_{\alpha\beta}^{\gamma}A_n^\beta A_i^\gamma 
  +\frac{1}{2}\sum_{m} (A_n^m A_i^\alpha -A_n^\alpha A_i^m -A_m^n A_i^\alpha)
   \alpha(h_{\alpha_m})=0, 
          \label{DAR2}\displaybreak[0]\\
 &x_i\partial_i r_{\alpha}
  -\frac{1}{2}A_i^\alpha \sum_j A_j^{-\alpha} \alpha(h_{\alpha_j})
  -\frac{1}{2}A_i^{-\alpha} \sum_j A_j^{\alpha} \alpha(h_{\alpha_j})=0,
       \label{DAR3}\displaybreak[0]\\
 &A_i^\alpha \sum_j A_j^\beta \alpha (h_{\alpha_j})
  -A_i^\beta \sum_j A_j^\alpha \beta(h_{\alpha_j})
   -2\Big(r_\alpha(x)-\frac{(\alpha,\alpha)}{(\beta,\beta)}r_{-\beta}(x)\!\Big)
    N_{-\alpha,\alpha+\beta}^{\beta}A_{i}^{\alpha+\beta}=0,
       \label{DAR4}
\end{align}
for $\alpha \ne -\beta$, where we have denoted 
$[e_{\alpha}, e_{\beta}]=N_{\alpha,\beta}^{\alpha+\beta}e_{\alpha+\beta}.$

Combining equations (\ref{flatness1},\ref{DAR1}) on the one hand, and 
equations (\ref{flatness2},\ref{DAR2}) on the other hand, 
one obtains respectively the following linear 
equations for $A_i^n$ and  $A_i^\alpha$:
\begin{align}
  &x_m\partial_m A^n_i-x_n\partial_n A^m_i=0,
             \label{rotA=0}\\
  &x_n\partial_n A_i^\alpha
   +\alpha(\zeta^{\alpha_n})F_{\alpha}(x)A_i^\alpha 
   -A_i^\alpha\sum_m A_m^n\alpha(h_{\alpha_m})=0,  
             \label{diffeqAalpha}
\end{align}
where $F_\alpha(x)=-(1+\prod_j x_j^{2\alpha(\zeta^{\alpha_j})})
(1-\prod_j x_j^{2\alpha(\zeta^{\alpha_j})})^{-1}.$
The general solution of the equation (\ref{rotA=0}) is 
\begin{equation}\label{A-der-phi}
   A_i^n=x_n\partial_n \phi_i,
\end{equation} 
where $\phi_i$ are arbitrary functions.
The  general solution of the equation (\ref{diffeqAalpha}) is then
\begin{equation}\label{A-alpha}
  A_i^\alpha=C_i^\alpha\, (1+F_\alpha)\,\prod_{j}x_j^{-\alpha(\zeta^{\alpha_j})}
      \exp\Big(\sum_m \phi_m \alpha(h_{\alpha_m})\Big),
\end{equation}
in terms of some (so far arbitrary) constants $C_i^{\alpha}$.

Let us define the weight $C^{\alpha}$ by $C^{\alpha}=\sum_i C_i^{\alpha}  \alpha_i^{\vee},$
the equations (\ref{DAR3},\ref{DAR4}) become algebraic equations:
\begin{align}
&\frac{(\alpha,\alpha)}{2}\alpha +C^\alpha
(C^{-\alpha},\alpha)+C^{-\alpha}(C^\alpha, \alpha)=0,\label{DAR3'}\\ 
&C^\alpha  (C^\beta, \alpha)-C^\beta (C^\alpha, \beta)
+C^{\alpha+\beta}N^{\alpha,\beta}_{\alpha+\beta}=0,\label{DAR4'}
\end{align}
where we have denoted $N^{\alpha,\beta}_{\alpha+\beta}
=\frac{(\alpha,\alpha)}{2}N^{\beta}_{-\alpha,\alpha+\beta}.$

The first equation is uniquely solved by decomposing  $C^{\alpha}=
\frac{c(\alpha)}{2}(\alpha+K^\alpha)$ for $\alpha>0$ and $K^\alpha \perp \alpha.$
As a result we obtain that $C^{-\alpha}=\frac{1}{2c(\alpha)}(-\alpha+K^\alpha)$ for $\alpha>0.$

It remains to show that the set of equations (\ref{DAR4'}) rules out all the simple Lie algebras except the $A_n$ series.

Proving this property is simplified by the following observation:
if (\ref{DAR4'}) admits a solution $(C^\alpha)$ for a Lie algebra associated 
to the Dynkin diagram $D$ labelling the simple roots $\alpha_1,\ldots,\alpha_r$, 
and if $D'$ is a connected subdiagram of $D$ associated to the roots 
$\alpha_{j},$ $j\in D'\subset D=\{1,\ldots, r\}$, then the orthogonal projection 
of $(C^\alpha)$ on the vector space generated by $\alpha_{j}$, $j\in D'$, 
is a solution of (\ref{DAR4'}) for the Lie algebra generated by $D'$.
As a result, one obtains that it is sufficient to show that the solution of (\ref{DAR4'}) is empty for $D_4$, $B_2$ and $G_2$ for ruling out all but the 
$A_n$ series.
The next observation comes from the theorem that if $\alpha,\beta$ are roots
such that   $\alpha+\beta$ and $\alpha-\beta$ are simultaneously non roots, 
then  $(\alpha,\beta)=0.$
As a result, in this case, by combining the equations (\ref{DAR4'}) for 
$\alpha+\beta$ and $\alpha-\beta$, one obtains that $K^\alpha\perp \beta.$
As a result, we obtains that $K^\alpha\perp V_\alpha$ where
$V_{\alpha}=\mathbb{C}\alpha+\sum_{\beta\in \Phi,\alpha\pm \beta\notin \Phi}\mathbb{C}\beta.$
We trivially have $V_{w\alpha}=w V_{\alpha}$ where $w$ is any element of the Weyl group.
As a result, $\mathrm{codim} (V_{\alpha})$ only depends on the length of $\alpha.$

An elementary analysis ot the root system of $D_4$ and $G_2$ proves that  
$\mathrm{codim} (V_{\alpha})=0$  for all roots.
Therefore $K_{\alpha}=0$, but the corresponding $C^{\alpha}$ is not a solution of (\ref{DAR4'}).

In the case of $B_2$  one obtains that $\mathrm{codim} (V_{\alpha})=0$ if 
$\alpha$ is long and $\mathrm{codim} (V_{\alpha})=1$ if $\alpha$ is short. 
The explicit study of the system (\ref{DAR4'}) shows that once again the set of solutions is empty.

This concludes the proof that the coboundary equation can admit solutions only 
in the case where $\mathfrak{g}=sl(n+1)$.

\bigskip

We will now show that, in the case where $\mathfrak{g}=sl(n+1)$, 
such a solution $r_J$ is unique up to an automorphism. 
As we know that $r_{\tau,s}$, with $s$ defined as (\ref{sCG}), is 
a solution of the coboundary equation in this case 
(this is a direct consequence of Theorem~\ref{lemme1} in the fundamental 
representation, and of Theorem~\ref{theo-Mexpr} at the universal level),
it means that, for any solution $r_J$ of the coboundary equation,  
there exists an automorphism $\phi$ of $sl(n+1)$ such that 
$r_J=(\phi\otimes \phi)(r_{\tau,s})$.
 
In order to prove this uniqueness property, let us caracterise completely all
the possible solutions for the connection $A$ in the $sl(n+1)$ case.
In this case, positive roots are of the form $u_i-u_j$, $i<j$, where 
$u_1,\ldots,u_{n+1}$ are orthonormal vectors. From the previous 
considerations, $K^{u_i-u_j}$ should be simultaneously orthogonal to 
$u_i-u_j$ and to all $u_k-u_l$ such that $\{k,l\}\cap \{i,j\}=\emptyset$.
It is therefore of the form
\begin{equation*}
   K^{u_i-u_j}=\epsilon_{ij}\Big(u_i+u_j-\frac{2}{n+1}\sum_{k=1}^{n+1}u_k\Big),
\end{equation*}
and it can be shown, using Eq.\eqref{DAR4'}, that all the constants 
$\epsilon_{ij}$ are equal to some common value $\epsilon\in\{+1,-1\}$. 
Furthermore, still
from  Eq.\eqref{DAR4'}, the constants $c(\alpha)$ associated to positive roots
have to satisfy
\begin{equation*}
   c(u_i-u_j)\, c(u_j-u_k)=\epsilon\, c(u_i-u_k),\quad \text{if\ } i<j<k, 
\end{equation*}
which means that there exists some constants $c_1,\ldots,c_n$ such that,
for all positive roots $\alpha$, we have
$c(\alpha)=\epsilon\, \exp\big(\sum_m c_m \alpha(h_{\alpha_m})\big)$. 
On the other hand, pluging \eqref{A-der-phi} and \eqref{A-alpha} into 
\eqref{flatness1}, we obtain the following conditions on the functions 
$\phi_i$:
\begin{equation*}
   \phi_j=\frac{1}{2}\sum_{\alpha\in\Phi^+} K_j^\alpha 
    \log\Big(\prod_l x_l^{\alpha(\zeta^{\alpha_l})}
            -\prod_l x_l^{-\alpha(\zeta^{\alpha_l})}\Big)
    +x_j \partial_j\psi,
\end{equation*}
where $\psi$ is an arbitrary function. Note at this stage that we may as well
fix $c(\alpha)=\epsilon$ and absorb the arbitrariness of this constant in $\psi$. 

Let us now prove that a modification of the function
$\psi$ does not affect $r_J$. 
Let us denote by $A_{(\psi)}$, $M_{(\psi)}$, $(r_J)_{(\psi)}$ 
(respectively by $A_{(0)}$, $M_{(0)}$, $(r_J)_{(0)}$ )  
the connection, the coboundary and the corresponding $r$-matrix associated
to a given $\psi$ (respectively to $\psi=0$).
We have:
\begin{align}
  &  M_{(\psi)}= M_{(0)}\cdot g_\psi,\\
  &  (A_{(\psi)})_j= x_j \partial_j g_\psi \cdot g_\psi^{-1} 
         + g_\psi\cdot (A_{(0)})_j\cdot g_\psi^{-1},\label{Apsi}
\end{align}
where $g_\psi=\exp (\sum_m x_m \partial_m\psi \; h_{\alpha_m} )$.
Note that, due to the specific form of $g_\psi$, only the second term in 
\eqref{Apsi} gives a non-zero contribution to 
$(A_{(\psi)})_j\wedge h_{\alpha_j}$, 
and therefore, 
\begin{align}
  (r_J)_{(\psi)} 
  & = (M_{(0)})_1^{-1} (M_{(0)})_2^{-1} 
                   (g_\psi)_1^{-1} (g_\psi)_2^{-1}
                  \Big( r(x) 
   %+\frac{1}{2}\sum_{j=1}^n x_j \partial_j g_\psi \cdot g_\psi^{-1}\wedge h_{\alpha_j}  
   +\frac{1}{2}\sum_{j=1}^n (g_\psi \,(A_{(0)})_j\, g_\psi^{-1})
                          \wedge h_{\alpha_j}
      \Big) \nonumber\\
  &\hspace{7cm}\times
         (g_\psi)_1\,(g_\psi)_2\, (M_{(0)})_1\, (M_{(0)})_2\nonumber\\
  &= (M_{(0)})_1^{-1} (M_{(0)})_2^{-1} 
    \Big( r(x)+\frac{1}{2}\sum_{j=1}^n (A_{(0)})_j
                          \wedge h_{\alpha_j}
      \Big)\,  (M_{(0)})_1\, (M_{(0)})_2=(r_J)_{(0)}.
\end{align}

Finally, the only arbitrariness in $r_J$ is due to gauge transformations of 
the form 
\begin{align*}
   &M(x) \rightarrow M(x)\cdot u,\\
   &r_J  \rightarrow (u\otimes u)^{-1} r_J \, (u\otimes u),
           \quad u\in\mathfrak{g},
\end{align*}
and to the automorphism $\alpha\rightarrow -\alpha,\ \forall\alpha\in\Phi$,
corresponding to the change $\epsilon\rightarrow -\epsilon$.
Thus, the solution $r_J$ is unique up to an automorphism of $sl(n+1)$, 
which concludes the proof. 

\newpage
%%%%%%%%%%%%%%%%%%%%%%%%%%%%%%%%%%%%%%%%%%%%%%%%%%%%%%%%%%%%%%%%%%%%%%%%%%%%%%%
\subsection{Miscellaneous lemmas}
\label{sec-lemmas}

\Lemma{\label{appendlemma1}}{
Under the hypothesis of Theorem~\ref{cocycleBCG}, we have, 
for all $p=1,\ldots,n-1$,
\begin{multline*}
  \prod_{k=p}^{n-1}\left\{ W^{[k]}_{13}\, W^{[k-p+1]}_{23}\,
   \widehat{J}{}^{[p-1,k-p+1]}_{1(2\mid 3} \,
   \widehat{J}^{[k-p+1]}_{23} \right\}
   \Big\{\prod_{k=p}^{n-1}J^{[k]}_{12}\Big\}
      \\
  =(id \otimes \Delta)(J^{[p]}) \; 
   \prod_{k=p+1}^{n-1}\Big\{ W^{[k]}_{13}\, W^{[k-p]}_{23} \, 
        \widehat{J}{}^{[p,k-p]}_{1 (2 \mid 3} \,
        \widehat{J}^{[k-p]}_{23} \Big\} 
     \Big\{ \prod_{k=p+1}^{n-1} J^{[k]}_{12} \Big\}\, J^{[n-p]}_{23},  
\end{multline*}
with $\widehat{J}{}^{[k,m]}_{1 (2 \mid 3}$ defined as in \eqref{defJ1}.
}

\proof
Reorganising the factors in the product and using successively 
\eqref{prop1}, \eqref{prop2},
\eqref{prop3}, \eqref{propJ3}, \eqref{prop11} and \eqref{propJ2}, we can 
reexpress the left hand side as
\begin{align*}
  LHS 
%  &=W^{[p]}_{13} \left(W^{[1]}_{23}\widehat{J}{}^{[p-1,1]}_{1(2\mid 3}W^{[1]\;-1}_{23} \right)  J^{[1]}_{23}\times\\
%  &\times \prod_{k=p+1}^{n-1}\left( W^{[k]}_{13}\; \left( W^{[k-p+1]}_{23}\widehat{J}{}^{[p-1,k-p+1]}_{1(2\mid 3}W^{[k-p+1]\;-1}_{23} \right) J^{[k-p+1]}_{23} \right) \;\left( \prod_{k=p}^{n-1}J^{[k]}_{12} \right)\\
 &=W^{[p]}_{13} \left(W^{[1]}_{23}\,\widehat{J}{}^{[p-1,1]}_{1(2\mid 3}\,
    (W^{[1]}_{23})^{-1} \right)\\
 &\quad\times
    \prod_{k=p+1}^{n-1}\Big\{ J^{[k-p]}_{23}\, W^{[k]}_{13}\, 
    \Big( W^{[k-p+1]}_{23}\,\widehat{J}{}^{[p-1,k-p+1]}_{1(2\mid 3}\,
           (W^{[k-p+1]}_{23})^{-1} \Big)  \Big\}\, J^{[n-p]}_{23}\,
    \Big\{ \prod_{k=p}^{n-1}J^{[k]}_{12} \Big\},
       \\
 &=W^{[p]}_{13} \left(W^{[1]}_{23}\,\widehat{J}{}^{[p-1,1]}_{1(2\mid 3}\,
   (W^{[1]}_{23})^{-1} \right)\,
  \prod_{k=p+1}^{n-1}\Big\{ J^{[k-p]}_{23}\, W^{[k]}_{13}\, \\ 
 &\quad\times
  \Big( W^{[k-p+1]}_{23}\,\widehat{J}{}^{[p-1,k-p+1]}_{1(2\mid 3}\,
  (W^{[k-p+1]}_{23})^{-1} \Big)  \Big\}\, J^{[p]}_{12} \,
  \Big\{ \prod_{k=p+1}^{n-1}J^{[k]}_{12} \Big\}\, J^{[n-p]}_{23},
       \\
 &=W^{[p]}_{13}\,W^{[p]}_{12} \left( 
   (W^{[p]}_{12})^{-1}\,W^{[1]}_{23}\,\widehat{J}{}^{[p-1,1]}_{1(2\mid 3}\, 
   (W^{[1]}_{23})^{-1} W^{[p]}_{12}  \right)\,
   \prod_{k=p+1}^{n-1} \Big\{(W^{[p]}_{12})^{-1} J^{[k-p]}_{23}\, 
     W^{[p]}_{12}\,W^{[k]}_{13} \\
 &\quad\times 
   \left( W^{[k-p+1]}_{23}\,(W^{[p]}_{12})^{-1} 
   \widehat{J}{}^{[p-1,k-p+1]}_{1(2\mid 3}\, W^{[p]}_{12}\,
   (W^{[k-p+1]}_{23})^{-1} \right)  \Big\}\, 
   \widehat{J}{}^{[p]}_{12} \,
   \Big\{ \prod_{k=p+1}^{n-1}J^{[k]}_{12} \Big\}\, J^{[n-p]}_{23},
       \\
 &=W^{[p]}_{13}\, W^{[p]}_{12} \left((W^{[p]}_{12})^{-1} W^{[1]}_{23}\,
    \widehat{J}{}^{[p-1,1]}_{1(2\mid 3}\, (W^{[1]}_{23})^{-1} W^{[p]}_{12}
      \right)\\
 &\quad\times
    \prod_{k=p+1}^{n-1}\left\{ W^{[k]}_{13} \,J^{[k-p]}_{23}
    \left( (W^{[p]}_{12})^{-1} W^{[k-p+1]}_{23}\,
    \widehat{J}{}^{[p-1,k-p+1]}_{1(2\mid 3}\, (W^{[k-p+1]}_{23})^{-1}
    W^{[p]}_{12} \right)  \right\}\,\\
 &\quad\times \widehat{J}{}^{[p]}_{12}\,
    \Big\{ \prod_{k=p+1}^{n-1}J^{[k]}_{12} \Big\}\, J^{[n-p]}_{23},
    \\
 &=W^{[p]}_{13}\,W^{[p]}_{12}\, \widehat{J}{}^{[p,0]}_{1 \mid 2) 3} \,
    \prod_{k=p+1}^{n-1}\left\{ W^{[k]}_{13}\,W^{[ k-p]}_{23}\,
    \widehat{J}{}^{[k-p]}_{23}\, \widehat{J}{}^{[p,k-p]}_{1 \mid 2) 3} \right\}
    \, \widehat{J}{}^{[p]}_{12}\,
    \Big\{ \prod_{k=p+1}^{n-1} J^{[k]}_{12} \Big\}    \,J^{[n-p]}_{23},
    \\
 &=W^{[p]}_{13}\, W^{[p]}_{12}\, \widehat{J}{}^{[p,0]}_{1 \mid 2) 3} \,
            \widehat{J}{}^{[p]}_{12}
   \prod_{k=p+1}^{n-1}\left\{ W^{[k]}_{13}\, W^{[k-p]}_{23} \, 
   \widehat{J}{}^{[p,k-p]}_{1 ( 2 \mid 3}\, \widehat{J}{}^{[k-p]}_{23}\right\} 
   \Big\{ \prod_{k=p+1}^{n-1} J^{[k]}_{12} \Big\} \, J^{[n-p]}_{23},
   \\
   &=RHS,
\end{align*}
which concludes the proof.
\qed

\bigskip

%%%%%%%%%%%%%%%%%%%%%%%%%%%%%
\Lemma{\label{appendlemma2}}{
With the hypothesis of Theorem \ref{maintheorem}, we have
\begin{align}  
  U_{12}(x)=S^{[1]}_{12}\prod_{k=1}^{n}\left( 
  {\mathfrak C}^{[+k]}_1(xq^{h_2})\, 
   (S^{[k+1]}_{12})^{-1}\,\widehat{J}{}^{[k]}_{12}\,S^{[k+1]}_{12} \right) ,
\end{align}
where $U_{12}(x)=\Delta\big(\mathcal{M}^{(+)}(x)\big)\, J\, 
\mathcal{M}^{(+)}_2(x){}^{-1}$.
}

\proof
Let us first mention some useful relations which can be derived from the properties of $\tau$ and $S^{[k]}$:
\begin{align}
   &\Delta({\mathfrak C}^{[+k]}(x))=S^{[1]}_{21}\,K_{12}\;
     {\mathfrak C}^{[+k]}_1(x)\;(S^{[1]}_{21})^{-1} K^{-1}_{12}\,
     S^{[1]}_{12}\;{\mathfrak C}^{[+k]}_2(x)\;(S^{[1]}_{12})^{-1},
        \label{Delta-plus-k}\\
   &(S^{[1]}_{12})^{-1} K_{12}\, S^{[1]}_{21}
    =q^{2\zeta^{(n)}\otimes \zeta^{(n)}}
     (\tau \otimes id)\big( (S^{[1]}_{12})^{-1} K^{-1}_{12} S^{[1]}_{21}\big),
         \label{tau-cart}\\
   &(\tau^p\otimes id)(R)=q^{\zeta^{(n-p+1)}\otimes \zeta^{(n)}}\,
    (S^{[p]})^{2} \,(S^{[p+1]})^{-1} (S^{[p-1]})^{-1}
     \widehat{J}{}^{[p]},\quad\forall p\ge 1,
        \label{tauR}\\	
%   &[ (S^{[k+1]}_{12})^{-1}\widehat{J}{}^{[k]}_{12}\,S^{[k+1]}_{12},
%          {\mathfrak C}^{[+(n-j)]}_2(x)]=0,\quad
%         \forall j<k,
%	 \label{cpluskcommut}\\
   &{\mathfrak C}^{[+k]}_1(xq^{h_2})
    =(\tau^{k-1} \otimes id)
      \big(K_{12}\,(S^{[1]}_{12})^{-1} S^{[1]}_{21}\,
         {\mathfrak C}^{[+]}_1(x)\, K_{12}^{-1} S^{[1]}_{12}\,
	 (S^{[1]}_{21})^{-1}\big)
	     \nonumber\\
   &\phantom{{\mathfrak C}^{[+k]}_1(xq^{h_2})}
    =(S_{12}^{[k]})^{-2} (S_{12}^{[k-1]})^{2} \; 
    {\mathfrak C}^{[+k]}_1(x)\; 
     (S_{12}^{[k]})^2 \,(S_{12}^{[k-1]})^{-2},\quad\forall k\geq 2.
     \label{cpluskshift}
\end{align}
Using \eqref{Delta-plus-k}, reorganising the factors in the product and using
\eqref{tau-cart}, we can rewrite $\Delta\big(\mathcal{M}^{(+)}(x)\big)$
as
\begin{align}
  \Delta\big(\mathcal{M}^{(+)}(x)\big)
  &=\prod_{k=1}^n \Delta\big(\mathfrak{C}^{[+k]}(x)\big),\nonumber\\
  &=  K_{12}\, S_{21}^{[1]}\;\mathfrak{C}^{[+]}_1(x)\; (S^{[0]})^{-1}\,
     q^{-\zeta^{(n)}\otimes\zeta^{(n)}}
     \nonumber\\
  &\hspace{2cm}\times (\tau\otimes id)\bigg(
            \prod_{k=1}^n \Delta'\big(\mathfrak{C}^{[+k]}(x)\big)\bigg)\ 
       (S^{[0]})^{-1} S_{12}^{[1]}\, q^{\zeta^{(n)}\otimes\zeta^{(n)}}.
\label{prod1}
\end{align}
On the other hand, from \eqref{tauR}, $J$ can be expressed as
\begin{align}
  J &=  \prod_{p=1}^n \big\{ S^{[p]}\,(S^{[p+1]})^{-1}\widehat{J}{}^{[p]}\big\},
       \nonumber\\
    &= q^{-\zeta^{(n)}\otimes\zeta^{(n)}} S^{[0]}\, (S^{[1]})^{-1}\;
    \prod_{p=1}^{n-1} \Big\{ (\tau^p\otimes id)(R)\;
    q^{-\zeta^{(n-p)}\otimes\zeta^{(n)}}\, S^{[p]}\, (S^{[p+1]})^{-1}\Big\}
 %      \nonumber\\  
 %  &\hspace{9cm}\times    q^{\zeta^{(1)}\otimes\zeta^{(n)}}(S^{[n-1]})^{-1}
 .
\label{prod2}
\end{align}
Therefore,
\begin{multline}
  \Delta\big(\mathcal{M}^{(+)}(x)\big)\, J 
  =  %K_{12}\, S_{21}^{[1]}\;\mathfrak{C}^{[+]}_1(x)\; (S^{[0]})^{-1}\,
     S^{[1]}\,\mathfrak{C}^{[+]}_1(x q^{h_2})\, S^{[0]}\,(S^{[1]})^{-2}\,
     q^{-\zeta^{(n)}\otimes\zeta^{(n)}}
     \ (\tau\otimes id)\bigg(
            \prod_{k=1}^n \Delta'\big(\mathfrak{C}^{[+k]}(x)\big)\bigg)\\
          \times \prod_{p=1}^{n-1} \Big\{ (\tau^p\otimes id)(R) \      
    q^{-\zeta^{(n-p)}\otimes\zeta^{(n)}}\, S^{[p]}\, (S^{[p+1]})^{-1}\Big\}
   %    \   q^{\zeta^{(1)}\otimes\zeta^{(n)}}(S^{[n-1]})^{-1}
   ,  
\label{prod3}
\end{multline}
where we have also used \eqref{axiomABRR1} to reexpress the first factor.
To reorganise factors in the product, we use recursively the following relation, derived from the quasitriangularity
property \eqref{quasitriangularity2}:
\begin{align*}
  (\tau^p\otimes id)\bigg(
   \prod_{k=1}^n \Delta'\big(\mathfrak{C}^{[+k]}(x)\big)\bigg)\
  (\tau^p\otimes id)(R)
  &= (\tau^p\otimes id)(R)\ (\tau^p\otimes id)\bigg(
   \prod_{k=1}^n \Delta\big(\mathfrak{C}^{[+k]}(x)\big)\bigg),
\end{align*}
which, using \eqref{prod1} and \eqref{tauR}, can be rewritten as,
\begin{multline}
  (\tau^p\otimes id)\bigg(
   \prod_{k=1}^n \Delta'\big(\mathfrak{C}^{[+k]}(x)\big)\bigg)\
  \Big\{(\tau^p\otimes id)(R)\ q^{\zeta^{(n-p)}\otimes\zeta^{(n)}}\, 
  S^{[p]}\, (S^{[p+1]})^{-1}\Big\}\\
  = \Big\{q^{\zeta^{(n-p+1)}\otimes \zeta^{(n)}}\,
    (S^{[p]})^{2} \,(S^{[p+1]})^{-1} (S^{[p-1]})^{-1}
     \widehat{J}{}^{[p]}\, (S^{[p]})^{2}\, (S^{[p+1]})^{-1} \\
  \times   
     \mathfrak{C}^{[+(p+1)]}_1(x)\; 
     (S^{[p]})^{-1}\, q^{-\zeta^{(n-p)}\otimes\zeta^{(n)}}\Big\}\
   (\tau^{p+1}\otimes id)\bigg(
   \prod_{k=1}^n \Delta'\big(\mathfrak{C}^{[+k]}(x)\big)\bigg).    
\label{prod4}
\end{multline}
This gives us
\begin{align}
  \Delta\big(\mathcal{M}^{(+)}(x)\big)\, J 
  &=  S^{[1]}\,\mathfrak{C}^{[+]}_1(x\, q^{h_2})\,(S^{[2]})^{-1}
   \widehat{J}{}^{[1]}\nonumber\\
   &\quad\times
   \prod_{p=2}^n\Big\{(S^{[p-1]})^{2}\, (S^{[p]})^{-1}
   \mathfrak{C}^{[+p]}_1(x)\;(S^{[p]})^{2}\,(S^{[p-1]})^{-2} (S^{[p+1]})^{-1}
   \widehat{J}{}^{[p]}\Big\}\nonumber\\
   &\quad\times %q^{-\zeta^{(1)}\otimes\zeta^{(n)}}S^{[n-1]}\
   (\tau^{n}\otimes id)\bigg(
   \prod_{k=1}^n \Delta'\big(\mathfrak{C}^{[+k]}(x)\big)\bigg)\ 
   %q^{\zeta^{(1)}\otimes\zeta^{(n)}}(S^{[n-1]})^{-1}
   ,
   \nonumber\\
  &=\prod_{p=1}^n\Big\{S^{[p]}\,\mathfrak{C}^{[+p]}_1(x\, q^{h_2})\,
   (S^{[p+1]})^{-1}\widehat{J}{}^{[p]}\Big\}\
    \prod_{k=1}^n \mathfrak{C}^{[+k]}_2(x),
  \label{prod5}
\end{align}
where we have used \eqref{cpluskshift} and the nilpotency of $\tau$.
We finally multiply this expression by $\mathcal{M}_2^{(+)}(x)^{-1}$, 
and this concludes the proof.
\qed

\bigskip

%%%%%%%%%%%%%%%%%%%%%%%%%%%%
\Lemma{\label{appendlemma3}}{
Under the hypothesis of Theorem~\ref{theo-Mexpr}, we have
\begin{align*}
   (\stackrel{f}{\pi}\otimes \stackrel{f}{\pi})
   ({\cal W}_{12}-\widetilde{\cal W}_{12})(x)=0,
\end{align*}
with
\begin{align*}
   &{\cal W}_{12}(x)= {\mathfrak C}^{[+]}_1(xq^{h_2})\,
   \{B_{2}(x)\,(S^{[2]}_{12})^{-1}\widehat{J}^{[1]}_{12}\,
     S^{[2]}_{12}\,B_{2}(x)^{-1}\}\,
   {\mathfrak C}_2^{[-]}(x)^{-1}, \\
   &\widetilde{\cal W}_{12}(x)=
    {\mathfrak C}^{[-]}_2(x)^{-1}
   \{(S^{[1]}_{12})^{-1}\widehat{R}_{12}\,S^{[1]}_{12}\}\,
   {\mathfrak C}_1^{[+]}(xq^{h_2}).
\end{align*}
}

\proof
A direct computation leads to
\begin{align*}
 &(\stackrel{f}{\pi}\otimes \stackrel{f}{\pi})
  ((S^{[1]}_{12})^{-1} \widehat{R}_{12}\,S^{[1]}_{12})
   = 1 \otimes 1+ (1-q^{-2})\sum_{i<j}q^{\frac{2(i-j)}{(n+1)}}\ 
                         E_{i,j}\otimes E_{j,i},\\
 &(\stackrel{f}{\pi}\otimes \stackrel{f}{\pi})
  (B_{2}(x)\,(S^{[2]}_{12})^{-1} \widehat{J}^{[1]}_{12}\,S^{[2]}_{12}\,
  B_{2}(x)^{-1})
  = 1 \otimes 1\\
  &\hspace{6.5cm}+ (1-q^{-2})\sum_{i<j}q^{\frac{3(j-i)}{(n+1)}}\,
    \frac{\nu_{j+1}}{\nu_{i+1}}\ E_{i,j}\otimes E_{j+1,i+1},\\
 &(\stackrel{f}{\pi}\otimes \stackrel{f}{\pi})
  ({\mathfrak C}^{[-]}_2(x)^{-1}) 
  =1 \otimes1+\sum_{k=1}^{n} q^{\frac{k-1}{n+1}}\, \nu_{k+1}
          \ 1 \otimes E_{k+1,k},\\
 &(\stackrel{f}{\pi}\otimes \stackrel{f}{\pi})
  ({\mathfrak C}^{[+]}_1(xq^{h_2}))=1 \otimes1+\sum_{i<j}\sum_k (-1)^{i-j}
   q^{-\frac{(j-i)(i+j-7)}{2(n+1)}-2 \delta_{i<k \leq j}}\,
   \prod_{l=i+1}^j\nu^{-1}_{l}\  E_{i,j}\otimes E_{kk}.
\end{align*}
Using these intermediary results, it is then
straightforward to check the announced result.
\qed

%%%%%%%%%%%%%%%%%%%%%%%%%%%%%%%%%%%%%%%%%%%%%%%%%%%%%%%%%%%%%%%%%%%%%%%%%%%%%%%
\subsection{A shortcut construction of $M(x)$ in the $U_q(sl(2))$ case  }
\label{sec-constrM}
A simpler alternative construction of $M(x)$ solution of the SQDBP can be done in the rank one case. This suggests a deeper relation between Primitive Representations of reflection algebras and  Dynamical coboundaries  but we unfortunately have not been able to generalize the present method for $U_q(sl(n+1))$ with $n\geq 2.$ The proof is however so simple that we have not resisted to include it  here.

\Lemma{}{
Let us  define a map $G:\mathbb{C}\rightarrow U_q(sl(2))^{\otimes 2}$ 
by
$G_{12}(x)=B_1(x)^{\onehalf}R_{12}(x)K_{12}^{-1}B_1(x)^{-\onehalf},$  we have the following 
dynamical quasitriangularity property 
\begin{equation}
(\Delta\otimes id)(G(x))=F_{12}(x)G_{13}(xq^{h_2})G_{23}(x)F_{12}(xq^{h_3})^{-1}.
\end{equation}
Moreover $G(x)$ is an element of ${\cal L}_q(sl(2))^{\otimes 2}.$ As a result we have also $F_{12}(xq^{h_3}) \in {\cal L}_q(sl(2))^{\otimes 3},$ and $(id \otimes id \otimes {\cal E})(F_{12}(xq^{h_3}))=1 \otimes 1$ for any primitive representation ${\cal E}$ of ${\cal L}_q(sl(2)).$
}
\Proof{It is straightforward, using the quasitriangularity property (\ref{quasitriangularity}) for $R$, the dynamical cocycle equation (\ref{eq:s-cocycle}) and the zero weight property of $F(x),$ to obtain that $G(x)$ satisfies the dynamical quasitriangularity equation.
{}From the explicit expression of $R(x)$ and from the isomorphism between
 $U_q(sl(2))$ and $\hat{{\cal L}}_q(sl(2))$ we have the property  
$G(x)\in {\cal L}_q(sl(2))^{\otimes 2}.$ (This is the step we are not able to generalize in the higher rank case). 
The other properties are trivial. 
}
\Proposition{}{Due to the previous lemmas, for any primitive representation
${\cal E}$ of ${\cal L}_q(sl(2)),$ it makes sense to
define the following map from ${\mathbb C}$ to $ {\cal L}_q(sl(2))$ by 
\begin{equation}
M^{({\cal E})}(x)=(id \otimes {\cal E})(G(x)).\label{M=EG}
\end{equation}
For any ${\cal E},$ $M^{({\cal E})}(x)$ verifies the QDBE. Its explicit expression is given by (\ref{MdeXsl2}) up to trivial gauge transformations.
}
\Proof{
Using previous lemmas and the fact that ${\cal E}$ is a morphism, we have
\begin{align*}
\Delta(M^{({\cal E})}(x))&=(\Delta \otimes {\cal E})(G(x))=(id \otimes id \otimes {\cal E})
\left( F_{12}(x)G_{13}(xq^{h_2})G_{23}(x)F_{12}(xq^{h_3})^{-1} \right)\\
&=F_{12}(x)M^{({\cal E})}_{1}(xq^{h_2})M^{({\cal E})}_{2}(x)
\end{align*}
which concludes the proof of the dynamical coboundary equation. The explicit expression (\ref{MdeXsl2})
is recovered from (\ref{M=EG}) by a trivial computation. 
}

\newpage 

%%%%%%%%%%%%%%%%%%%%%%%%%%%%%%%%%%%%%%%%%%%%%%%%%%%%%%%%%%%%%%%%%%%%%%%%%%%%%%%
%%%%%%%%%%%%%%%%%%%%%%%%%%%%%%%%%%%%%%%%%%%%%%%%%%%%%%%%%%%%%%%%%%%%%%%%%%%%%%%

\newpage
\end{document}